\documentclass[a4paper,10pt]{article}

\usepackage{amsmath}
\usepackage{amssymb}
\usepackage{amsfonts}
\usepackage{amsthm}
\usepackage{bbm}
\usepackage[mathscr]{eucal}
\usepackage{a4wide}
\usepackage{authblk}
\usepackage{float}

\usepackage{hyperref}
\hypersetup{
    colorlinks = true,
    urlcolor = {magenta},
    citecolor = {red},
    linkcolor = {blue}
}
\usepackage{doi}
\usepackage[square,numbers]{natbib}

\usepackage{xcolor}
\usepackage{todonotes}

\usepackage{tikz}
\usetikzlibrary{calc,fadings,decorations.pathreplacing}
\usetikzlibrary{shapes,snakes}

\newcommand\EncircleNode[2][0.3]{
\draw ($ (#2) + (#1,0) $) arc (0:360:#1);
}

\tikzset{%
  >=latex, 
  inner sep=0pt,%
  outer sep=2pt,%
  rootnode/.style={inner sep=0pt,outer sep=0pt,minimum size=15pt,fill=black,star,star points=7},
  root2/.style={inner sep=0pt,outer sep=0pt,minimum size=8pt,fill=black,circle},
  neighborhood2/.style={inner sep=0pt,outer sep=0pt,minimum size=7pt,draw=black,diamond},
  notneighborhood2/.style={inner sep=0pt,outer sep=0pt,minimum size=5pt,fill=black,circle},
  insidecircle/.style={inner sep=0pt,outer sep=0pt,minimum size=5pt,fill=black,circle}%
}

\newcommand{\RNum}[1]{\uppercase\expandafter{\romannumeral #1\relax}}

\newcommand{\pr}{\mathbb{P}}								

\newcommand{\Prob}[1]{\pr\left(#1\right)}					
																
\newcommand{\CProb}[2]{\pr\left(#1 \bigg| #2\right)}	

\newcommand{\e}{\mathbb{E}}								

\newcommand{\Exp}[1]{\e\left[#1\right]}					
															
\newcommand{\CExp}[2]{\e\left[\left.#1\right | #2\right]}	

\newcommand{\plim}{\ensuremath{\stackrel{\pr}{\rightarrow}}}	
\newcommand{\dlim}{\ensuremath{\stackrel{d}{\rightarrow}}}		

\newcommand{\1}{\mathbbm{1}}								
\newcommand{\ind}[1]{\1_{\left\{#1\right\}}}	
\newcommand{\indE}[1]{\1_{#1}}					

\newcommand\numberthis{\addtocounter{equation}{1}\tag{\theequation}}

\allowdisplaybreaks

\newtheorem{theorem}{Theorem}[section]
\newtheorem{lemma}[theorem]{Lemma}
\newtheorem{proposition}[theorem]{Proposition}
\newtheorem{corollary}[theorem]{Corollary}
\newtheorem{conjecture}[theorem]{Conjecture}
\newtheorem{definition}[theorem]{Definition}
\newtheorem{assumption}[theorem]{Assumption}
\newtheorem{remark}[theorem]{Remark}

\newcommand{\eqan}[1]{\begin{align}#1\end{align}}

\numberwithin{equation}{section}

\title{Local limits of spatial inhomogeneous random graphs}
\author[1]{Remco van der Hofstad, Pim van der Hoorn, Neeladri Maitra}

\affil[1]{Department of Mathematics and Computer Science, Eindhoven University of Technology}

\begin{document}

\maketitle

\begin{abstract}
    Consider a set of $n$ vertices, where each vertex has a location in $\mathbb{R}^d$ that is sampled uniformly from the unit cube in $\mathbb{R}^d$, and a weight associated to it. Construct a random graph by placing edges independently for each vertex pair with a probability that is a function of the distance between the locations, and the vertex weights.
    
    Under appropriate integrability assumptions on the edge probabilities that imply sparseness of the model, after appropriately blowing up the locations, we prove that the local limit of this random graph sequence is the (countably) infinite random graph on $\mathbb{R}^d$ with vertex locations given by a homogeneous Poisson point process, having weights which are i.i.d.\ copies of limiting vertex weights.  Our setup covers many sparse geometric random graph models from the literature, including Geometric Inhomogeneous Random Graphs (GIRGs), Hyperbolic Random Graphs, Continuum Scale-Free Percolation and Weight-dependent Random Connection Models. 
    
    We prove that the limiting degree distribution is mixed Poisson, and the typical degree sequence is uniformly integrable, and obtain convergence results on various measures of clustering in our graphs as a consequence of local convergence. Finally, as a by-product of our argument, we prove a doubly logarithmic lower bound on typical distances in this general setting.
\end{abstract}

\vspace{5 pt}



\section{Introduction and main results }\label{sec:Introduction}

Random graphs with underlying geometry are becoming the model of choice when it comes to modeling and understanding real-world networks \cite{SFP, KPKVB_HRG, KB_Avg_dist_16, WDRCM_Rec_trans} (see also \cite[Section 9.5]{RGCN_2}). The presence of an ambient geometric space enables one to model the natural tendency of forming connections between entities that are \emph{close} to each other, where \emph{closeness} is measured in terms of the underlying geometry. The power of these models is the inherent diversity of the geometric component. It can encode actual physical distance, such as two servers in adjacent cities, as well as a more abstract form of similarity, e.g., users with similar interests or hobbies. Spatial networks have been employed to study social networks \cite{Spatial_social_nets}. Empirical research shows spatial positions of individuals often play a role in formation of friendships among them  \cite{spatial_friendship_nets,Nbr_nets_B_W,wellman1988networks,diff_strokes_friendhsip}. 

The prototypical random graph model in this setting is the random geometric graph, first introduced by Gilbert~\cite{Gilbert_61} and later popularised by Penrose~\cite{Penrose_RGGs}. Here the graph is formed by placing edges between pairs of points of some Poisson point process on a Euclidean space, if and only if their metric distance is smaller than some fixed parameter. In general, one can consider models in which the connection probabilities are a decreasing function of the distance between pairs of vertices. Since nearby vertices are more likely to be connected in these spatial models, they naturally exhibit {\em clustering}, which captures the tendency of the existence of a connection between two entities having a common neighbor, a feature that is often observed in real-world networks.

Other than clustering, experimental studies suggest that most real-world networks are {\em sparse}, i.e., the number of connections is often of the same order as the number of individuals, and that they are small worlds, i.e., for most pairs of nodes, it takes only a small number of connection steps to reach one from the other. Finally, many networks have highly inhomogeneous degree distributions, i.e., many vertices have only a few connections, while a small proportion of vertices have a lot of connections. The most notable spatial random graph models capable of capturing all these features are scale-free percolation~\cite{SFP}, geometric inhomogeneous random graphs~\cite{KB_GIRGs_19} and hyperbolic random graphs~\cite{KPKVB_HRG}.

When studying network models, one is often interested in the limits of specific network measures, i.e., degree distribution, average path length or clustering coefficients. There is plenty of literature analyzing limits of such measures for a wide variety of models. Here, however, one first must prove the convergence. Suppose that instead we would have a limit object for our graph models, which implies convergence of the network measures of interest. Then we no longer have to worry about convergence and can instead focus on actually analyzing the limit measures. For sparse graphs such a framework exists under the name of local weak convergence~\cite{BSc01,AS04}. At a high level, if a model has a local weak limit, then any local property converges to an associated measure on the limit graph. Therefore, if we know the local graph limit, we immediately obtain the limits of a wide variety of local network measures \cite[Chapter 2]{RGCN_2}. Furthermore, sometimes conclusions about sufficiently global properties of the graph can also be obtained from its local behaviour \cite{Giant_local}, \cite[Chapter 2]{RGCN_2}.


In this paper, we study local convergence of a general class of spatial random graphs, which cover as particular cases many well-known models that are sparse, have inhomogeneous degree distributions, short path lengths, and exhibit non-vanishing clustering. We establish that the typical local behaviour in these ensembles matches with the expected local behaviour of the natural infinite version of the model. Our results immediately imply convergence of several local network measures, such as the degree distribution of a random vertex and the local clustering coefficient.

\paragraph{\bf Organization of the paper.}
In the rest of this section, we review basics of local convergence of graphs, introduce our model, and state our main results. Next, in Section~\ref{sec:discussion}, we present examples of models that are covered, as well as results on degrees and clustering in our graphs, and a small discussion on typical distances in our models. All the proofs are deferred to Section~\ref{sec:proofs}.

\subsection{The space of rooted graphs and local convergence}

The notion of local convergence of graphs was first introduced by Benjamini and Schramm~\cite{BSc01} and independently by Aldous and Steele~\cite{AS04}. Intuitively, this notion studies the asymptotic local properties of a graph sequence, as observed from a typical vertex.  

A (possibly infinite) graph $G=(V(G),E(G))$ is called {\em locally finite} when every vertex has finite degree. A rooted locally finite graph is a tuple $(G,o)$, where $G=(V(G),E(G))$ is a locally finite graph, with a designated vertex $o \in V(G)$ called \emph{the root}.

\begin{definition}[Rooted isomorphism]\label{defn:rooted_iso}
For rooted locally finite graphs $(G_1,o_1)$ and $(G_2,o_2)$ where $G_1=(V(G_1),E(G_1))$ and $G_2=(V(G_2),E(G_2))$, we say that $(G_1,o_1)$ is {\em rooted isomorphic} to $(G_2,o_2)$ when there exists a graph isomorphism between $G_1$ and $G_2$, mapping $o_1$ to $o_2$, i.e., when there exists a bijective function $\phi\colon V(G_1) \to V(G_2)$ such that 
\[\{u,v\}\in E(G_1) \iff \{\phi(u),\phi(v)\}\in E(G_2),\]
and satisfying $\phi(o_1)=o_2$.
\end{definition}
We use the notation
\[
(G_1,o_1)\cong (G_2,o_2)
\]
to denote that there is a rooted isomorphism between $(G_1,o_1)$ and $(G_2,o_2)$.

\medskip

Let $\mathcal{G}_{\star}$ be the space of all rooted isomorphic equivalence classes of locally finite rooted graphs, i.e., $\mathcal{G}_{\star}$ consists of the equivalence classes $[(G,o)]$ of rooted locally finite graphs, where two rooted locally finite graphs $(G_1,o_1)$ and $(G_2,o_2)$ belong to the same class if they are rooted isomorphic to each other. We often omit the equivalence class notation, and just write $(G,o)$ for an element of $\mathcal{G}_{\star}$. 

Fix a graph $G=(V(G),E(G))$. We denote the graph distance in $G$ by $d_G$, i.e., for any two vertices $u,v \in V(G)$, $d_G(u,v)$ equals the length of the shortest path in $G$ from $u$ to $v$, and we adopt the convention that $d_G(u,v)=\infty$ whenever $u$ and $v$ are not connected by a sequence of edges. 

\begin{definition}[Neighbourhood of the root]\label{defn:nbhds}
For any $R\geq 1$ and $(G,o) \in \mathcal{G}_{\star}$ where $G=(V(G),E(G))$, we call the element $\left(B^G_o(R),o\right)$ of $\mathcal{G}_{\star}$, the \emph{$R$-neighbourhood of $o$ in $G$}, where $B^G_o(R)$ is the subgraph of $G$ induced by
\[\{v \in V(G)\colon d_G(o,v) {\leq} R\}.\]
We sometimes abbreviate $(B^G_o(R),o)$ as $B^G_o(R)$.
\end{definition}

The space $\mathcal{G}_{\star}$ is usually endowed with the \emph{local topology}, which is the smallest topology that makes the functions of the form $(G,o) \mapsto {\ind{B^G_o(K) \cong (G',o')}}$ for $K \geq 1$ and $(G',o') \in \mathcal{G}_{\star}$ continuous. This topology is metrizable with an appropriate metric $d_{\star}$ (see \cite[Definition 2.4]{RGCN_2}), and $(\mathcal{G}_{\star},d_{\star})$ is a Polish space, which enables one to do probability on it. We omit further details on this.

\medskip

For the next two definitions, we assume that $(G_n)_{n \geq 1}$ with $G_n=(V(G_n),E(G_n))$ is a sequence of (possibly random) graphs that are (almost surely) finite, i.e., $|V(G_n)| \stackrel{a.s.}{<} \infty$ for all $n \geq 1$, and, conditionally on $G_n$, $U_n$ is a uniformly chosen random vertex of $G_n$. Note then that $(G_n,U_n)$ is a random variable taking values in $\mathcal{G}_{\star}$. Local weak convergence of such random variables is defined as follows:

\begin{definition}[Local weak convergence]\label{defn:LWC}
The sequence of graphs $(G_n)_{n \geq 1}$ is said to {\em converge locally weakly} to the random element $(G,o) \in \mathcal{G}_{\star}$ having law $\mu_{\star}$, as $n \to \infty$, when, for every $r>0$, and for every $(G_{\star},o_{\star}) \in \mathcal{G}_{\star}$,
\[\Prob{B^{G_n}_{U_n}(r) \cong (G_{\star},o_{\star})} \to \mu_{\star}\left({B^{G}_o(r)\cong (G_{\star},o_{\star})}\right),\]
as $n \to \infty$. 
\end{definition}

Definition \ref{defn:LWC} is in fact equivalent to saying that the sequence $((G_n,U_n))_{n \geq 1}$ of random variables taking values in $\mathcal{G}_{\star}$ converges in distribution to the random variable $(G,o)$ taking value in $\mathcal{G}_{\star}$ as $n \to \infty$, under the topology induced by the metric $d_{\star}$ on $\mathcal{G}_{\star}$ (for a proof of this fact, see \cite[Definition 2.10, and Theorem 2.13]{RGCN_2}).

This concept of \emph{local convergence} can be adapted to the setting of convergence in probability. For a sequence of random variables $Z, (Z_n)_{n \geq 1}$, we write $Z_n \plim Z$ to indicate that $Z_n$ converges in probability to $Z$, as $n \to \infty$.

\begin{definition}[Convergence locally in probability]\label{defn:LPC}
The sequence of graphs $(G_n)_{n \geq 1}$ is said to {\em converge locally in probability} to the (possibly) random element $(G,o) \in \mathcal{G}_{\star}$ having (possibly random) law $\mu$ if for any $r>0$ and any $(G_{\star},o_{\star}) \in \mathcal{G}_{\star}$,
\[
    \CProb{B^{G_n}_{U_n}(r) \cong (G_{\star},o_{\star})}{G_n}=\frac{1}{|V(G_n)|}\sum_{i \in V(G_n)}\ind{B^{G_n}_i(r) \cong (G_{\star},o_{\star})} \plim \mu\left(B^G_o(r) \cong (G_{\star},o_{\star}) \right),
\]
as $n \to \infty$.
\end{definition}
For more on local convergence, we refer the reader to \cite[Chapter 2]{RGCN_2}.

\subsection{Model description and assumptions}
\label{sec-SIRG}

Let $S$ be either $[n]=\{1,\dots,n\}$ or $\mathbb{N}\cup \{0\}$.

Consider a sequence $\mathbf{X}=(X_i)_{i \in S}$ of (possibly random) points in $\mathbb{R}^d$, a sequence of (possibly random) reals   $\mathbf{W}=(W_i)_{i \in S}$, and a function $\kappa: \mathbb{R}_+ \times \mathbb{R} \times \mathbb{R} \to [0,1]$, which is symmetric in its second and third arguments. 

Conditionally on $\mathbf{X}$ and $\mathbf{W}$, form the (undirected) random graph $G=(V(G), E(G))$ with vertex set $V(G)=S$ and each possible edge $\{i,j\} \in E(G) \subset S \times S$ is included independently with probability 
\begin{equation}\label{eq:defn_edge_probabilities_SIRG}
    \kappa(\|X_i-X_j\|,W_i,W_j).
\end{equation}
Note that the requirement of symmetry is necessary since the probability of including the edge $\{i,j\}$ has to be the same as the probability of including the edge $\{j,i\}$. For each vertex $i$, we think of $X_i$ as the {\em spatial location} of the vertex, and $W_i$ as the {\em weight} associated to it. 

We call such a graph $G$ a \textbf{spatial inhomogeneous random graph}, or \textbf{SIRG} for short, and we use the notation
\begin{equation}\label{notn_G(X,W,kappa)}
    G(\mathbf{X},\mathbf{W},\kappa)
\end{equation}
to denote the SIRG corresponding to the location sequence $\mathbf{X}$, weight sequence $\mathbf{W}$, and connection function $\kappa$.

\begin{remark}[Non-spatial inhomogeneous random graphs]
\label{rem:model_difference_BRJ}
In the setting of (non-spatial) inhomogeneous random graph as proposed by Bollob\'as, Janson and Riordan in \cite{IRGs_BRJ}, the vertex locations of the SIRGs can be thought of as \emph{types} associated to each vertex. However, the type space in \cite{IRGs_BRJ} was taken to be an abstract space without any structure. In our model, these are locations in a Euclidean space, with connections depending upon distances between vertex locations, incorporating the underlying metric structure. In particular, in our models, the edge probabilities are large for closeby vertices, while in \cite{IRGs_BRJ}, edge probabilities are typically of the order $1/n$. This allows for non-trivial clustering to be present in our models, unlike in \cite{IRGs_BRJ}, where the local limits are multitype branching processes, see \cite[Chapter 3]{RGCN_2}.  
\end{remark}

Consider a sequence $G_n=(V(G_n)=[n],E(G_n))=G(\mathbf{X}^{(n)}, \mathbf{W}^{(n)}, \kappa_n)$ of SIRGs of size $n$. Then $\mathbf{X}^{(n)}=(X^{(n)}_i)_{i\in[n]}$ and $\mathbf{W}^{(n)}=(W^{(n)}_i)_{i\in[n]}$ are the location and weight sequences, respectively, while $\kappa_n$ is the sequence of connection functions. If $\mathbf{X}^{(n)}$, $\mathbf{W}^{(n)}$ and $\kappa_n$ have limits $\mathbf{X}$, $\mathbf{W}$ and $\kappa$, respectively, as $n \to \infty$ and in some appropriate sense, it is natural to expect the sequence of random graphs $G(\mathbf{X}^{(n)},\mathbf{W}^{(n)},\kappa_n)$ to have the graph $G(\mathbf{X},\mathbf{W},\kappa)$ as its limit. Our main contribution formalizes this intuition using local convergence, under appropriate convergence assumptions on $\mathbf{X}^{(n)}$, $\mathbf{W}^{(n)}$ and $\kappa_n$, which we discuss next.

\begin{assumption}[Law of vertex locations]
\label{sssec:assumption_location}
Define the box 
\begin{equation}\label{eq:defn_box_I}
        I := 
        \left[ - \tfrac{1}{2},\tfrac{1}{2} \right]^d.
\end{equation} 
For each $n \in \mathbb{N}$, the collection $(X^{(n)}_i)_{i\in[n]}$ is a collection of independent and uniformly distributed random variables on the box $I$.
\end{assumption}


\begin{assumption}[Weight distributions]\label{sssec:assumption_weights}
Let $\mathbf{W}^{(n)}=(W^{(n)}_i)_{i\in[n]}$ be the sequence of weights associated to the vertices of $G_n$, 
and assume that $\mathbf{W}^{(n)}$ is independent of $\mathbf{X}^{(n)}$. Further, we assume that there exists a random variable $W$, with distribution function $\mathrm{F}_W(x):=\Prob{W \leq x}$, such that for every continuity point $x$ of $\mathrm{F}_W(x)$,
 
 \begin{equation}\label{eq:conv_in_law_Wn_to_W}
    \frac{1}{n}\sum_{i=1}^n \ind{W^{(n)}_i\leq x} \plim \mathrm{F}_W(x),
\end{equation} 
as $n \to \infty$. That is, we assume convergence in probability of the (possibly random) empirical distribution function to the (deterministic) distribution function $\mathrm{F}_W$ of some limiting weight random variable $W$.
\end{assumption}

The above assumption on the weights in our model is natural. Note that if we let $U_n$ be uniform on $[n]=\{1,\dots,n\}$, then  (\ref{eq:conv_in_law_Wn_to_W}) is equivalent to the convergence in distribution of the \emph{typical} weight $W^{(n)}_{U_n}$, conditionally on the entire weight sequence $(W^{(n)}_i)_{i \in [n]}$, to $W$. Such regularity conditions are essential to understand behaviour of the \emph{typical} vertex, which is what local convergence aims to understand. Similar regularity assumptions on the weight of a typical vertex are also required for generalized random graphs (see e.g.\ \cite[Condition 6.4]{RGCN_1}) to make sure the typical degree distribution has a limit.

In most of the spatial random graph models, one starts with an infinite weight sequence, and assigns the first $n$ entries of this sequence as weights to the vertices of $G_n$. This infinite weight sequence is usually deterministic, or an independent and identically distributed (i.i.d.) sequence of weights from some given distribution.  Note that for the latter, the above assumption is satisfied using the strong law of large numbers, and for the former case, (\ref{eq:conv_in_law_Wn_to_W}) ensures that we have convergence in distribution of the typical weight.

\medskip
We also need some convergence assumption for the sequence of connection functions $\kappa_n$. 

 
\begin{assumption}[Connection functions]\label{sssec:assumption_connections}
\hfill
\begin{itemize}
        
    \item[1.] There exists a function $\kappa: \mathbb{R}_+ \times \mathbb{R} \times \mathbb{R} \to [0,1]$ such that {for every pair of real sequences $x_n\to x$ and $y_n\to y$, and for almost every $t \in \mathbb{R}_+$,} 
    \begin{equation}\label{eq:pw_conv_strong_cty_kappa_n}
      \kappa_n(n^{1/d}t,x_n,y_n)\to \kappa(t,x,y),  
    \end{equation}
        as $n \to \infty$.
        
    
    \item[2.] Let $W^{(1)}$ and $W^{(2)}$ be two i.i.d.\ copies of the limiting weight random variable $W$. Then there exists $t_0>0$, $\alpha>0$, 
    such that, for any $t \in \mathbb{R}_+$ with $t>t_0$, 
    \begin{equation}\label{eq:power_law_tail_kappa}
        \Exp{\kappa(t,W^{(1)},W^{(2)})}\leq t^{-\alpha}, 
    \end{equation}
    \end{itemize}
\end{assumption}

Before stating our results, let us reflect on Assumption \ref{sssec:assumption_connections} a little. The term $n^{1/d}$ in the first argument of $\kappa_n$ is required, because to obtain a local limit, we need to blow up the vertex locations from uniform on $I$ (recall (\ref{eq:defn_box_I})) to uniform on
\begin{equation}\label{eq:defn_box_In}
    I_n := \left[-\frac{n^{1/d}}{2}, \frac{n^{1/d}}{2} \right]^d.
\end{equation}
This is done via the transformation $x \mapsto n^{1/d}x$, and the $n^{1/d}$ term in~\eqref{eq:pw_conv_strong_cty_kappa_n} ensures convergence of $\kappa_n$ after this transformation.

Next, the `strong' form of continuity in (\ref{eq:pw_conv_strong_cty_kappa_n}) is required so that the connection function sequence $\kappa_n$ is continuous with respect to convergence of weights as formulated in (\ref{eq:conv_in_law_Wn_to_W}). To explain this, let $j \in \mathbb{N}$, and consider a sequence $(h_n)_{n \geq 1}$ of bounded continuous functions $h_n\colon \mathbb{R}^j \to \mathbb{R}$, converging to a bounded continuous function $h\colon \mathbb{R}^j \to \mathbb{R}$ in the following `strong' sense analogous to (\ref{eq:pw_conv_strong_cty_kappa_n}): for any collection $\{(x^i_n)_{n \geq 1}\colon 1\leq i \leq j\}$ of $j$ real sequences satisfying $x^i_n \to x_i$ as $n \to \infty$ for $1 \leq i \leq j$, $h_n(x^1_n,\ldots,x^j_n) \to h(x_1,\ldots,x_j)$. Then, if for each $n \in \mathbb{N}$ we let $U_{n,1},\ldots,U_{n,j}$ be $j$ independent uniformly distributed random variables on $[n]$ and let $W^{(1)},\ldots,W^{(j)}$ be $j$ i.i.d.\ copies of the limiting weight random variable $W$, (\ref{eq:conv_in_law_Wn_to_W}) implies that
\begin{equation}\label{eq:conv_bdd_func_weghts}
  \Exp{h_n(W^{(n)}_{U_{n,1}},\ldots,W^{(n)}_{U_{n,j}})} \to \Exp{h(W^{(1)},\ldots,W^{(j)})},  
\end{equation}
as $n \to \infty$. In particular, as a consequence of (\ref{eq:pw_conv_strong_cty_kappa_n}) and (\ref{eq:conv_in_law_Wn_to_W}), one can take $h_n$ in (\ref{eq:conv_bdd_func_weghts}) to be $\kappa_n$, or products of $\kappa_n$, which will enable us to compute probabilities that paths exist in our finite graphs $G_n$. Note that (\ref{eq:conv_bdd_func_weghts}) is not true in general if we just assume that the function $h_n$ is continuous on $\mathbb{R}^j$ for each fixed $n \in \mathbb{N}$. 

The above strong sense of continuity is in line with the assumption of \emph{graphical kernels} in \cite{IRGs_BRJ} (see \cite[(2.10)]{IRGs_BRJ}). In the setting of \cite{IRGs_BRJ}, if we now consider the weights (instead of vertex locations as in Remark \ref{rem:model_difference_BRJ}) associated to each vertex as vertex types, then we merely require the connection kernels to be continuous with respect to convergence of types, which is one of the key properties of \emph{graphical kernels} in \cite{IRGs_BRJ}.

{Finally, let us discuss the technicality of having the convergence (\ref{eq:pw_conv_strong_cty_kappa_n}) hold for almost every $t \in \mathbb{R}_+$, instead of all $t$. This is a purely technical condition, to avoid pathological examples. For example, consider the sequence of connection functions $\kappa_n(t,x,y):=\ind{t<n^{1/d}xy}$. We naturally want the limit of this sequence of connection functions to be $\ind{t < xy}$. It is easily verified that for any pair of real sequences $x_n\to x$ and $y_n \to y$, (\ref{eq:pw_conv_strong_cty_kappa_n}) holds, except possibly for $t\in \{xy\}$, a set of measure zero, with limiting connection function $\kappa(t,x,y)=\ind{t < xy}$. In particular, Assumption \ref{sssec:assumption_connections} (1) holds. Typically, $t$ will be replaced by $\|X^{(n)}_i-X^{(n)}_j\|$, with $X^{(n)}_i,X^{(n)}_j$ being the respective locations of two vertices $i$ and $j$. Since under Assumption \ref{sssec:assumption_location}, the random variable $\|X^{(n)}_i-X^{(n)}_j\|$ is a continuous random variable, it will avoid sets of measure zero, and Assumption \ref{sssec:assumption_connections} (1) will continue to hold in our probabilistic statements.}


\subsection{Statement of main results}
\label{ssec:statements}

Let us now formally discuss the limit of the SIRG model. As mentioned before, we shall consider the `blown up' SIRG sequence $\mathbb{G}_n=G(\mathbf{Y}^{(n)},\mathbf{W}^{(n)},\kappa_n)$, where $\mathbf{Y}^{(n)}=(Y^{(n)}_i)_{i\in[n]}$, with $Y^{(n)}_i=n^{1/d}X^{(n)}_i$ for all $i\in[n]$, and $\mathbf{X}^{(n)}=(X^{(n)}_i)_{i\in[n]}$ satisfies Assumption \ref{sssec:assumption_location}, i.e., the locations of $\mathbb{G}_n$ are now i.i.d.\ uniform variables on $I_n$ (recall (\ref{eq:defn_box_In})).

For convenience, we also define the point process on $\mathbb{R}^d$, corresponding to the spatial locations of $\mathbb{G}_n$, as
\begin{equation}\label{eq:defn_Gamma_n}
    \Gamma_n(\cdot):= \sum_{i=1}^n \delta_{Y^{(n)}_i}(\cdot),
\end{equation}
where $\delta_x$ denotes the Dirac measure at $x \in \mathbb{R}^d$.

If $\mathbb{G}_n$ has a local limit, then it is natural to expect that the weights associated to each vertex in this limit are i.i.d.\ copies of the limiting weight random variable $W$ because of Assumption \ref{sssec:assumption_weights}, and the limiting connection function to be $\kappa$. For the locations, it is standard that the sequence $\Gamma_n$, viewed as random measures on $\mathbb{R}^d$, converges to a unit-rate homogeneous Poisson point process on $\mathbb{R}^d$ (see \cite{Kallenberg_RMTA}), whose points will serve as the locations of the vertices of the limiting graph.

To this end, let us consider a unit-rate homogeneous Poisson point process $\Gamma$, and write its atoms as $\Gamma=\{Y_i\}_{i \in \mathbb{N}}$ (such an enumeration is possible by \cite[Corollary 6.5]{Last_Penrose_LPP}). Define the point process 
\begin{equation}\label{eq:defn_Gamma_infty}
    \Gamma_{\infty}:=\Gamma \cup \{\mathbf{0}\},
\end{equation} 
(where $\mathbf{0}=(0,0,\dots,0) \in \mathbb{R}^d$) which is the Palm version of $\Gamma$, and write the sequence of its atoms as $\mathbf{Y}=(Y_i)_{i \in \mathbb{N}\cup \{0\}}$, where $Y_0 = \mathbf{0}$. Also let $\mathbf{W}=(W_i)_{i \in \mathbb{N}\cup\{\mathbf{0}\}}$ be a collection of i.i.d.\ copies of the limiting weight random variable $W$. Then we define the \emph{infinite SIRG}, whose vertex set is given by $\mathbb{N\cup \{\mathbf{0}\}}$, as $\mathbb{G}_{\infty} =G(\mathbf{Y}, \mathbf{W}, \kappa)$.

Our first result establishes that $\mathbb{G}_n$ converges locally weakly to $\mathbb{G}_\infty$:

\begin{theorem}[Convergence of SIRGs in the local weak sense]\label{thm:main_LWC_SIRGs}
Consider the sequence $(\mathbb{G}_n)_{n\geq 1}=(G(\mathbf{Y}^{(n)},\mathbf{W}^{(n)},\kappa_n))_{n \geq 1}$ of SIRGs, where $Y^{(n)}_i=n^{1/d}X^{(n)}_i$ for each $n$ and each $i\in [n]$, with $\mathbf{X}^{(n)}=(X^{(n)}_i)_{i\in[n]}$, $\mathbf{W}^{(n)}=(W^{(n)}_i)_{i\in[n]}$ and $\kappa_n$  satisfying Assumptions \ref{sssec:assumption_location}, \ref{sssec:assumption_weights} and \ref{sssec:assumption_connections}, respectively, with {the parameter $\alpha$ in (\ref{eq:power_law_tail_kappa}) satisfying} 
\[\alpha > d.\]
Then $(\mathbb{G}_n)_{n \geq 1}$ converges locally weakly to the infinite rooted {\rm SIRG} $(\mathbb{G}_{\infty},0)$, rooted at vertex $0$, where $\mathbb{G}_{\infty}=G(\mathbf{Y},\mathbf{W},\kappa)$.
\end{theorem}

\medskip

\begin{remark}[Regularly varying connection functions] Note that if $\Exp{\kappa(t,W^{(1)},W^{(2)})}$ {as a function of $t$, is} either itself, or is dominated by, a regularly varying function in $t$, with some exponent greater than $d$, then using Potter's theorem (e.g., see \cite{Potter}), Assumption \ref{sssec:assumption_connections} (2) is satisfied, with some $\alpha >d$. Hence for these kind of connection functions, Theorem \ref{thm:main_LWC_SIRGs} goes through. 
\end{remark}

Theorem \ref{thm:main_LWC_SIRGs} is equivalent to the statement that if $U_n$ is uniformly distributed on $[n]$, then the random rooted graph $(\mathbb{G}_n, U_n)$ converges in distribution to the random rooted graph $(\mathbb{G}_{\infty},0)$, in the space $(\mathcal{G}_{\star},d_{\star})$.

The condition $\alpha > d$ is required to avoid non-integrability of $z\mapsto\Exp{\kappa(\|z\|,W^{(1)},W^{(2)})}$ as a function on $\mathbb{R}^d$. Integrability of $\Exp{\kappa(\|z\|,W^{(1)},W^{(2)})}$ implies that the degrees in $\mathbb{G}_{\infty}$ have finite mean, which ensures that our random graph model is sparse. It also ensures that $\mathbb{G}_{\infty}$ is locally finite almost surely, which we need for local weak convergence to make sense.

There are two main challenges in proving Theorem \ref{thm:main_LWC_SIRGs}. The first is to formulate the probability that the rooted subgraph induced by vertices whose locations fall inside a Euclidean ball with some fixed radius centered at the location $Y^{(n)}_{U_n}$ of $U_n$, is isomorphic to a given graph, in terms of suitable functionals of the spatial locations of the vertices. The latter have nice limiting behaviour owing to vague convergence of the spatial locations to a homogeneous Poisson process. The second is a careful path-counting analysis to conclude that the expected number of paths starting at $U_n$, containing vertices with spatial locations far away from the location of $U_n$, can be made arbitrarily small.           

The result of Theorem \ref{thm:main_LWC_SIRGs} can be improved to local convergence in probability:

\begin{theorem}[Convergence of SIRGs locally in probability]\label{thm:main_LCinP_SIRGs} 
Under the assumptions of Theorem \ref{thm:main_LWC_SIRGs}, the sequence of SIRGs $(\mathbb{G}_n)_{n \geq 1}$ converges locally in probability to the infinite rooted {\rm SIRG} $(\mathbb{G}_{\infty},0)$, where $\mathbb{G}_{\infty}=G(\mathbf{Y},\mathbf{W},\kappa)$.
\end{theorem}

The improvement in Theorem \ref{thm:main_LCinP_SIRGs} is achieved via a second moment analysis on neighbourhood counts. The required independence essentially follows from the fact that the spatial locations of two uniformly chosen vertices of $\mathbb{G}_n$ are with high probability far apart from each other.

\section{Consequences and discussion}\label{sec:discussion}

In this section, we discuss implications of Theorems \ref{thm:main_LWC_SIRGs} and \ref{thm:main_LCinP_SIRGs}. We first discuss some standard examples that are covered under our setting. We then discuss how local convergence implies convergence of interesting graph functionals. We focus our attention on the degree and clustering structure of our SIRGs. Finally we provide a lower bound on typical distances in our graphs.

\subsection{Examples}\label{ssec:examples_covered}
In this section, we discuss examples of spatial random graph models that are covered under our setup.

\subsubsection{Product SIRGs}

We begin by discussing a particular type of SIRGs, having a product form in their connections:
\begin{definition}
[Product SIRGs]
For a SIRG $G(\mathbf{X},\mathbf{W},\kappa)$, if the connection function $\kappa$ has the product form
\begin{equation}\label{eq:PSIRG_connection_func}
   \kappa(t,x,y)= 1 \wedge f(t)g(x,y), 
\end{equation}
for some non-negative functions $f\colon \mathbb{R}_+ \to \mathbb{R}_+$, and $g\colon\mathbb{R}^2 \to \mathbb{R}_+$, where $g$ is symmetric, then we call such a SIRG a {\em product SIRG}, or PSIRG for short, and we call $\kappa$ a product kernel. 

\end{definition}

Theorem \ref{thm:main_LWC_SIRGs} can be directly adopted to the PSIRG setting under appropriate conditions on the product kernel.

\begin{assumption}[PSIRG connection functions]\label{sssec:assump_PSIRG_connections}
Let $\kappa\colon \mathbb{R}_+\times \mathbb{R} \times \mathbb{R} \to [0,1]$ be a product kernel $\kappa(t,x,y) = 1\wedge f(t)g(x,y)$ such that 
\begin{itemize}
    \item[1.]there exists a $\alpha_p > 1$ and $t_1 \in \mathbb{R}_+$ such that for all $t>t_1$
    \[f(t)\leq t^{-\alpha_p},\] 
    \item[2.]there exists a $\beta_p > 0$ and $t_2 \in \mathbb{R}_+$ such that for all $t> t_2$
    \[\Prob{g(W^{(1)},W^{(2)})>t}\leq t^{-\beta_p},\] where $W^{(1)}, W^{(2)}$ are i.i.d.\ copies of the limiting weight random variable $W$ in (\ref{eq:conv_in_law_Wn_to_W}).
\end{itemize}
\end{assumption}

Let $W^{(1)}$ and $W^{(2)}$ be two i.i.d.\ copies of $W$ (see Assumption \ref{sssec:assumption_weights}). The following lemma derives a useful bound on product kernels that satisfy Assumption \ref{sssec:assump_PSIRG_connections}:
\begin{lemma}[Bound on product kernel]
\label{lem:PSIRG_regular_variation}
    Under Assumption \ref{sssec:assump_PSIRG_connections}, {for any $\epsilon>0$, there exists $t_0=t_0(\epsilon)>0$,} such that whenever $t>t_0$, 
    \[
    \Exp{\kappa(t,W^{(1)},W^{(2)})} \leq t^{-\min\{\alpha_p,\alpha_p\beta_p\}+\epsilon}.
    \]
\end{lemma}

In other words, if the limiting connection function $\kappa$ in Theorem \ref{thm:main_LCinP_SIRGs} satisfies Assumption \ref{sssec:assump_PSIRG_connections}, it also satisfies Assumption \ref{sssec:assumption_connections} (3), with {$\alpha=\min\{\alpha_p,\alpha_p\beta_p\}-\epsilon$ for any $\epsilon>0$}. So, if 
\[
\gamma_p:=\min\{\alpha_p,\alpha_p\beta_p\}>d,
\]
then using Lemma \ref{lem:PSIRG_regular_variation}, we have that $\kappa$ satisfies Assumption \ref{sssec:assumption_connections} (3) with some $\alpha>d$ {, by choosing $\epsilon>0$ sufficiently small.} Hence, we obtain the following direct corollary to Theorem \ref{thm:main_LCinP_SIRGs}, whose proof we omit. Recall the vector $\mathbf{W}=(W_i)_{i \in \mathbb{N} \cup \{0\}}$ of i.i.d. copies of the limiting weight variable $W$ (see (\ref{eq:conv_in_law_Wn_to_W})), and $\mathbf{Y}=(Y_i)_{i \in \mathbb{N} \cup \{0\}}$ the atoms of $\Gamma_{\infty}$ (see (\ref{eq:defn_Gamma_infty})), with $Y_0=\mathbf{0}$.

\begin{corollary}[Convergence of PSIRGs locally in probability] \label{cor:PSIRG}
    Consider the sequence $(\mathbb{G}_n)_{n\geq 1}=(G(\mathbf{Y}^{(n)},\mathbf{W}^{(n)},\kappa_n))_{n \geq 1}$ of SIRGs, where for each $n$ and each $i\in[n]$, $Y^{(n)}_i=n^{1/d}X^{(n)}_i$, with $\mathbf{X}^{(n)}=(X^{(n)}_i)_{i\in[n]}$, $\mathbf{W}^{(n)}=(W^{(n)}_i)_{i\in[n]}$ satisfying Assumptions \ref{sssec:assumption_location}, \ref{sssec:assumption_weights} and $\kappa_n$ satisfying Assumption \ref{sssec:assumption_connections} (1), with $\kappa$ satisfying Assumption \ref{sssec:assump_PSIRG_connections}, with
    \[
    \gamma_p:=\min\{\alpha_p, \alpha_p\beta_p\} > d.
    \]
Then $(\mathbb{G}_n)_{n \geq 1}$ converges locally in probability to the infinite rooted SIRG $(\mathbb{G}_{\infty},0)$, rooted at vertex $0$, where $\mathbb{G}_{\infty}=G(\mathbf{Y},\mathbf{W},\kappa)$.
\end{corollary}

\begin{remark}[Regularly varying product forms]\label{rem:reg_var_psirg} Note that if $f(t)$ and $\Prob{g(W^{(1)},W^{(2)})>t}$ are regularly varying functions of $t$ outside some compact sets, with respective exponents $\alpha_p>0$ and $\beta_p>0$ with $\min\{\alpha_p,\alpha_p\beta_p\}>d$, then by Potter's theorem (see \cite{Potter}), they respectively satisfy Assumptions \ref{sssec:assump_PSIRG_connections} (1) and (2) with some exponents $\overline{\alpha_p}>0$ and $\overline{\beta_p}>0$, with $\min\{\overline{\alpha_p},\overline{\alpha_p}\overline{\beta_p}\}>d$. Hence for these kind of connection functions, Lemma \ref{lem:PSIRG_regular_variation}, and hence Corollary \ref{cor:PSIRG}, continues to be true.

\end{remark}

\begin{remark}[Dominance by PSIRG connections]\label{rem:dominance_PSIRG_connections} Note that even if $\kappa$ is not of product form, but instead is dominated by $1 \wedge f(t)g(x,y)$, with $f$ and $g$ respectively satisfying Assumptions \ref{sssec:assump_PSIRG_connections} (1) and (2), with $\gamma_p>d$, then Lemma \ref{lem:PSIRG_regular_variation}, and hence Corollary \ref{cor:PSIRG}, continues to be true.   

\end{remark}

Next, we discuss several known models which are all examples of PSIRGs, or SIRGs with connection functions dominated by PSIRG connection functions. {The results that follow, are presented as corollaries of Theorem \ref{thm:main_LCinP_SIRGs}, and their proofs are in Section \ref{ssec:proof_examples}.}

\subsubsection{Geometric Inhomogeneous Random Graphs}\label{sssec:examples_GIRGs}
Geometric Inhomogeneous Random Graphs (GIRGs) \cite{KB_Avg_dist_16, KB_SamplingGIRGs_17,KB_GIRGs_19} were motivated as spatial versions of the classic Chung-Lu random graphs \cite{CL02_1,CL02_2}. Although Chung-Lu random graphs with suitable parameters are scale free and exhibit small-world properties, they fail to capture clustering, a ubiquitous property of real-world networks.

GIRGs have four parameters: the number of vertices $n$, $\alpha_G \in (1,\infty]$, $\beta_G>2$ and the dimension $d\geq 1$. To each vertex $i \in \mathbb{N}$, one associates an independent, uniformly distributed random location $X^{(n)}_i$ on $I$ (recall (\ref{eq:defn_box_I})), and a real weight $W^{(n)}_i$ (possibly random), such that (\ref{eq:conv_in_law_Wn_to_W}) holds. Here the limiting weight variable $W$ has a power law tail with exponent $\beta_G$: there exists $t_G \in \mathbb{R}_+$ such that 
\begin{equation}\label{eq:power_law_limit_wt_GIRG}
    c_Gz^{1-\beta_G} \leq \Prob{W>z}\leq C_Gz^{1-\beta_G},
\end{equation}
whenever $z > t_G$, for some absolute constants $c_G, C_G>0$. In particular, we have $\Exp{W}< \infty$. Conditionally on $(X^{(n)}_i)_{i \in [n]}$ and $(W^{(n)}_i)_{i \in [n]}$, each edge $\{i,j\}$ is included independently with probability $p_{i,j}$ given by
\begin{equation}\label{GIRG_connection_fn}
    \begin{split}
        p_{i,j}=
        \begin{cases}
            1 \wedge \left(\frac{W^{(n)}_iW^{(n)}_j}{\sum_{i\in[n]} W^{(n)}_i}\right)^{\alpha_G}\frac{1}{\|X^{(n)}_i-X^{(n)}_j\|^{d\alpha_G}} , &\text{if}\; 1<\alpha_G<\infty; \\
            \ind{\left(\frac{W^{(n)}_iW^{(n)}_j}{\sum_{i\in[n]} W^{(n)}_i}\right)^{1/d}>\|X^{(n)}_i-X^{(n)}_j\|}, & \text{if}\; \alpha_G = \infty.
        \end{cases}
    \end{split}
\end{equation}
We denote the resulting random graph by
\[
\mathrm{GIRG}_{n,\alpha_G,\beta_G,d}.
\]

\begin{remark}[Relation to GIRGs in \cite{KB_GIRGs_19}]
Our formulation of GIRGs is closer to the formulation adopted in \cite{JK_BL_explosions_HRGs_19} than \cite{KB_GIRGs_19}. In the original definition of GIRGs (see \cite{KB_GIRGs_19}), the connection function is only assumed to be bounded above and below by multiples of (\ref{GIRG_connection_fn}), and vertex locations are assumed to be uniform on the torus $\mathbb{T}_n$, which is obtained by identifying the boundaries of $I$. However, to define a local limit, we need the connection function to converge to a limiting function, for which we have taken the explicit form (\ref{GIRG_connection_fn}). Finally, using the observation that only a negligible proportion of vertex locations fall near the boundary of $I$, our results can be easily transferred to the torus setting. 
\end{remark}

Now consider the following infinite SIRG $G(\mathbf{Y},\mathbf{W},\kappa^{\alpha_G})$, where $\mathbf{Y}=(Y_i)_{i \in \mathbb{N}\cup \{0\}}$ is the sequence of atoms of $\Gamma \cup \{\mathbf{0}\}$ (with $Y_0=\mathbf{0}$), $\Gamma$ is a unit-rate homogeneous Poisson point process on $\mathbb{R}^d$, $\mathbf{W}=(W_i)_{i \in \mathbb{N}\cup \{0\}}$ is an i.i.d.\ collection of limiting weight random variables $W$, and $\kappa^{(\alpha_G)}$ is the connection function
\begin{equation}\label{eq:GIRG_limiting_connection_fn}
    \begin{split}
        \kappa^{(\alpha_G)}(t,x,y)=
        \begin{cases}
            1 \wedge \left(\frac{xy}{\Exp{W}}\right)^{\alpha_G}{t^{-d\alpha_G}}, &\text{if}\; 1<\alpha_G<\infty; \\
            \ind{\left(\frac{xy}{\Exp{W}}\right)^{1/d}>t} , & \text{if}\; \alpha_G = \infty.
        \end{cases}
    \end{split}
\end{equation}

As a corollary to Theorem \ref{thm:main_LCinP_SIRGs}, we establish the local limit of the GIRG sequence to be $(G(\mathbf{Y},\mathbf{W},\kappa^{\alpha_G}),0)$. This answers a question posed in \cite{JK_BL_explosions_HRGs_19} (see \cite[Section 2.1]{JK_BL_explosions_HRGs_19}) in the affirmative. We call the above infinite SIRG, the \emph{infinite GIRG}, and denote it as $\mathrm{GIRG}_{\infty,\alpha_G,\beta_G,d}$.

\begin{corollary}[Convergence of GIRGs locally in probability] \label{cor:LWP_GIRGs}
As $n\rightarrow \infty$, the sequence \\ $(\mathrm{GIRG}_{n,\alpha_G,\beta_G,d})_{n \geq 1}$ converges locally in probability to the rooted infinite GIRG $(\mathrm{GIRG}_{\infty,\alpha_G,\beta_G,d},0)$, rooted at $0$, where $\alpha_G \in (1,\infty]$, $\beta_G > 2$, $d \in \mathbb{N}$.
\end{corollary}

\subsubsection{Hyperbolic Random Graphs}\label{sssec:examples_HRG}

Hyperbolic Random Graphs (HRGs) were first proposed by Krioukov et al.\ in 2010 \cite{KPKVB_HRG}, as a model that captures the three main properties of most real-world networks: scale free, small distances, and non-vanishing clustering coefficient. 

\medskip

HRGs have three parameters, namely the number of vertices $n$,  $\alpha_H>\frac{1}{2}$, and $\nu>0$, which are fixed constants. Let
\begin{align*}
    R_n:=2\log{\frac{n}{\nu}}. \numberthis \label{defn:HRG_radius}
\end{align*}
The vertex set of the graph is the set of $n$ i.i.d.\ points $u^{(n)}_1,\ldots,u^{(n)}_n$ on the hyperbolic plane $\mathbb{H}$, where $u^{(n)}_i=(r^{(n)}_i,\theta^{(n)}_i)$ is the polar representation of $u^{(n)}_i$. The angular component vector $(\theta^{(n)}_i)_{i=1}^n$ is a vector with i.i.d.\ coordinates, each coordinate having the uniform distribution on $[-\pi,\pi]$. The radial component vector $(r^{(n)}_i)_{i=1}^n$ is independent of $(\theta_i)_{i=1}^n$, and has i.i.d.\ coordinates, with cumulative distribution function
\begin{equation}\label{eq:radial_CDF_HRG}
    F^{(n)}_{\alpha_H,\nu}(r)=
    \begin{cases}
        0 &\text{if}\;\;r<0, \\
        \frac{\cosh{\alpha_H r-1}}{\cosh{\alpha_H R_n -1}} &\text{if}\;\;0\leq r \leq R_n, \\
        1 & \text{if}\;\; r> R_n.
    \end{cases}
\end{equation}

Given $(u^{(n)}_i)_{i=1}^n=((r^{(n)}_i,\theta^{(n)}_i))_{i=1}^n$, one forms the \emph{threshold} HRG (THRG) by {placing edges between all pairs of vertices $u^{(n)}_i$ and $u^{(n)}_j$ with conditional probability
\begin{align*}
    p^{(n)}_{\mathrm{THRG}}(u^{(n)}_i,u^{(n)}_j):=\ind{d_{\mathbb{H}}(u^{(n)}_i,u^{(n)}_j)<R_n}, \numberthis \label{eq:THRG_conn}    
\end{align*}
where $d_{\mathbb{H}}$ denotes the distance in the hyperbolic plane $\mathbb{H}$, i.e., the edge between $u^{(n)}_i$ and $u^{(n)}_j$ is included if and only if $d_{\mathbb{H}}(u^{(n)}_i,u^{(n)}_j)<R_n$.  }

 {Similarly,} one forms a parametrized version of the THRG (see \cite[Section \RNum{6}]{KPKVB_HRG}) which we call the \emph{parametrized} HRG (PHRG) by placing edges independently between $u^{(n)}_i$ and $u^{(n)}_j$ with conditional probability
 \begin{align*}
    {p^{(n)}_{\mathrm{PHRG}}(u^{(n)}_i,u^{(n)}_j)}:=\left(1 + \exp{\left(\frac{d_{\mathbb{H}}(u^{(n)}_i,u^{(n)}_j)-R_n}{2T_H}\right)}\right)^{-1}, \numberthis \label{eq:PHRG_conn} 
 \end{align*}
 {where $T_H>0$ is another parameter.}

We denote the THRG model with parameters $n$, $\alpha_H$ and $\nu$  {by} $\mathrm{THRG}_{n,\alpha_H,\nu},$
and the PHRG model with parameters $n,\alpha_H,T_H,\nu$  {by} $\mathrm{PHRG}_{n,\alpha_H,T_H,\nu}.$
\medskip



    

{Both the THRG and PHRG models can be seen as finite SIRGs, which gives us the local limits for these models, as a corollary to Theorem \ref{thm:main_LCinP_SIRGs}.}


\begin{corollary}[Convergence of HRGs locally in probability]\label{cor:LWP_HRGs}
Let $\alpha_H> {\tfrac{1}{2}}$, $0<T_H<1$, and $n \in \mathbb{N}$. {Let $\mathbf{Y}$ be the sequence of atoms of (\ref{eq:defn_Gamma_infty}). Then there exists a random variable $W$ having a power-law distribution with exponent $2\alpha_H+1$, such that if $\mathbf{W}=(W_i)_{i \in \mathbb{N}\cup\{0\}}$ is a sequence of i.i.d.\ copies of $W$,}
\begin{itemize}
\item[(a)] the sequence $(\mathrm{THRG}_{n,\alpha_H,\nu})_{n \geq 1}$ converges locally in probability to {the infinite SIRG \\
$G(\mathbf{Y},\mathbf{W},\kappa_{\mathrm{THRG},\infty})$, where
\begin{align*}
    \kappa_{\mathrm{THRG},\infty}(t,x,y):= \ind{t \leq \frac{\nu xy}{\pi}}, 
\end{align*}
} 
\item[(b)] the sequence $(\mathrm{PHRG}_{n,\alpha_H,T_H,\nu})_{n \geq 1}$ converges locally in probability to {the infinite SIRG\\
$G(\mathbf{Y},\mathbf{W},\kappa_{\mathrm{PHRG},\infty})$, where
\begin{align*}
    \kappa_{\mathrm{PHRG},\infty}(t,x,y):= \left(1+\left(\frac{ \pi t}{\nu xy} \right)^{1/T_H} \right)^{-1}.
\end{align*}}
\end{itemize}
\end{corollary}


\subsubsection{Continuum Scale-Free Percolation}
The continuum scale-free percolation (CSFP) model \cite{CSFP_D_W} was introduced as a continuum analogue of the discrete scale-free percolation (SFP) model \cite{SFP}, a model motivated by capturing power-law degree distributions, while preserving non-zero clustering and logarithmic typical distances. 

We now formally define the model following \cite{CSFP_D_W}. The vertex set is the set of points of a homogeneous Poisson point process $(Y_i)_{i \in \mathbb{N}}$, marked with i.i.d.\ weights $(W_i)_{i \in \mathbb{N}}$, which  {have} a Pareto distribution with power-law tail parameter $\beta>0$ and scale parameter $1$:
\begin{align*}
\Prob{W>w}=w^{- \beta}, \numberthis \label{eq:wt_law_CSFP}
\end{align*}
whenever $w>1$.  {Conditionally} on $(Y_i)_{i \in \mathbb{N}}$ and $(W_i)_{i \in \mathbb{N}}$, each edge $\{Y_i,Y_j\}$ is included independently with probability
\[
1-\exp{\left(-\frac{\lambda W_iW_j}{\|Y_i-Y_j\|^{\alpha}}\right)},
\]
where $\lambda>0$ is a parameter. Here we remark that in the original definition, instead of a homogeneous Poisson process $(Y_i)_{i \in \mathbb{N}}$, a Poisson process with some constant intensity $\nu>0$ was considered in \cite{CSFP_D_W}. By standard scaling arguments, this does not make any difference in our results.

Considering the Palm version $(Y_i)_{i \in \mathbb{N}\cup \{0\}}$ of $(Y_i)_{i \in \mathbb{N}}$, where $Y_0=\mathbf{0} \in \mathbb{R}^d$, marking $Y_0$ with  {an independent} weight $W_0$, and rooting the resulting graph at $0$, it is immediate that the resulting rooted infinite CSFP model is the rooted SIRG $(G(\mathbf{Y},\mathbf{W},\kappa),0)$, where $\mathbf{Y}=(Y_i)_{i \in \mathbb{N}\cup \{0\}}$, $\mathbf{W}=(W_i)_{i \in \mathbb{N}\cup\{0\}}$, and $\kappa(t,x,y):=1-\exp{\left(-\frac{\lambda xy}{t^{\alpha}}\right)}
$. 

 {For each $n \geq 1$, let} $\mathbf{X}^{(n)}=(X^{(n)}_i)_{i \in [n]}$ satisfy Assumption \ref{sssec:assumption_location}, and consider the `blown up' finite CSFP model $G(\mathbf{Y}^{(n)},\mathbf{W}^{(n)},\kappa_n)$, where $Y^{(n)}_i=n^{1/d}X^{(n)}_i$ for $i \in [n]$, $\mathbf{W}^{(n)}=(W^{(n)}_i)_{i \in [n]}$ is a vector of i.i.d.\ weight variables having law (\ref{eq:wt_law_CSFP}), and $\kappa_n(t,x,y)=1-\exp{\left(-\frac{\lambda xy}{n^{-\frac{\alpha}{d}}t^{\alpha}}\right)}$. Then as corollary to Theorem \ref{thm:main_LCinP_SIRGs}, we have that the infinite CSFP is the local limit of finite CSFPs under suitable assumptions. 

\begin{corollary}[Convergence of CSFPs locally in probability]
\label{cor:LCinP_CSFP}
    Let 
    \[
    \min\{\alpha,\alpha\beta\}>d.
    \]
    Then, as $n \to \infty$, the graph sequence $G(\mathbf{Y}^{(n)},\mathbf{W}^{(n)},\kappa_n)$ converges locally in probability to the infinite rooted CSFP $(G(\mathbf{Y},\mathbf{W},\kappa),0)$. 
\end{corollary}

\subsubsection{Weight dependent Random Connection Models}
Another very general class of spatial random graph model called the weight dependent random connection model (WDRCM) was first introduced in \cite{WDRCM_Rec_trans}, motivated by the study of recurrence and transience properties of general geometric graphs, which we briefly discuss.

To construct the graph, one takes a unit {-rate} Poisson process on $\mathbb{R}^d \times[0,1]$, conditionally on which, edges between pairs of vertices $(\mathbf{x},s)$ and $(\mathbf{y},t)$ are placed independently with probability 
\[
\rho(h(s,t,\|\mathbf{x}-\mathbf{y}\|)),
\]
for some \emph{profile} function $\rho:\mathbb{R}_+ \to [0,1]$, and a suitable \emph{kernel} $h:[0,1]\times[0,1]\times \mathbb{R}_+ \to \mathbb{R}_+$. The vertex $(\mathbf{x},s)$ is thought of being located at $\mathbf{x} \in \mathbb{R}^d$, and having a weight of $s^{-1}$ associated to it.

Including the point $\mathbf{r}=(\mathbf{0},U) \in \mathbb{R}^d \times [0,1]$, where $U$ is uniformly distributed on $[0,1]$ and independent of the weights and locations of other vertices, and rooting the graph at $\mathbf{r}$, we obtain the infinite rooted SIRG $(G(\mathbf{Y},\mathbf{W},\kappa),0)$, where $\mathbf{Y}=(Y_i)_{i \in \mathbb{N}\cup \{0\}}$ are the atoms of (\ref{eq:defn_Gamma_infty}), $\mathbf{W}=(W_i)_{i \in \mathbb{N}\cup\{0\}}$ is a sequence of i.i.d.\ uniform on $[0,1]$ random weights, and $\kappa = \rho \circ h$.

 {For} $n \geq 1$, let $(X^{(n)}_i)_{i \in [n]}$ satisfy Assumption \ref{sssec:assumption_location}, $(W^{(n)}_i)_{i \in [n]}$ be a collection of $n$ i.i.d.\ uniform on $[0,1]$ weights, and let $\kappa_n:\mathbb{R}\times \mathbb{R} \times \mathbb{R}_+ \to [0,1]$ be defined as $\kappa_n(t,x,y)=\rho(h(x,y,n^{-1/d}t))$, with $\kappa(t,x,y)=\rho(g(x,y,t))$ satisfying Assumption \ref{sssec:assumption_connections} (2). Then as a direct consequence of Theorem \ref{thm:main_LCinP_SIRGs} we have the following corollary whose proof we omit:

\begin{corollary}[Convergence of WDRCMs locally in probability] Let $\mathbf{Y}^{(n)}=(Y^{(n)}_i)_{i \in [n]}=(n^{1/d}X^{(n)}_i)_{i \in [n]}$. Then, as $n \to \infty$, the sequence of SIRGS $G(\mathbf{Y}^{(n)},\mathbf{W}^{(n)},\kappa_n)$ converges locally in probability to the infinite rooted WDRCM $(G(\mathbf{Y},\mathbf{W},\kappa),0)$. 

\end{corollary}

\subsection{Consequences of local convergence: Degrees}
Theorem \ref{thm:main_LCinP_SIRGs} is equivalent to the statement that for any subset $A \subset \mathcal{G}_{\star}$,     
\[\frac{1}{n}\sum_{i=1}^n \ind{(\mathbb{G}_n,i)\in A}\]
converges in probability to $\Prob{(\mathbb{G}_{\infty},0)\in A}$ (see e.g., \cite[(2.3.5)]{RGCN_2}).

In particular, for fixed $k\in \mathbb{N}$, one can take $A_k$ to be the subset of those rooted graphs $(G,o)$ for which the root $o$ has degree $k$ in $G$, to conclude that
\[
\frac{N_k(\mathbb{G}_n)}{n} \plim \Prob{D=k},
\]
where $N_k(\mathbb{G}_n)$ is the number of vertices with degree $k$ in $\mathbb{G}_n$, and where $D$ is the degree of $0$ in $\mathbb{G}_{\infty}$.

Taking expectation and applying dominated convergence, we have $\Prob{D_n=k}\to \Prob{D=k}$, as $n \to \infty$, where $D_n$ is the degree of $U_n$ in $\mathbb{G}_n$, for any $ {k} \in \mathbb{N}$, which implies
\begin{equation}\label{eq:D_n_dlim_D}
    D_n \dlim D,
\end{equation}
(where $\dlim$ means convergence in distribution) as $n \to \infty$.
We next give a description of the random variable $D$.  {For this, let $W_0$ be the weight of $0$ in $\mathbb{G}_{\infty}$, and $W^{(1)}$ be an independent copy of $W_0$.}

\begin{proposition}[Degree distribution]\label{prop:degree_law}
Under the assumptions of Theorem \ref{thm:main_LWC_SIRGs}, the random variable $D$ is 
\[
\text{\rm Poi}\left(\int_{\mathbb{R}^d} \CExp{\kappa(\|z\|,W_0,W^{(1)})}{W_0}dz \right)
\] 
distributed, i.e. it has a mixed Poisson distribution with mixing parameter 
\[
\int_{\mathbb{R}^d} \CExp{\kappa(\|z\|,W_0,W^{(1)})}{W_0}dz.
\]   
\end{proposition}
 {Note that} Proposition \ref{prop:degree_law}, Assumption \ref{sssec:assumption_connections} (3) and Fubini's theorem imply that
\[
\Exp{D}=\int_{\mathbb{R}^d}\Exp{\kappa(\|x\|,W_0,W^{(1)})}dx,
\]
which is finite. Since the expectation of the mixing parameter is finite, the mixing parameter  {is} finite almost surely, and hence Proposition \ref{prop:degree_law} makes sense.

In particular, when the mixing parameter $\int_{\mathbb{R}^d} \CExp{\kappa(\|z\|,W_0,W^{(1)})}{W_0}dz$ is regularly varying with some exponent $\zeta>0$, the random variable $D$ is also regularly varying with the same exponent $\zeta>0$. This allows for the existence of power-law degree distributions in spatial random graphs.

Recall that $D_n$ is the degree of the uniform vertex $U_n$ in $\mathbb{G}_n$. 

\begin{proposition}[Uniform integrability of typical degree sequence]\label{prop:UI_typical_degs}
Under the assumptions of Theorem \ref{thm:main_LWC_SIRGs}, the sequence $(D_n)_{n \geq 1}$ is a uniformly integrable sequence of random variables.
\end{proposition}

 {The proof is given in Section~\ref{ssec:proof_degree}. This} result is of independent interest. Uniform integrability of the typical degree sequence does in general not follow from local convergence, even when the limiting degree distribution has finite mean, see for example \cite[Exercise 2.14]{RGCN_2}.


Combining Proposition \ref{prop:UI_typical_degs} with (\ref{eq:D_n_dlim_D}), we note that \[
\Exp{D_n} \to \Exp{D},
\] 
as $n \to \infty$. Note that one cannot directly conclude this from Theorem \ref{thm:main_LWC_SIRGs}, because the function
$\mathrm{D}:\mathcal{G}_{\star}\to \mathbb{R}_+$ defined by 
\[
    \mathrm{D}((G,o)):=\text{degree of}\;\;o\;\;\text{in}\;\;G,
\]
is continuous, but not necessarily bounded.

\subsection{Consequences of local convergence: Clustering}

In this section, we discuss convergence of various clustering measures of SIRGs.

For any graph $G=(V(G),E(G))$, we let 
\begin{equation}\label{eq:2_wedge}
   \mathcal{W}_{G}:= \sum_{v_1,v_2,v_3 \in V(G)}\ind{\{v_1,v_2\},\{v_2,v_3\}\in E(G)}=\sum_{v \in V(G)} d_v(d_v-1), 
\end{equation}
(where $d_v$ is the degree of $v$ in $G$) be twice the number of \emph{wedges} in the graph $G$, and 
\begin{equation} \label{eq:3_triangles}
    \Delta_{G}:= \sum_{v_1,v_2,v_3 \in V(G)} \ind{\{v_1,v_2\},\{v_2,v_3\},\{v_3,v_1\} \in E(G)}
\end{equation}
be {six times} the number of \emph{triangles} in the graph $G$, {where the sums in (\ref{eq:2_wedge}) and (\ref{eq:3_triangles}) are over distinct vertices $v_1,v_2,v_3 \in V(G)$.}

Then the {\em global clustering coefficient} $\mathrm{CC}_{G}$ of the graph $G$ is defined as
\begin{equation}\label{eq:defn_GCC}
    \mathrm{CC}_{G}:= \frac{\Delta_{G}}{\mathcal{W}_{G}}.
\end{equation}
 {We next discuss a local notion of clustering. Define for $v \in V(G)$,
\[
    \mathrm{CC}_{G}(v):= \begin{cases}
        \frac{\Delta_v(G)}{d_v(d_v-1)} &\mbox{if } d_v \ge 2, \\
        0 &\mbox{else,}
    \end{cases}
\]
}
where 
\[
\Delta_v(G)=\sum_{v_1,v_2 \in V(G)}\ind{\{v_1,v\},\{v_2,v\},\{v_1,v_2\}\in E(G)}
\] 
is twice the number of triangles in $G$ containing the vertex $v$, and $d_v$ is as before the degree of $v$ in $G$.
 {The {\em local clustering coefficient} $\overline{\mathrm{CC}}_{G}$ of $G$ is then defined as}
\begin{equation}\label{eq:defn_LCC}
    \overline{\mathrm{CC}}_{G}:= \frac{1}{n}\sum_{v \in V(G)}\mathrm{CC}_{G}(v).
\end{equation}

Finally, we discuss a notion of clustering contribution from only vertices of certain degree. For $k \in \mathbb{N}$, define the {\em clustering function} to be
\[
k \mapsto \mathrm{CC}_{G,k},
\]
where $\mathrm{CC}_{G,k}$ is defined as
\begin{equation}\label{eq:defn_ClF}
    \mathrm{CC}_{G,k}:= \begin{cases}
        \frac{1}{N_k(G)}\sum_{v \in v(G),d_v=k} {\frac{\Delta_v(G)}{k(k-1)}} &\text{if}\;N_k(G)>0,\\
        0&\text{otherwise},
    \end{cases}
\end{equation}
where $N_k(G)$ is the total number of vertices in $G$ with degree $k$.  {Thus}, $\mathrm{CC}_{G,k}$ measures the proportion of wedges that are triangles, where one of the participant vertices has degree $k$.

We now present the results on convergence of these various clustering measures for SIRGs:

\begin{corollary}[Convergence of clustering coefficients of SIRGs]\label{cor:conv_clustering}
Under the assumptions of Theorem \ref{thm:main_LWC_SIRGs}, as $n \to \infty$,

\begin{itemize}
    \item[1.] if $\alpha> 2d$, then 
    \begin{equation}\label{{eq:conv_GCC}}
        \mathrm{CC}_{\mathbb{G}_n} \plim \frac{\Exp{\Delta_0}}{\Exp{D(D-1)}},
    \end{equation}
    where $\Delta_0:=\sum_{i,j \in \mathbb{N}}\ind{\{0,i\},\{0,j\},\{i,j\}\in E(\mathbb{G}_{\infty})}$ is twice the number of triangles containing $0$ in $\mathbb{G}_{\infty}$, and $D$ is the degree of $0$ in $\mathbb{G}_{\infty}$.
    
    \item[2.]
    \begin{equation}\label{eq:conv_LCC}
        \overline{\mathrm{CC}}_{\mathbb{G}_n} \plim \Exp{\frac{\Delta_0}{D(D-1)}}.
    \end{equation}
    
    \item[3.] for any $k \in \mathbb{N}$,
    \begin{equation}\label{eq:conv_ClF}
        \mathrm{CC}_{G,k} \plim \frac{1}{\binom{k}{2}}\CExp{\Delta_0}{D=k}.
    \end{equation}
    
\end{itemize}
\end{corollary}

 {Corollary~\ref{cor:conv_clustering} (2) and (3) are direct consequences of local convergence (see \cite[Section 2.4.2]{RGCN_2}). For Corollary~\ref{cor:conv_clustering} (1), we need an additional uniform convergence property of the {\em square} of the degree of a uniform vertex, which we prove is implied by the condition $\alpha>2d$.}

\begin{remark}[Condition on $\alpha$]
The condition $\alpha>2d$ in Corollary \ref{cor:conv_clustering} (1) is not optimal, as will be evident in the proof. Our purpose is not to find the optimal conditions under which the global clustering coefficient converges, but to demonstrate how local convergence of graphs  {implies convergence} of the global clustering coefficient.
\end{remark}

Recently, precise results about convergence of clustering coefficients, and scaling of the clustering function as $k$ grows to infinity, for Hyperbolic Random Graphs has been obtained in \cite{FHMS_clustering_HRGs}. Also, it was shown in \cite{CSFP_clustering} that under suitable conditions, the CSFP model has non-zero clustering in the limit.

\subsection{Consequences of our local convergence proof: Distance  {lower bound}}\label{ssec:disussion_dist}
Finally, we provide a result on typical distances in our graphs. Let $U_{n,1}$ and $U_{n,2}$ be two i.i.d.\ uniformly distributed vertices of $\mathbb{G}_n$, so that $d_{\mathbb{G}_n}(U_{n,1},U_{n,2})$ is the graph distance in $\mathbb{G}_n$ between $U_{n,1}$ and $U_{n,2}$. Recall  {that} by convention we let $d_{\mathbb{G}_n}(U_{n,1},U_{n,2})=\infty$  {when} $U_{n,1}$ and $U_{n,2}$ are not in the same connected component of $\mathbb{G}_n$, so that $d_{\mathbb{G}_n}(U_{n,1},U_{n,2})$ is a well-defined random variable.  

\begin{theorem}[Lower bound on typical distances]\label{thm:main_typical_dist}
Under the assumptions of Theorem \ref{thm:main_LWC_SIRGs}, for any $C \in \left( 0, \frac{1}{\log(\frac{\alpha}{\alpha-d})}\right)$,
\[ \Prob{d_{\mathbb{G}_n}(U_{n,1},U_{n,2})> C \log \log n} \to 1,\]
as $n \to \infty$.
\end{theorem}

We in fact believe the limit in the above display holds with $C$ replaced by $\frac{1}{\log(\frac{\alpha}{\alpha-d})}$. But our proof method does not allow us to establish this improvement. As we will see, the proof is a direct by-product of the proof of the local weak limit  {in} Theorem \ref{thm:main_LWC_SIRGs}. Note that as $\alpha$ approaches $d$, the lower bound in Theorem \ref{thm:main_typical_dist} becomes trivial.

If instead of a regularly varying domination as in Assumption \ref{sssec:assumption_connections} (2), $\Exp{\kappa(t,W_0,W_1)}$ itself is regularly varying in $t$ with exponent $\alpha$, then it follows from Proposition \ref{prop:degree_law} that the expectation of the limiting degree distribution $D$ is infinite in the regime $\alpha \in (0,d)$. We conjecture the distances are of constant order in this regime:
\begin{conjecture}[Constant distances for {$\alpha \in (0,d)$}]
Let $\mathbb{G}_n=G(\mathbf{Y}^{(n)},\mathbf{W}^{(n)},\kappa_n)$ satisfy the assumptions of Theorem \ref{thm:main_LWC_SIRGs} except Assumption \ref{sssec:assumption_connections} (2), where instead we assume the limiting connection function $\kappa$ is such that $\Exp{\kappa(t,W_0,W_1)}$ is regularly varying with exponent $\alpha \in (0,d)$. Then if $U_{n,1}$ and $U_{n,2}$ be two uniformly chosen vertices in the SIRG $\mathbb{G}_n=G(\mathbf{Y}^{(n)},\mathbf{W}^{(n)},\kappa_n)$, conditionally on the event that $U_{n,1}$ and $U_{n,2}$ are connected in $\mathbb{G}_n$,
\[
d_{\mathbb{G}_n}(U_{n,1},U_{n,2}) \plim K(\alpha,d),
\]
as $n \to \infty$, where $K(\alpha,d)$ is a constant depending only on the exponent $\alpha$ and dimension $d$.
\end{conjecture}

 {For the case $\alpha=d$, we {do} not expect universal behaviour, and the question then becomes model dependent.} 
 {Results for constant distances when the limiting degree distribution has infinite mean are known for} lattice models such as Long range Percolation, see \cite[Example 6.1]{Geom_USF_BKPS}, for Scale Free Percolation, see \cite[Theorem 2.1]{SC_SFP_JHH},  {and for the configuration model \cite{RvdH_GH_vdE_CM_const_dist}, which is a model without geometry.}

{Theorem \ref{thm:main_typical_dist} poses the question, when are the typical distances exactly doubly logarithmic? Interestingly, distances can be larger than doubly logarithmic, even when the limiting degree distribution has infinite second moment, $\Exp{D^2}=\infty$, as was shown in \cite{chem_dist}, see for example \cite[Theorem 1.1 (a)]{chem_dist}. We conjecture that certain lower order moments below a (model dependent) critical threshold being infinite imply ultra-small distances. This is also the behaviour that the authors of \cite{chem_dist} observe, for a special class of models, but we believe this behaviour is universal:}

\begin{conjecture}[Ultra-small distances]
\label{conj:UB_typ_dist}
Under the assumptions of Theorem \ref{thm:main_typical_dist}, {there is a constant $\varepsilon_{\ast} \in (0,1)$ depending on the model parameters, such that for any $\varepsilon>\varepsilon_{\ast}$, $\Exp{D^{2-\varepsilon}}=\infty$ implies that} there is a constant $C(\alpha,d)>0$ such that
\[
\CProb{d_{\mathbb{G}_n}(U_{n,1},U_{n,2})< C(\alpha,d) \log \log n}{U_{n,1}\;\text{and}\;U_{n,2}\;\text{are connected in}\;\mathbb{G}_n} \to 1,
\]
as $n \to \infty$.
\end{conjecture}

\section{Proofs}\label{sec:proofs}

In this section we give all the proofs. We first start  {in} Section \ref{ssec:proof_notations}  {by defining the notation} that we use throughout this section, and  {by outlining the} general proof strategy of the main Theorems \ref{thm:main_LWC_SIRGs} and \ref{thm:main_LCinP_SIRGs}. Sections \ref{ssec:Proof_prop_loc_trunc} and \ref{ssec:proof_path_counting} contain proofs of some of the key tools  {that} we employ to prove our main results. The proofs of Theorems \ref{thm:main_LWC_SIRGs} and \ref{thm:main_LCinP_SIRGs} can be found in Sections \ref{ssec:proof_LWC} and \ref{ssec:proof_LCinP} respectively.  {The proof} of Theorem \ref{thm:main_typical_dist} can be found in Section \ref{ssec:proof_typ_dist}. Proofs of results on examples covered under our setup are in Section \ref{ssec:proof_examples}. The proofs of degree and clustering results can be found respectively in Sections \ref{ssec:proof_degree} and \ref{ssec:proof_clustering}.

\subsection{Notations and general proof strategy for Theorems \ref{thm:main_LWC_SIRGs} and \ref{thm:main_LCinP_SIRGs}}\label{ssec:proof_notations}

Recall the SIRGs $\mathbb{G}_n=G(\mathbf{Y}^{(n)},\mathbf{W}^{(n)},\kappa_n)$ and $\mathbb{G}_{\infty}=G(\mathbf{Y},\mathbf{W},\kappa)$ from Theorem \ref{thm:main_LWC_SIRGs}.
We first define some notations which we will use throughout. Recall $\mathbb{G}_n$ has vertex set $V(\mathbb{G}_n)=[n]$, and $\mathbb{G}_{\infty}$ has vertex set $V(\mathbb{G}_{\infty})=\mathbb{N}\cup\{0\}$.

For $r>0$, define the set
 \begin{equation}\label{eq:defn_set_Ar_n}
   A^r_n
   :=\left[-\frac{n^{1/d}}{2}+r, \frac{n^{1/d}}{2}-r\right]^d.  
 \end{equation}
 {Thus,} $A^r_n$ is a sub-box of the box $I_n$,  {such that for any point in $A^r_n$ the open Euclidean ball of radius $r$ around that point is contained in the box $I_n$. Hence, the number of points of the binomial process $\Gamma_n$ (recall (\ref{eq:defn_Gamma_n})) falling in this open ball}, has the same distribution as the number of points of $\Gamma_n$ falling in the open ball of radius $r$ around the origin $\mathbf{0} \in \mathbb{R}^d$. 
We will use this property of $A^r_n$ in a suitable manner, which we formally explain next.
 
To this end, for any $x \in \mathbb{R}^d$, define the ball
 \begin{equation}\label{eq:defn_ball_L_r_x}
   \mathscr{B}^r_x:=\{y\in \mathbb{R}^d: \|x-y\|<r\}.
 \end{equation}
Then if we let $\partial(I_n)=I_n\setminus\text{int}(I_n)$ denote the boundary of the set $I_n$, where $\text{int}(I_n)$ is the interior of $I_n$, i.e. the union of all open subsets of $I_n$, we note that for any  vertex $j \in V(\mathbb{G}_n)$ with its location $Y^{(n)}_j \in A^r_n$, {the ball $\mathscr{B}^r_{Y^{(n)}_j}$ does not intersect the boundary $\partial(I_n)$ of $I_n$, i.e. $\mathscr{B}^r_{Y^{(n)}_j} \subset \text{int}(I_n)$}.   {As a result,} the distribution of the number of vertices of $\mathbb{G}_n$ (other than $j$) with locations in $\mathscr{B}^r_{Y^{(n)}_j}$ does not depend on $Y^{(n)}_j$, and follows a $\text{Bin}\left(n-1, \frac{\lambda_d(\mathscr{B}^r_{\mathbf{0}})}{n}\right)$ distribution (where $\lambda_d$ denotes the Lebesgue measure on $\mathbb{R}^d$). 

In particular, since $\frac{\lambda_d(I_n \setminus A^r_n)}{n}\to 0$ as $n \to \infty$, the location $Y^{(n)}_{U_n}$ of the uniformly chosen vertex $U_n$ of $\mathbb{G}_n$ will with high probability fall in $A^r_n$. We will condition on this good event, under which the number of points of $\Gamma_n$ in $\mathscr{B}^r_{Y^{(n)}_{U_n}}$ follows a $\text{Bin}\left(n-1, \frac{\lambda_d(\mathscr{B}^r_{\mathbf{0}})}{n}\right)$ distribution, and this will  {simplify} our computations.

\begin{definition}[Euclidean graph neighborhoods  {around} a vertex]\label{defn:F_graphs}
    For a vertex $i \in V(\mathbb{G}_n)$, 
we define $(F^{\mathbb{G}_n}_i(r),i)$
to be the rooted subgraph of $\mathbb{G}_n$ rooted at $i$, induced by those vertices $j$, whose locations $Y^{(n)}_j$ satisfy $Y^{(n)}_j \in \mathscr{B}^r_{Y^{(n)}_i}$.
 
Similarly, we define $(F^{\mathbb{G}_{\infty}}_0(r),0)$
to be the rooted subgraph of $G_{\infty}$ rooted at $0$, induced by the vertices $j \in V(\mathbb{G}_{\infty})$, whose locations $Y_j$ satisfy $Y_j \in \mathscr{B}^r_{Y_0}=\mathscr{B}^r_{\mathbf{0}}$.
 
\end{definition} 
For $i \in [n]$, we will sometimes abbreviate the rooted graph $(F^{\mathbb{G}_n}_i(r),i)$ as simply $F^{\mathbb{G}_n}_i(r)$, and similarly for $F^{\mathbb{G}_{\infty}}_0(r)$.

For any graph $G=(V(G), E(G))$, edge $e=\{v_1,v_2\} \in E(G)$, and vertex $v \in V(G)$, by the \emph{ graph distance} of the edge $e {=\{v_1,v_2\}}$ from $v$, we mean the number 
\begin{equation}\label{eq:defn_dist_edge_vertex}
     \min\{d_{G}(v_1,v),d_{G}(v_2,v)\},
\end{equation}
where $d_G$ is the graph distance on $G$.
\medskip

Having introduced the main notations, we next discuss the main ingredients and the proof strategy for Theorems \ref{thm:main_LWC_SIRGs} and \ref{thm:main_LCinP_SIRGs}.

\paragraph{Local convergence of Euclidean graph neighborhoods.} Recall the graphs $F^{\mathbb{G}_n}_{U_n}(r)$ and $F^{\mathbb{G}_{\infty}}_{0}(r)$ from Definition \ref{defn:F_graphs}.  {In Section \ref{ssec:Proof_prop_loc_trunc}, we will prove} that the typical local graph structure in any deterministic Euclidean ball around the root location is asymptotically what it should be, i.e., for any $r$, $(F^{\mathbb{G}_n}_{U_n}(r),U_n)$ is close in distribution to  $(F^{\mathbb{G}_{\infty}}_0(r),0)$:

\begin{proposition}[Local convergence of Euclidean graph neighborhoods]\label{prop:loc_trunc}
For any fixed rooted finite graph $H_*=(H ,h) \in \mathcal{G}_{\star}$, and for any $r>0$, 
\[\Prob{(F^{\mathbb{G}_n}_{{U_n}}(r),U_n) \cong  (H,h)}\to \Prob{(F^{\mathbb{G}_{\infty}}_{0}(r),0)\cong  (H,h)},\] 
as $n \to \infty$, where $U_n$ is uniformly distributed on $V(\mathbb{G}_n)=[n]$. 
\end{proposition}


\paragraph{Path-counting analysis.} Next, in Section \ref{ssec:proof_path_counting}, we do a path-counting analysis. We begin by proving a technical lemma  {that will help us in implementing this path-counting analysis}. To state this lemma, we first introduce some more notations to keep things neat. Recall  {that} $U_n$ is uniformly distributed on $[n]$.

For $n,j \in \mathbb{N}, v_1,\ldots,v_j \in [n]$ and $\vec{x}= (x_0,\ldots,x_j) \in (\mathbb{R}^d)^{j+1}$, we denote 
 \begin{equation}\label{eq:defn_W_n}
    \mathbb{W}^{v_1,\ldots,v_j}_{n}(\vec{x})
    :=\Exp{\kappa_n\left(\|x_1-x_0\|,W^{(n)}_{U_n},W^{(n)}_{v_1}\right)\cdots\kappa_n\left(\|x_j-x_{j-1}\|,W^{(n)}_{v_{j-1}},W^{(n)}_{v_j}\right)},
\end{equation}
and for $j \in \mathbb{N}$, $u_1,\dots,u_j \in \mathbb{N}$, $x_0,\ldots,x_j \in \mathbb{R}^d$, we denote 
\begin{equation}\label{eq:defn_W}
    \mathbb{W}_j(\vec{x}):=\Exp{\kappa\left(\|x_{1}-x_{0}\|,W_0,W_{u_1} \right)\cdots\kappa\left(\|x_{{j}}-x_{{j-1}}\|,W_{u_{j-1}},W_{u_j}\right)}.
\end{equation}
 {Note that the values of the expectations on the RHS of {(\ref{eq:defn_W_n}) and} (\ref{eq:defn_W}) does not depend on the values of the {$v_i$'s and} $u_i$'s {respectively}.}

Then,  {our main path-counting tool is the following lemma:}
\begin{lemma}[Path counting estimate]
\label{lem:bad_paths}
For any $j \geq 1$, $a>1$,

    \begin{equation}\label{eq:lem_bad_path_nj}
    \begin{split}
        \lim_{m \to \infty} \limsup_{n \to \infty} \frac{1}{n}\int_{I_n}\cdots\int_{I_n}& \frac{1}{n^j} \sum_{v_1,\ldots,v_j \in [n]}\mathbb{W}^{v_1,\ldots,v_j}_{n}(\vec{x})\\&  {\times}{\left(\prod_{i=0}^{j-2}\ind{\|x_{i+1}-x_i\|<a^{m^{i+1}}}\right)}\ind{\|x_{j}-x_{j-1}\|>a^{m^j}} dx_{0}\cdots dx_j =0,
    \end{split}
    \end{equation}
and
    \begin{equation}\label{eq:lem_bad_path_j}
        \begin{split}
            \lim_{m \to \infty} \int_{\mathbb{R}^d}\ldots \int_{\mathbb{R}^d} &\mathbb{W}_j(\mathbf{0},x_1\ldots,x_j)\\&   {\times}\ind{\|x_1\|<a^m}{\left( \prod_{i=1}^{j-2}\ind{\|x_{i+1}-x_i\|<a^{m^{i+1}}}\right)}\ind{\|x_j-x_{j-1}\|>a^{m^j}} dx_1\cdots dx_j=0.
        \end{split} 
    \end{equation}
\end{lemma}

As a corollary to Lemma \ref{lem:bad_paths}, we will show for any fixed $K \in \mathbb{N}$  {and $a,m>1$}, if we choose 
\begin{equation}\label{eq:choice_of_r}
    r= r(a,m,K)= a^m+a^{m^2}+a^{m^3}+\dots+a^{m^K},
\end{equation}
then as $m \to \infty$, with high probability, the $K$-neighborhood of the rooted graph $(F^{\mathbb{G}_n}_{U_n}(r),U_n)$ will be the $K$-neighborhood $B^{\mathbb{G}_n}_{U_n}(K)$ of $(\mathbb{G}_n,U_n)$, 
and a similar result for $(\mathbb{G}_{\infty},0)$. This  {will} help us in proving Theorems \ref{thm:main_LWC_SIRGs} and \ref{thm:main_LCinP_SIRGs}. To state  {the} result, we again introduce some shorthand notations:

\begin{equation}
    \label{eq:imp_notations}
    \begin{split}
    B_n&:=B^{\mathbb{G}_n}_{{U_n}}(K),\quad B:=B^{\mathbb{G}_{\infty}}_{0}(K),\qquad 
    F_{n,r}:=F^{\mathbb{G}_n}_{{U_n}}(r),\quad F_r:=F^{\mathbb{G}_{\infty}}_{0}(r),\\
    BF_{n,r}&:=B^{F_{n,r}}_{{U_n}}(K),\quad BF_{r}:=B^{F_{r}}_{0}(K).
    \end{split}
\end{equation}

\begin{figure}[H]
    \centering
    \resizebox{0.5\textwidth}{!}{
        \begin{tikzpicture} 

\draw (8,4) arc (0:360:4);
\coordinate[rootnode] (O) at (4,4);
\coordinate[neighborhood2] (path1node1) at (8,8);
\coordinate[neighborhood2] (path1node2) at (10,5);
\coordinate[neighborhood2] (path2node1) at (4,1);
\coordinate[neighborhood2] (path2node2) at (1,0);
\coordinate[neighborhood2] (path3node1) at (1,3);
\coordinate[neighborhood2] (path4node1) at (1,8);
\coordinate[neighborhood2] (path5node1) at (5,6);
\coordinate[neighborhood2] (path5node2) at (6,7);
\coordinate[neighborhood2] (path5branchto1) at (3,7);
\coordinate[notneighborhood2] (otherinside1) at (6,5);
\coordinate[notneighborhood2] (otherinside2) at (3,3);
\coordinate[notneighborhood2] (otherinside3) at (7,3);
\coordinate[notneighborhood2] (otherinside4) at (1,5);
\coordinate[notneighborhood2] (otherinside5) at (6,2);
\coordinate[notneighborhood2] (otheroutside1) at (8,0);
\coordinate[notneighborhood2] (otheroutside2) at (0,1);
\coordinate[notneighborhood2] (otheroutside3) at (9,3);
\coordinate[notneighborhood2] (otheroutside4) at (3,9);
\coordinate[notneighborhood2] (otheroutside5) at (-2,5);
\draw (O) -- (path1node1) -- (path1node2);
\draw (O) -- (path2node1) -- (path2node2);
\draw (O) -- (path3node1);
\draw (O) -- (path4node1);
\draw (O) -- (path5node1) -- (path5node2);
\draw (path5node1) -- (path5branchto1) -- (path2node2);

\EncircleNode{path2node1};
\EncircleNode{path3node1};
\EncircleNode{path5node1};
\EncircleNode{path5node2};
\EncircleNode{path5branchto1};
\EncircleNode{otherinside1};
\EncircleNode{otherinside2};
\EncircleNode{otherinside3};
\EncircleNode{otherinside4};
\EncircleNode{otherinside5};

\end{tikzpicture}
    }
    \caption{Illustration to distinguish between the two kinds of neighbourhoods. The star in the middle is the location $Y^{(n)}_{U_n}$ of the root $U_n$, and the big circle around it is the boundary of the Euclidean ball centered at $Y^{(n)}_{U_n}$ of radius $r$, in $\mathbb{R}^d$. Diamonds are the vertices of the graph neighbourhood $B^{\mathbb{G}_n}_{U_n}(2)$ of radius $2$ around the root. Black dots are the vertices which are not in $B^{\mathbb{G}_n}_{U_n}(2)$. The encircled vertices are the vertices of $F^{\mathbb{G}_n}_{U_n}(r)$. {The encircled diamonds are the vertices of $BF_{n,r}$, the graph neighbourhood of radius $2$ about the root $U_n$, of the graph $F^{\mathbb{G}_n}_{U_n}(r)$.}}
    \label{fig:Fig_two_types_nbhds}
\end{figure}
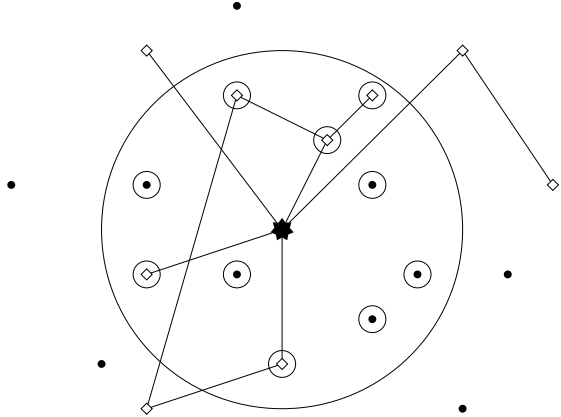

\begin{remark}[Spatial and graph neighborhoods]
At this point, we emphasize that  {we rely on} two kinds of \emph{neighborhoods}  {around} the root $U_n$ (respectively $0$) of $\mathbb{G}_n$ (respectively, $\mathbb{G}_{\infty}$).  {These are} the graph $F_{n,r}$ (respectively $F_r$), which is the subgraph induced by those vertices whose  {{\em spatial locations}} are within  {Euclidean distance $r$} of the root location $Y^{(n)}_{U_n}$ (respectively $Y_0$), and the graph $B_n$ (respectively $B$), which is the {\em graph neighborhood} of radius $K$  {around} the root $U_n$ (respectively $0$) (see Figure \ref{fig:Fig_two_types_nbhds}). The difference between these two kinds of \emph{neighborhoods} is to be understood clearly. For example, the graph $F_{n,r}$ may possibly be disconnected, while $B_n$ is always a connected graph.   Moreover, $BF_{n, r}$ (respectively $BF_{r}$) is the graph neighborhood of radius $K$, of the rooted graph $(F_{n,r},U_n)$ (respectively $(F_r,0)$). 
\end{remark}

We have the following corollary to Lemma \ref{lem:bad_paths}:

\begin{corollary}[{Coupling spatial and graph neighborhoods}]
\label{cor:error_probs_small}
Let $BF_{n,r}$, $B_n$, $BF_r$, $B$ be as in (\ref{eq:imp_notations}), where $r=r(a,m,K)$ is as in (\ref{eq:choice_of_r}). Then,
    \begin{equation}\label{eq:error_n_small_lem}
        \lim_{m \to \infty} \limsup_{n \to \infty} \Prob{BF_{n,r}\neq B_n} =0,
    \end{equation} 
and
    \begin{equation}\label{eq:error_small_lem}
        \lim_{m \to \infty} \Prob{BF_{r}\neq B} =0.
    \end{equation}
\end{corollary}

In proving Corollary \ref{cor:error_probs_small}, we  {perform} a careful path-counting analysis to bound the expected number of $K$-paths of $\mathbb{G}_n$ which are not a $K$-path in $F^{\mathbb{G}_n}_{U_n}(r)$, by the integral expression (\ref{eq:lem_bad_path_nj}), and the similar bound (\ref{eq:lem_bad_path_j}) for $\mathbb{G}_{\infty}$. Corollary \ref{cor:error_probs_small} then follows directly using the technical Lemma \ref{lem:bad_paths}. 

 {In the course of} proving Corollary \ref{cor:error_probs_small}, we develop a general path estimate, where for $r=r(a,m,K)$ as in (\ref{eq:choice_of_r}), we can have $a=a_n$, $K=K_n$ to be $n$ dependent sequences, while $m$  {does not depend} on $n$. This general estimate will be used in the proof of Theorem \ref{thm:main_typical_dist}.


\paragraph{Proof strategy of Theorem \ref{thm:main_LWC_SIRGs}.} In Section \ref{ssec:proof_LWC}, we prove Theorem \ref{thm:main_LWC_SIRGs}. Recall  {that}, to conclude Theorem \ref{thm:main_LWC_SIRGs}, we need to show  {that} for any $K \in \mathbb{N}$, the $K$-neighborhoods of $(\mathbb{G}_n,U_n)$ and $(\mathbb{G}_{\infty},0)$ are close in distribution in $\mathcal{G}_{\star}$. Consequently, to conclude Theorem \ref{thm:main_LWC_SIRGs}, using Corollary \ref{cor:error_probs_small}, it will be enough to show  {that} the $K$-neighborhoods of $(F^{\mathbb{G}_n}_{U_n}(r),U_n)$ and $(F^{\mathbb{G}_{\infty}}_{0}(r),0)$ are close in distribution. This we will observe to be an easy consequence of Proposition \ref{prop:loc_trunc}.


\paragraph{Proof strategy of Theorem \ref{thm:main_LCinP_SIRGs}.} In Section \ref{ssec:proof_LCinP},  we prove Theorem \ref{thm:main_LCinP_SIRGs}. The first step is to show  {that} the empirical Euclidean graph neighborhood distribution concentrates. That is, for $H_*=(H,h) \in \mathcal{G}_{\star}$, where $h \in V(H)$, and for any $r>0$, we define the random variables  

\begin{equation}\label{eq:defn_C_rn_rvs}
    C_{r,n}(H,h):=\CProb{F^{\mathbb{G}_n}_{U_n}(r)\cong (H,h)}{\mathbb{G}_n}=\frac{1}{n}\sum_{i=1}^n \ind{F^{\mathbb{G}_n}_i(r)\cong (H,h)},
\end{equation}
and show these random variables concentrate:
\begin{lemma}[Concentration of empirical Euclidean graph neighborhood measure]\label{Lem:sec_mom_C_rn}
For any $r>0$ and a locally finite rooted graph $H_*= (H,h) \in \mathcal{G}_{\star}$, the variance of the random variable $C_{r,n} (H,h)$ converges to $0$ as $n \to \infty$, i.e. 
\[
\left|\Exp{C_{r,n} (H,h)^2}-\Exp{C_{r,n} (H,h)}^2\right|\to0,
\] 
as $n \to \infty$.
\end{lemma}
The key observation in proving this lemma is that the Euclidean graph neighborhoods  {around} two uniformly chosen vertices of $\mathbb{G}_n$  {are} asymptotically independent, which is a consequence of the fact that the distance between the locations of two uniformly chosen vertices of $\mathbb{G}_n$ in $\mathbb{R}^d$  {diverges} in probability, as $n \to \infty$.

Next, we combine Lemma \ref{Lem:sec_mom_C_rn} with Corollary \ref{cor:error_probs_small} to show that the empirical neighborhood distribution of $\mathbb{G}_n$ also concentrate. That is, for any $K \in \mathbb{N}$ and $G_*=(G,g) \in \mathcal{G}_{\star}$, if we define random variables 
\begin{equation}\label{eq:defn_B_n_rvs}
    B_n(G,g):=\CProb{B^{\mathbb{G}_n}_{U_n}(K)\cong (G,g)}{\mathbb{G}_n}=\frac{1}{n}\sum_{i=1}^n \ind{B^{\mathbb{G}_n}_i(K)\cong (G,g)},
\end{equation}
for $n \in \mathbb{N}$, then these random variables also concentrate. This is achieved by first using Lemma \ref{Lem:sec_mom_C_rn} and taking a sum over all rooted graphs $H_*= (H,h)$ with $B^{H}_h(K) \cong (G,g)$, to show that the random variables
\[
\frac{1}{n}\sum_{i =1}^n\ind{ B^{F^{\mathbb{G}_n}_i(r)}_i(K)\cong (G,g)}
\]
concentrate, and then employ Corollary \ref{cor:error_probs_small} to obtain the same conclusion for the variables $B_n(G,g)$.
Theorem \ref{thm:main_LCinP_SIRGs} is then a direct consequence of these observations, combined with Proposition \ref{prop:loc_trunc}.

\subsection{Proof of Proposition \ref{prop:loc_trunc}}\label{ssec:Proof_prop_loc_trunc}
Let us first discuss the proof strategy informally. We make use of the following two key observations:
\begin{itemize}
    \item The number of vertices in $F^{\mathbb{G}_n}_{U_n}(r)$ converges in distribution to the number of vertices in $F^{\mathbb{G}_{\infty}}_0(r)$.
    \item Conditionally on the number of vertices of $F^{\mathbb{G}_n}_{U_n}(r)$ (respectively $F^{\mathbb{G}_{\infty}}_0(r)$) other than the root, their locations are uniform on $\mathscr{B}^r_{Y^{(n)}_{U_n}}$ (respectively $\mathscr{B}^r_{\mathbf{0}}$). It will then follow that any \emph{`not so bad'} translation-invariant function evaluated at the locations of the vertices of $F^{\mathbb{G}_n}_{U_n}(r)$, should have nice limiting behaviour. 
\end{itemize}
    In particular, we define functions $\mathcal{F}^{n}_{ (H,h)}$ and $\mathcal{F}^{\infty}_{ (H,h)}$ (see (\ref{eq:defn_fnG})), which for a given rooted graph $ (H,h) \in \mathcal{G}_{\star}$,  {count} the number of rooted isomorphisms between $ (H,h)$ and $F^{\mathbb{G}_n}_{U_n}(r)$, and between $ (H,h)$ and $F^{\mathbb{G}_{\infty}}_0(r)$.  {We show} that for any such rooted graph $ (H,h) \in \mathcal{G}_{\star}$ , the probability of the event $\{\mathcal{F}^{n}_{ (H,h)}>0\}$ converges to  {that} of the event $\{\mathcal{F}^{\infty}_{ (H,h)}>0\}$.  {In turn, this} implies that the random rooted graph $F^{\mathbb{G}_n}_{U_n}(r)$ converges in distribution to the random rooted graph $F^{\mathbb{G}_{\infty}}_{0}(r)$, on the space $\mathcal{G}_{\star}$. We now go into the details.

\begin{proof}[Proof of Proposition \ref{prop:loc_trunc}]
 {The} random variable $Y^{(n)}_{U_n}$ is uniformly distributed on $I_n$. We write
\begin{align*}
    &\Prob{(F^{\mathbb{G}_n}_{{U_n}}(r),U_n) \cong  (H,h)} \\
      &\qquad\quad= \Prob{(F^{\mathbb{G}_n}_{{U_n}}(r),U_n) \cong  (H,h), Y^{(n)}_{U_n} \in A^r_n}
      +\Prob{(F^{\mathbb{G}_n}_{{U_n}}(r),U_n) \cong  (H,h), Y^{(n)}_{U_n} \notin A^r_n},
\end{align*}
 and observe that, as $n \to \infty$, 
 \begin{align*}
     \Prob{(F^{\mathbb{G}_n}_{{U_n}}(r),U_n) \cong  (H,h), Y^{(n)}_{U_n} \notin A^r_n}
     & \leq \Prob{Y^{(n)}_{U_n} \notin A^r_n} 
        \leq \frac{2d\left(r\left(n^{(d-1)/d}\right)\right)}{n}\to0. \numberthis \label{eq:aas_typical_in_Arn}
 \end{align*}  
 Therefore, it is enough to show that, as $n \to \infty$,
\begin{equation}\label{eq:suff_1}
  \Prob{(F^{\mathbb{G}_n}_{{U_n}}(r),U_n) \cong  (H,h), Y^{(n)}_{U_n} \in A^r_n}\to \Prob{(F^{\mathbb{G}_{\infty}}_{0}(r),0) \cong  (H,h)}.
\end{equation}

Note that
\begin{align*}
  &\Prob{(F^{\mathbb{G}_n}_{{U_n}}(r),U_n) \cong  (H,h), Y^{(n)}_{U_n} \in A^r_n}
  \\& \qquad\quad =\Prob{(F^{\mathbb{G}_n}_{{U_n}}(r),U_n) \cong  (H,h), Y^{(n)}_{U_n} \in A^r_n, |V(F^{\mathbb{G}_n}_{U_n}(r))|=|V(H)|},
\end{align*}
and so we can repeatedly condition to rewrite
\begin{equation}\label{eq:F_prop_2}
    \begin{split}
    &\Prob{(F^{\mathbb{G}_n}_{{U_n}}(r),U_n) \cong  (H,h), Y^{(n)}_{U_n} \in A^r_n}
    \\& \qquad= \CProb{(F^{\mathbb{G}_n}_{{U_n}}(r),U_n) \cong  (H,h)}{|V(F^{\mathbb{G}_n}_{{U_n}}(r))|=|V(H)|, Y^{(n)}_{U_n} \in A^r_n}
    \\& \qquad\quad \times \CProb{|V(F^{\mathbb{G}_n}_{{U_n}}(r))|=|V(H)|}{Y^{(n)}_{U_n} \in A^r_n}\Prob{Y^{(n)}_{U_n} \in A^r_n}.
    \end{split}
\end{equation}

Using (\ref{eq:aas_typical_in_Arn}) the last term in the RHS of (\ref{eq:F_prop_2})  {tends} to $1$ as $n \to \infty$.  {Observe} that         
    \[\CProb{|V(F^{\mathbb{G}_n}_{{U_n}}(r))|
    =|V(H)|}{Y^{(n)}_{U_n} \in A^r_n}=\Prob{\mathcal{Y}_n=|V(H)|-1},
    \] 
where $\mathcal{Y}_n$ follows a $\text{Bin}(n-1, \lambda_d(\mathscr{B}^r_{\mathbf{0}})/n)$ distribution.

Since $\mathcal{Y}_n$ converges in distribution to $\mathcal{Y}\sim \text{Poi}(\lambda_d(\mathscr{B}^r_{\mathbf{0}}))$, and since $\mathcal{Y}$ is equal in distribution to $\Gamma(\mathscr{B}^r_{\mathbf{0}})$ (recall $\Gamma$ from (\ref{eq:defn_Gamma_infty})), it follows that $\mathcal{Y}_n \dlim \Gamma(\mathscr{B}^r_{\mathbf{0}})$, as $n \to \infty$.

Observe that $\Gamma(\mathscr{B}^r_{\mathbf{0}})\stackrel{d}{=}|V(F^{\mathbb{G}_{\infty}}_{\mathbf{0}}(r))|-1$, so that
\begin{align*}
  \lim_{n \to \infty} \CProb{|V(F^{\mathbb{G}_n}_{{U_n}}(r))|=|V(H)|}{Y^{(n)}_{U_n} \in A^r_n}
  &=\lim_{n \to \infty} \Prob{\mathcal{Y}_n=|V(H)|-1} 
  \\&=  \Prob{|V(F^{\mathbb{G}_{\infty}}_{0}(r))|-1=|V(H)|-1}
  \\&=\Prob{|V(F^{\mathbb{G}_{\infty}}_{0}(r))|=|V(H)|}.
\end{align*}
Hence from (\ref{eq:F_prop_2}), we note that to conclude (\ref{eq:suff_1}), it is enough to show that
    \begin{equation}\label{eq:F_prop_3}
    \begin{split}
    \lim_{n \to \infty}\CProb{(F^{\mathbb{G}_n}_{{U_n}}(r),U_n) \cong  (H,h)}{|V(F^{\mathbb{G}_n}_{{U_n}}(r))|=|V(H)|, Y^{(n)}_{U_n} \in A^r_n}  \\&\hspace{- 200 pt}=\CProb{(F^{\mathbb{G}_{\infty}}_{0}(r),0) \cong  (H,h)}{|V(F^{\mathbb{G}_{\infty}}_{0}(r))|=|V(H)|}.
    \end{split}
    \end{equation}

For the remainder of the proof, we assume $|V(H)|=l+1$. We continue by making some observations on the locations and weights of the vertices of the graph $F^{\mathbb{G}_n}_{U_n}(r)$ (respectively $F^{\mathbb{G}_\infty}_{\mathbf{0}}(r)$), conditionally on $\{|V(F^{\mathbb{G}_n}_{U_n}(r))|=l+1, Y^{(n)}_{U_n} \in A^r_n\}$ (respectively \{$|V(F^{\mathbb{G}_\infty}_{\mathbf{0}}(r))|=l+1$\}) .

\paragraph{Locations of $F^{\mathbb{G}_n}_{U_n}(r)$.}
Since $(Y^{(n)}_i)_{i\in[n]}$ is an i.i.d.\ collection of uniform random variables on $I_n$, conditionally on the event $\{Y^{(n)}_{U_n}\in A^r_n, |V(F^{\mathbb{G}_n}_{{U_n}}(r))|=l+1\}$, the locations $P_1,\ldots,P_l$ of the $l$ vertices of $\mathbb{G}_n$ (in some order) falling in $\mathscr{B}^r_{Y^{(n)}_{U_n}}$, other than $Y^{(n)}_{U_n}$, are at independent uniform locations in the ball $\mathscr{B}^r_{Y^{(n)}_{U_n}}$, given $Y^{(n)}_{U_n}$. Hence, conditionally on $\{Y^{(n)}_{U_n}\in A^r_n, |V(F^{\mathbb{G}_n}_{{U_n}}(r))|=l+1\}$, the random variables $P_1-Y^{(n)}_{U_n},\ldots,P_l-Y^{(n)}_{U_n}$ are independently, uniformly distributed on the ball $\mathscr{B}^r_{\mathbf{0}}$.   

Conditionally on $\{Y^{(n)}_{U_n}\in A^r_n, |V(F^{\mathbb{G}_n}_{{U_n}}(r))|=l+1\}$, let the locations of all the vertices of $\mathbb{G}_n$ falling in $\mathscr{B}^r_{Y^{(n)}_{U_n}}$ (including $Y^{(n)}_{U_n}$) be $P_0, P_1, \ldots, P_l$, where $Y^{(n)}_{U_n}=P_0$. So, conditionally on $\{Y^{(n)}_{U_n}\in A^r_n, |V(F^{\mathbb{G}_n}_{{U_n}}(r))|=l+1\}$, the random matrix $(\|P_i-P_j\|)_{0\leq i,j \leq l; i\neq j}$ is equal in distribution to $(\|Y_i-Y_j\|)_{0\leq i,j \leq l; i\neq j}$, where the set $\{Y_i\colon 1\leq i\leq l\}$ consists of $l$ i.i.d.\ uniform points in $\mathscr{B}^r_{\mathbf{0}}$ (independent of $\mathbf{Y}^{(n)}$, $\mathbf{W}^{(n)}$, $\mathbf{Y}$, $\mathbf{W}$) and $Y_0=\mathbf{0}$.

\paragraph{Locations of $F^{\mathbb{G}_{\infty}}_0(r)$.}
Again, conditionally on the event $\{\Gamma_{\infty}(\mathscr{B}^r_{\mathbf{0}})=l+1\}=\{|V(F^{\mathbb{G}_{\infty}}_{0}(r))|=l+1\}$, the locations of the $l$ vertices of $\mathbb{G}_{\infty}$ in $\mathscr{B}^r_{\mathbf{0}}$ other than $0$ are i.i.d.\ uniform on $\mathscr{B}^r_{\mathbf{0}}$ (since $\Gamma$ is a homogeneous Poisson point process). So if $Z_0,\ldots,Z_l$ are the locations of the $l+1$ vertices of $\mathbb{G}_{\infty}$ in $\mathscr{B}^r_{\mathbf{0}}$ (in some order) where $Z_0=\mathbf{0}$, the random matrix $(\|Z_i-Z_j\|)_{0\leq i,j \leq l; i\neq j}$ is also equal in distribution to $(\|Y_i-Y_j\|)_{0\leq i,j \leq l; i\neq j}$, where the set $\{Y_i\colon 1\leq i\leq l\}$ consists of $l$ i.i.d.\ uniform points in $\mathscr{B}^r_{\mathbf{0}}$ (independent of $\mathbf{Y}^{(n)}$, $\mathbf{W}^{(n)}$, $\mathbf{Y}$, $\mathbf{W}$) and $Y_0=\mathbf{0}$.

We conclude that
\begin{equation}\label{eq:F_prop_4}
    \begin{split}
        (\|P_i-P_j\|)_{0\leq i,j \leq l; i\neq j}\bigg|{\{Y^{(n)}_{U_n}\in A^r_n, |V(F^{\mathbb{G}_n}_{{U_n}}(r))|=l+1\}}
       \\& \hspace{- 200 pt} \stackrel{d}{=}(\|Y_i-Y_j\|)_{1\leq i,j \leq l+1; i\neq j}
       \\& \hspace{- 200 pt} \stackrel{d}{=} (\|Z_i-Z_j\|)_{0\leq i,j \leq l; i\neq j}\bigg|{\{\Gamma_{\infty}(\mathscr{B}^r_{\mathbf{0}})=l+1\}}. 
    \end{split}
\end{equation}

\paragraph{Weights of $F^{\mathbb{G}_n}_{U_n}(r)$.}
Next, conditionally on $\{Y^{(n)}_{U_n}\in A^r_n, |V(F^{\mathbb{G}_n}_{{U_n}}(r))|=l+1\}$,  {for $0 \leq i \leq l$}, let $W_{n,i}$ denote the weight of the vertex of $F^{\mathbb{G}_n}_{U_n}(r)$ with location $P_i$. Then, $(W_{n,0},\ldots,W_{n,l})$ has the distribution of $l+1$ uniformly chosen weights from the weight set $\{W^{(n)}_1,\ldots,W^{(n)}_n\}$ without replacement, in some arbitrary order. This is because for any $i_0,\ldots,i_l \in [n]$,
\begin{align*}
    & \CProb{(W_{n,0},\ldots,W_{n,l})=(W^{(n)}_{i_0},\ldots,W^{(n)}_{i_l})}{Y^{(n)}_{U_n}\in A^r_n, |V(F^{\mathbb{G}_n}_{{U_n}}(r))|=l+1}\\
    & = \frac{\CProb{(W_{n,0},\ldots,W_{n,l})=(W^{(n)}_{i_0},\ldots,W^{(n)}_{i_l})}{Y^{(n)}_{U_n}\in A^r_n}}{\CProb{|V(F^{\mathbb{G}_n}_{{U_n}}(r))|=l+1}{Y^{(n)}_{U_n}\in A^r_n}},
     \numberthis \label{eq:weight_UAR_WR_1}
\end{align*}
where in the second step we use the fact that the event that the vector of weights of the vertices of $F^{\mathbb{G}_n}_{U_n}(r)$ is of length $l+1$, is contained in the event that $\{|V(F^{\mathbb{G}_n}_{{U_n}}(r))|=l+1\}$.
The numerator in (\ref{eq:weight_UAR_WR_1}) is 
\[\frac{1}{n} \times \CProb{Y^{(n)}_{i_p}\in \mathscr{B}^{r}_{Y^{(n)}_{i_0}}\;\forall\;p \in [l],\; Y^{(n)}_{i_q}\notin \mathscr{B}^{r}_{Y^{(n)}_{i_0}}\;\forall\; i_q \in [n]\setminus\{i_0,i_1,\ldots,i_l\}}{Y^{(n)}_{i_0}\in A^r_n}\times \frac{1}{l!},\]
(where the term $ {1/n}$ is just the probability that $U_n = i_0$, and the term $\frac{1}{l!}$ accounts for the choice of the ordering $P_k=Y^{(n)}_{i_k}$ for $1 \leq k \leq l$ among all possible $l!$ labelings).  {This} evaluates to 
    \[
    \frac{1}{n} \times \left(\frac{\lambda_d(\mathscr{B}^r_{\mathbf{0}})}{n}\right)^{l}\times \left(1-\frac{\lambda_d(\mathscr{B}^r_{\mathbf{0}})}{n}\right)^{n-1-l} \times \frac{1}{l!}.
    \]
The denominator in (\ref{eq:weight_UAR_WR_1}) is just the probability that a $\text{Bin}\left(n-1, \frac{\lambda_d(\mathscr{B}^r_{\mathbf{0}})}{n}\right)$ random variable takes the value $l$, which evaluates to $\binom{n-1}{l}\times \left(\frac{\lambda_d(\mathscr{B}^r_{\mathbf{0}})}{n}\right)^{l}\times \left(1-\frac{\lambda_d(\mathscr{B}^r_{\mathbf{0}})}{n}\right)^{n-1-l}$. Combining, we get that
\begin{equation}\label{eq:wt_UAR_WR_2}
  \CProb{(W_{n,0},\ldots,W_{n,l})=(W^{(n)}_{i_0},\ldots,W^{(n)}_{i_l})}{Y^{(n)}_{U_n}\in A^r_n, |V(F^{\mathbb{G}_n}_{{U_n}}(r))|=l+1}=\frac{1}{n(n-1)\cdots(n-l)}.  
\end{equation}

We also make the observation that the random vector $(W_{n,0},\ldots,W_{n,l})$ is permutation invariant: for any $\sigma \in S(l+1)$, where $S(l+1)$ denotes the set of all permutations of $\{0,\ldots,l\}$,
\begin{equation}\label{eq:weight_FGnUnr_permutation_inv}
 (W_{n,0},\ldots,W_{n,l}) \stackrel{d}{=} (W_{n,\sigma(0)},\ldots,W_{n,\sigma(l)}).   
\end{equation}

\paragraph{Weights of $F^{\mathbb{G}_{\infty}}_0(r)$.}
Again, conditionally on $\{\Gamma_{\infty}(\mathscr{B}^r_{\mathbf{0}})=l+1\}$,  {for $0 \leq i \leq l$}, let $W_{\infty,i}$ denote the weight of the vertex of $F^{\mathbb{G}_{\infty}}_{0}(r)$ with location $Z_i$. Since the weights in the limiting graph $\mathbb{G}_{\infty}$ are i.i.d., it immediately follows that the random weight vector $(W_{\infty,0},\ldots,W_{\infty,l})$ is permutation invariant: for any $\sigma \in S(l+1)$, 
\begin{equation}\label{eq:weight_F_infty_perm_inv}
    (W_{\infty,0},\ldots,W_{\infty,l})\stackrel{d}{=}(W_{\infty,\sigma(0)},\ldots,W_{\infty,\sigma(l)}),
\end{equation}
and that given $\{\Gamma_{\infty}(\mathscr{B}^r_{\mathbf{0}})=l+1\}$, the vector $(W_{\infty,0},\ldots,W_{\infty,l})\stackrel{d}{=}(\mathrm{W}_0,\ldots,\mathrm{W}_l)$, where each entry of the vector $(\mathrm{W}_0,\ldots,\mathrm{W}_l)$ is an i.i.d.\ copy of the limiting weight variable $W$ (recall (\ref{eq:conv_in_law_Wn_to_W})).

 {We now continue with the proof.}

\paragraph{Convergence of weights of $F^{\mathbb{G}_n}_{U_n}(r)$ to the weights of $F^{\mathbb{G}_{\infty}}_0(r)$.}
We also note that
\begin{equation}\label{eq:wt_FGn_conv_wt_FGinfty}
 (W_{n,0},\ldots,W_{n,l}) \bigg| \{Y^{(n)}_{U_n}\in A^r_n, |V(F^{\mathbb{G}_n}_{{U_n}}(r))|=l+1\} \dlim (W_{\infty,0},\ldots,W_{\infty,l}) \bigg| \{\Gamma_{\infty}(\mathscr{B}^r_{\mathbf{0}})=l+1\}.  
\end{equation}
This is because for any continuity set $A_0 \times\dots\times A_l$ of $(\mathrm{W}_0,\ldots,\mathrm{W}_l)$, by (\ref{eq:wt_UAR_WR_2}),
\begin{align*}
    & \CProb{(W_{n,0},\ldots,W_{n,l})\in A_0 \times\dots\times A_l}{Y^{(n)}_{U_n}\in A^r_n, |V(F^{\mathbb{G}_n}_{{U_n}}(r))|=l+1} \\
    & \geq \Exp{\sum_{j_0,\ldots,j_l}\frac{1}{n^l}\prod_{k=0}^l{\ind{W^{(n)}_{j_k}\in A_k}}}= \Exp{\prod_{k=0}^l \left(\frac{1}{n}\sum_{i=1}^n {\ind{W^{(n)}_i \in A_k}} \right)},
\end{align*}
and 
\begin{align*}
    & \CProb{(W_{n,0},\ldots,W_{n,l})\in A_0 \times\cdots\times A_l}{Y^{(n)}_{U_n}\in A^r_n, |V(F^{\mathbb{G}_n}_{{U_n}}(r))|=l+1} \\
    & \leq \Exp{\sum_{j_0,\ldots,j_l}\frac{1}{(n-l)^l}\prod_{k=0}^l{\ind{W^{(n)}_{j_k}\in A_k}}} =\Exp{\prod_{k=0}^l \left(\frac{1}{n-l}\sum_{i=1}^n {\ind{W^{(n)}_i \in A_k}} \right)},
\end{align*}
where the sums in the last two displays are taken over all cardinality $l+1$ subsets $\{j_0,\ldots,j_l\}$ of $[n]$.
Since the RHS of the last two displays are  {bounded} from above by $1$, we use dominated convergence and apply (\ref{eq:conv_in_law_Wn_to_W}), to conclude that the RHS in the last two displays both converge to
\begin{align*}
  \prod_{k=0}^l \Prob{W \in A_k}
    & = \Prob{(\mathrm{W}_0,\ldots,\mathrm{W}_l)\in A_0\times \cdots \times A_l}\\
    & = \CProb{(W_{\infty,0},\ldots,W_{\infty,l})\in A_0\times \cdots \times A_l}{\Gamma_{\infty}(\mathscr{B}^r_{\mathbf{0}})=l+1}.  
\end{align*}

\paragraph{Functions counting rooted isomorphisms, and their symmetry properties.} 
Now we proceed to define the functions we use to count the number of isomorphisms between the given rooted graph $ (H,h)$, and the Euclidean graph neighborhoods $F^{\mathbb{G}_n}_{U_n}(r)$, $F^{\mathbb{G}_{\infty}}_0(r)$. Let the vertices of $ (H,h)$ be $v_0,\ldots,v_l$, where $h=v_0$. Let us denote the subset of the set of permutations $S(l+1)$ of $\{0,1,\ldots,l\}$ that fix $0$ by $S_0(l+1)$, i.e., for all $\sigma \in S_0(l+1)$, $\sigma(0)=0$.

Let $M^S_{l+1}(\mathbb{R})$ denote the space of all square symmetric matrices of order $l+1$ with entries in $\mathbb{R}$. For each $n \in \mathbb{N}$, define the function $\mathcal{F}^n_{ (H,h)}\colon \left(\mathbb{R}^d\right)^{l+1}\times \left(\mathbb{R}^d\right)^{l+1} \times M^S_{l+1}(\mathbb{R}) \to \mathbb{R}$ as

\begin{equation}\label{eq:defn_fnG}
    \begin{split}
    \mathcal{F}^n_{ (H,h)}\left(\vec{x},\vec{y},(a_{ij})_{i,j=0}^l\right)\\
    :=\sum_{\pi \in S_0(l+1)} 
    & \prod_{\{v_i,v_j\}\in E(H)} \left({\ind{a_{\pi(i)\pi(j)}<\kappa_n\left(\|x_{\pi(i)}-x_{\pi(j)}\|,y_{\pi(i)},y_{\pi(j)}\right)}}\right)
    \\& \times \prod_{\{v_i,v_j\}\notin E(H)}\left({\ind{a_{\pi(i)\pi(j)}>\kappa_n\left(\|x_{\pi(i)}-x_{\pi(j)}\|,y_{\pi(i)},y_{\pi(j)}\right)}} \right),    
    \end{split}
\end{equation}
for $\vec{x}=(x_0,\dots,x_l), \vec{y}=(y_0,\dots,y_l) \in (\mathbb{R}^d)^{l+1}$, and $(a_{ij})_{i,j=0}^l \in M^S_{l+1}(\mathbb{R})$, and similarly define the function $\mathcal{F}^{\infty}_{ (H,h)}\colon \left(\mathbb{R}^d\right)^{l+1}\times \left(\mathbb{R}^d\right)^{l+1} \times M^S_{l+1}(\mathbb{R})\to \mathbb{R}$, where $\mathcal{F}^{\infty}_{ (H,h)}$ is just $\mathcal{F}^n_{ (H,h)}$ with $\kappa_n$ replaced by $\kappa$. 

{Heuristically, we want the indicators $\ind{a_{\pi(i)\pi(j)}<\kappa_n\left(\|x_{\pi(i)}-x_{\pi(j)}\|,y_{\pi(i)},y_{\pi(j)}\right)}$ to be the indicators of the events $\{\text{the edge}\; \{\pi(i),\pi(j)\}\; \text{is present}\}$. Since in our graphs these events occur independently each with probability $\kappa_n\left(\|x_{\pi(i)}-x_{\pi(j)}\|,y_{\pi(i)},y_{\pi(j)}\right)$ when the locations and weights of the vertices $\pi(i)$ and $\pi(j)$ are respectively $(x_{\pi(i)},y_{\pi(i)})$ and $(x_{\pi(j)},y_{\pi(j)})$, we will take the matrix $(a_{ij})_{i,j,=0}^l$ to be a symmetric i.i.d.\ uniform matrix. Before that, we first discuss some symmetry properties of the functions $\mathcal{F}^n_{ (H,h)}$ and $\mathcal{F}^{\infty}_{ (H,h)}$.}

Observe the following symmetry: for any permutation $\pi \in S_0(l+1)$, and for $ {\bullet}$ being either $n$ or $\infty$,
\begin{align*}
 & \mathcal{F}^{ {\bullet}}_{ (H,h)}\left((x_0,\ldots,x_l),(y_0,\ldots,y_l),(a_{ij})_{i,j=0}^l\right)   
  \\ &
 = \mathcal{F}^{ {\bullet}}_{ (H,h)}\left((x_{\pi(0)},\ldots,x_{\pi(l)}),(y_{\pi(0)},\ldots,y_{\pi(l)}),(a_{\pi(i)\pi(j)})_{i,j=0}^l\right). 
\end{align*}

Recall the permutation invariance of the weights of $F^{\mathbb{G}_n}_{U_n}(r)$ and $F^{\mathbb{G}_{\infty}}_{0}(r)$ from (\ref{eq:weight_FGnUnr_permutation_inv}) and (\ref{eq:weight_F_infty_perm_inv}), and note that for $(x_0,\ldots,x_l)$, $(z_0,\ldots,z_l)\in \left(\mathbb{R}^d\right)^{l+1}$, if there exists a permutation $\pi \in S_0(l+1)$ such that

    \[
    (\|x_i-x_j\|)_{i,j=0}^l=(\|z_{\pi(i)}-z_{\pi(j)}\|)_{i,j=0}^l\;\; \text{(entry-wise)},
    \]
then for any $(a_{ij})_{i,j=0}^l\in M^S_{l+1}(\mathbb{R})$, conditionally on $\{Y^{(n)}_{U_n}\in A^r_n, |V(F^{\mathbb{G}_n}_{{U_n}}(r))|=l+1\}$,       
    \begin{equation}\label{eq:symm_fnG}
    \begin{split}
        & \mathcal{F}^n_{ (H,h)}\left((x_0,\ldots,x_l),(W_{n,0},\ldots,W_{n,l}),(a_{ij})_{i,j=0}^l\right)
        \\& \stackrel{d}{=}\mathcal{F}^n_{ (H,h)}\left((z_0,\ldots,z_l),(W_{n,0},\ldots,W_{n,l}),(a_{ij})_{i,j=0}^l)\right),    
    \end{split}
    \end{equation}
and a similar distributional equality for $\mathcal{F}^{\infty}_{ (H,h)}$, conditionally on $\{\Gamma_{\infty}(\mathscr{B}^r_{\mathbf{0}})=l+1\}$.

\paragraph{Simplifying the events $\{F^{\mathbb{G}_n}_{{U_n}}(r)\cong  (H,h)\}$ and $\{F^{\mathbb{G}_{\infty}}_{0}(r) \cong  (H,h)\}$.}
Now, let $(U_{ij})_{i,j=0}^l$ and $(U'_{ij})_{i,j=0}^l$ be two i.i.d.\ random elements of $M^S_{l+1}(\mathbb{R})$, which are also independent from $\mathbf{X}^{(n)}$, $\mathbf{W}^{(n)}$, $\mathbf{X}$, $\mathbf{W}$ and $\{Y_i:0\leq i\leq l\}$, where for each $i<j$, $U_{ij}\sim \text{U}([0,1])$, all the entries above the diagonal of the random matrix $(U_{ij})_{i,j=1}^{l+1}$ are independent, and for each $i$, $U_{ii}:=0$ (the diagonal elements will not come into picture and can be defined arbitrarily).

Let 
\begin{equation}\label{eq:F_prop_RHS}
    \mathcal{A}_n
    :=
    {\ind{F^{\mathbb{G}_n}_{{U_n}}(r)\cong  (H,h)}}\bigg|{\{Y^{(n)}_{U_n} \in A^r_n, |V(F^{\mathbb{G}_n}_{{U_n}}(r))|=l+1\}},
\end{equation}
and 
\begin{equation}\label{eq:F_prop_LHS}
    \mathcal{A}_{\infty}:=
    {\ind{F^{\mathbb{G}_{\infty}}_{0}(r) \cong  (H,h)}}\bigg|{\{|V(F^{\mathbb{G}_{\infty}}_{0}(r))|=l+1\}}.    
\end{equation}
Note from (\ref{eq:F_prop_3}) that our target is to show $\Exp{\mathcal{A}_n}\to \Exp{\mathcal{A}_{\infty}}$, as $n \to \infty$.

Recall (\ref{eq:F_prop_4}).
Observe that
\begin{equation}\label{eq:F_prop_5}
   \mathcal{A}_n
   \stackrel{d}{=}
   {\ind{\mathcal{F}^n_{(H,h)}\left((P_0,\ldots,P_l), (W_{n,0},\ldots,W_{n,l}), (U_{ij})_{i,j=0}^l\right)>0}}\bigg|{\{Y^{(n)}_{U_n} \in A^r_n, |V(F^{\mathbb{G}_n}_{{U_n}}(r))|=l+1\}}, 
\end{equation}
since, conditionally on the event $\{Y^{(n)}_{U_n} \in A^r_n, |V(F^{\mathbb{G}_n}_{{U_n}}(r))|=l+1\}$, if for some $\pi \in S_0(l+1)$, the corresponding term in the sum in $\mathcal{F}^n_{ (H,h)}$ is positive, then $v_i \mapsto P_{\pi(i)}$ gives a rooted isomorphism (note that the fact that every permutation in $S_0(l+1)$ fixes $0$ ensures that the isomorphism is rooted), and if there is a rooted isomorphism $\phi\colon G\to F^{G_n}_{U_n}(r)$, then the term in the sum corresponding to the permutation $\sigma \in S_0(l+1)$ is positive, where $\phi(v_i)=P_{\sigma(i)}$.

Using the symmetry (\ref{eq:symm_fnG}) with (\ref{eq:F_prop_4}), we obtain 
\begin{equation}\label{eq:F_prop_6}
   \mathcal{A}_n 
   \stackrel{d}{=} 
   {\ind{\mathcal{F}^n_{ (H,h)}\left((Y_0,\ldots,Y_l), (W_{n,0},\ldots,W_{n,l}), (U_{ij})_{i,j=0}^l\right)>0}}\bigg|{\{Y^{(n)}_{U_n} \in A^r_n, |V(F^{\mathbb{G}_n}_{{U_n}}(r))|=l+1\}}=:\mathscr{A}_n.
\end{equation}

Using exactly similar arguments for the random variable $\mathcal{A}_{\infty}$ and the function $\mathcal{F}^{\infty}_{ (H,h)}$, we obtain
\begin{equation}\label{eq:F_prop_7}
    \mathcal{A}_{\infty} \stackrel{d}{=} 
    {\ind{\mathcal{F}^{\infty}_{ (H,h)}\left((Y_0,\ldots,Y_l), (W_{\infty,0},\ldots,W_{\infty,l}), (U'_{ij})_{i,j=0}^l\right)>0}}\bigg|{\{|V(F^{\mathbb{G}_{\infty}}_{0}(r))|=l+1\}}.
\end{equation}

Finally, using that $(U'_{ij})_{i,j=0}^l\stackrel{d}{=}(U_{ij})_{i,j=0}^l$, and that both the matrices $U_{ij}$ and $U'_{ij}$ are independent of everything else, it is easy to see that 
\begin{equation}\label{eq:F_prop_8}
    \mathcal{A}_{\infty} \stackrel{d}{=} 
    {\ind{\mathcal{F}^{\infty}_{ (H,h)}\left((Y_0,\ldots,Y_l), (W_{\infty,0},\ldots,W_{\infty,l}), (U_{ij})_{i,j=0}^l\right)>0}}\bigg|{\{|V(F^{\mathbb{G}_{\infty}}_{0}(r))|=l+1\}}=:\mathscr{A}_{\infty}.
\end{equation}

\paragraph{Conclusion.}
Using (\ref{eq:wt_FGn_conv_wt_FGinfty}), and using the convergence of connection functions (\ref{eq:pw_conv_strong_cty_kappa_n}), for any $(x_0,\ldots,x_l) \in \left( \mathbb{R}^d \right)^{l+1}$ {such that the collection of positive reals $\{\|x_i-x_j\|:1\leq i,j, \leq l,i\neq j\}$ avoid some set of measure zero, and for any} $(a_{ij})_{i,j=0}^l \in  M^S_{l+1}([0,1])$, as $n\rightarrow \infty$, 
\begin{equation}\label{eq:conv_law_fnG_to_fG}
\begin{split}
    & \mathcal{F}^n_{ (H,h)}\left((x_0,\ldots,x_l),(W_{n,0},\ldots,W_{n,l}),(a_{ij})_{i,j=0}^l \right)\bigg|{\{Y^{(n)}_{U_n} \in A^r_n, |V(F^{\mathbb{G}_n}_{{U_n}}(r))|=l+1\}}
    \\ & \dlim 
    \mathcal{F}^{\infty}_{ (H,h)}\left((x_0,\ldots,x_l), (W_{\infty,0},\ldots,W_{\infty,l}), (a_{ij})_{i,j=0}^l\right)\bigg|{\{|V(F^{\mathbb{G}_{\infty}}_{0}(r))|=l+1\}},
\end{split}
\end{equation} 
which implies that, as $n\rightarrow \infty$, 
    \begin{equation}\label{eq:fnG_positive_to_fG_positive}
    \begin{split}
    & \CProb{\mathcal{F}^n_{ (H,h)}\left((x_0,\ldots,x_l),(W_{n,0},\ldots,W_{n,l}),(a_{ij})_{i,j=0}^l \right)>0}{\{Y^{(n)}_{U_n} \in A^r_n, |V(F^{\mathbb{G}_n}_{{U_n}}(r))|=l+1\}}
\\ &    \to 
    \CProb{\mathcal{F}^{\infty}_{ (H,h)}\left((x_0,\ldots,x_l), (W_{\infty,0},\ldots,W_{\infty,l}), (a_{ij})_{i,j=0}^l\right)>0}{\{|V(F^{\mathbb{G}_{\infty}}_{0}(r))|=l+1\}}.
    \end{split}
    \end{equation}

Combining (\ref{eq:fnG_positive_to_fG_positive}) with the fact that both $(Y_0,\ldots, Y_l)$ and $(U_{ij})_{i,j=0}^l$ are independent of $\mathbf{W}^{(n)}$ and $\mathbf{W}$, {and the fact that the random variables $\|Y_i-Y_j\|$ for $1\leq i,j \leq l$ with $i \neq j$ are continuous random variables and hence almost surely avoid sets of measure $0$, yield}  (recall $\mathscr{A}_n$ and $\mathscr{A}_{ {\infty}}$ from (\ref{eq:F_prop_6}) and (\ref{eq:F_prop_8}))
\[\CExp{\mathscr{A}_n}{(Y_0,\ldots, Y_l),(U_{ij})_{i,j=0}^l}\stackrel{\text{a.s.}}{\to} \CExp{\mathscr{A}_{ {\infty}}}{(Y_0,\ldots, Y_l),(U_{ij})_{i,j=0}^l},\] 
which implies, using dominated convergence, (note {that} domination by $1$ works)
\[\Exp{\mathscr{A}_n} \to \Exp{\mathscr{A}_{ {\infty}}},\] as $n \to \infty$.

Hence, again using (\ref{eq:F_prop_6}) and (\ref{eq:F_prop_8}), $\Exp{\mathcal{A}_n} \to \Exp{\mathcal{A}_{ {\infty}}}$, which is just (\ref{eq:F_prop_3}) by recalling the definition of $\mathcal{A}_n$ from (\ref{eq:F_prop_RHS}) and the definition of $\mathcal{A}_{ {\infty}}$ from (\ref{eq:F_prop_LHS}). This completes the proof of Proposition \ref{prop:loc_trunc}.
\end{proof}

\subsection{Proofs of path-counting results: Lemma \ref{lem:bad_paths} and Corollary \ref{cor:error_probs_small}}\label{ssec:proof_path_counting}

\begin{proof}[Proof of Lemma \ref{lem:bad_paths}]
Recall the notation $\mathscr{B}^r_x$ from (\ref{eq:defn_ball_L_r_x}), which denotes the open Euclidean ball of radius $r$ in $\mathbb{R}^d$ centered at $x \in \mathbb{R}^d$. 
We first prove (\ref{eq:lem_bad_path_j}).

We first bound $\mathbb{W}_j(\mathbf{0},x_1,\ldots,x_j)$ from above by $\Exp{\kappa(\|x_j-x_{j-1}\|,W_0,W_1)}$, where $W_0$ and $W_1$ are i.i.d. copies of the limiting weight distribution $W$ (recall (\ref{eq:conv_in_law_Wn_to_W})). Then we apply the change of variables
\[
z_i=x_i-x_{i-1},\; 1\leq i \leq j,
\]
where $x_0=\mathbf{0}$, and apply Fubini's theorem, to obtain

\begin{align*}
    & \int_{\mathbb{R}^d}\cdots \int_{\mathbb{R}^d} \mathbb{W}_j(\mathbf{0},x_1\ldots,x_j)  \ind{\|x_1\|<a^m}\prod_{i=1}^{j-2}\ind{\|x_{i+1}-x_i\|<a^{m^{i+1}}}\ind{\|x_j-x_{j-1}\|>a^{m^j}} dx_1\cdots dx_j 
        \\ & \leq \int_{\mathscr{B}^{a^m}_{\mathbf{0}}}\cdots \int_{\mathscr{B}^{a^{m^{j-1}}}_{\mathbf{0}}}\int_{\mathbb{R}^d \setminus \mathscr{B}^{a^{m^j}}_{\mathbf{0}}} \Exp{\kappa(\|z_j\|,W_0,W_1)} dz_jdz_{j-1}\cdots dz_1
        \\ & \leq C_0 (a^m)^d \cdots (a^{m^{j-1}})^d \int_{\mathbb{R}^d \setminus \mathscr{B}^{a^{m^j}}_{\mathbf{0}}} \Exp{\kappa(\|z\|,W_0,W_1)} dz, \numberthis \label{eq:int_wj_UB_1}
\end{align*}
for some constant $C_0>0$.
Recall the polynomial domination from (\ref{eq:power_law_tail_kappa}), and the assumption $\alpha>d$ in Theorem \ref{thm:main_LWC_SIRGs}.

We note, by first making a change of variables to bring the integral
{\[
\int_{\mathbb{R}^d \setminus \mathscr{B}^{a^{m^j}}_{\mathbf{0}}} \Exp{\kappa(\|z\|,W_0,W_1)} dz
\]}
down to an integral on $\mathbb{R}$, and then using (\ref{eq:power_law_tail_kappa}), that for any $m$ sufficiently large such that $a^{m^j}>t_0$,  

\begin{align*}
\int_{\mathbb{R}^d \setminus \mathscr{B}^{a^{m^j}}_{\mathbf{0}}} \Exp{\kappa(\|z\|,W_0,W_1)} dz \leq C_1 \frac{1}{a^{(\alpha-d)m^j}},\numberthis \label{eq:int_Expec_kappa_UB} 
\end{align*}
for some constant $C_1>0$.

Combining (\ref{eq:int_Expec_kappa_UB}) with (\ref{eq:int_wj_UB_1}), we note that for some constant $C_2>0$,

\begin{align*}
    &\int_{\mathbb{R}^d}\cdots \int_{\mathbb{R}^d} \mathbb{W}_j(\mathbf{0},x_1,\ldots,x_j)  \ind{\|x_1\|<a^m}\prod_{i=1}^{j-2}\ind{\|x_{i+1}-x_i\|<a^{m^{i+1}}}\ind{\|x_j-x_{j-1}\|>a^{m^j}} dx_1\cdots dx_j
    \\& \leq C_2 \frac{1}{a^{dm^j\left(\frac{\alpha}{d}-1-\frac{1}{m}-\cdots-\frac{1}{m^{j-1}} \right)}}
    \to 0, \numberthis \label{eq:bad_path_almost_exact_error_term}
\end{align*}
as $m \to \infty$, since $\frac{\alpha}{d}>1$. This finishes the proof of (\ref{eq:lem_bad_path_j}).

We next go into the proof of (\ref{eq:lem_bad_path_nj}). Recall the definition of $\mathbb{W}^{v_1,\ldots,v_j}_{n}(x_0,x_1,\ldots,x_j)$ from (\ref{eq:defn_W_n}), where $v_0=U_n$ is the uniformly chosen vertex of $\mathbb{G}_n$.

We use the notations
\[
\vec{v}=v_1,\dots,v_j; \;\; \vec{x}=(x_0,x_1,\dots,x_j).
\]

For fixed $x_0,\ldots,x_j \in \mathbb{R}^d$, we define the function $\mathcal{W}^{\vec{x}}_n:\mathbb{R}^{j+1} \to \mathbb{R}$ as
\begin{align*}
    \mathcal{W}^{\vec{x}}_n(\vec{t}):= \kappa_n\left(\|x_1-x_0\|,t_0,t_1\right)\cdots\kappa_n\left(\|x_j-x_{j-1}\|,t_{j-1},t_j\right),\numberthis \label{eq:defn_calW_n_prod_path}
\end{align*}
for $\vec{t}=(t_0,t_1,\dots,t_j) \in \mathbb{R}^{j+1}$.

Note that 
\begin{align*}
    & \frac{1}{n^j}\sum_{v_1,\ldots,v_j \in [n]} \mathbb{W}^{\vec{v}}_{n}(\vec{x}) = \Exp{\frac{1}{n^{j+1}} \sum_{i_0,i_1,\ldots,i_j} \mathcal{W}^{\vec{x}}_{n}(W^{(n)}_{i_0},\ldots,W^{(n)}_{i_j})}\numberthis \label{eq:bb_Wn_to_cal_W_n}
\end{align*}
where we have used (\ref{eq:defn_W_n}) and the fact that $v_0=U_n$ is uniformly distributed over $n$.

Since clearly 
\begin{align*}
    \Exp{\frac{1}{n^{j+1}} \sum_{i_0,i_1,\ldots,i_j} \mathcal{W}^{\vec{x}}_{n}(W^{(n)}_{i_0},\ldots,W^{(n)}_{i_j})} = \Exp{\mathcal{W}^{\vec{x}}_{n}(W^{(n)}_{U_{n,0}},\ldots,W^{(n)}_{U_{n,j}})},
\end{align*}
where $U_{n,0},\ldots,U_{n,j}$ is an i.i.d.\ collection of uniformly distributed random variables on $[n]=\{1,\ldots,n\}$, we can take $h_n$ in (\ref{eq:conv_bdd_func_weghts}) to be $\mathcal{W}^{\vec{x}}_{n}$ to conclude that (recall the definition of $\mathbb{W}_j(\vec{x})=\mathbb{W}_j(x_0,x_1\cdots,x_j)$ from (\ref{eq:defn_W}))
\begin{align*}
    \frac{1}{n^j}\sum_{v_1,\ldots,v_j \in [n]} \mathbb{W}^{\vec{v}}_{n}(\vec{x}) = \Exp{\mathcal{W}^{\vec{x}}_{n}(W^{(n)}_{U_{n,0}},\ldots,W^{(n)}_{U_{n,j}})} \to \mathbb{W}_j(\vec{x}),\numberthis \label{eq:strong_cty_paths}
\end{align*}
as $n \to \infty$. Now, (\ref{eq:lem_bad_path_nj}) can be concluded using (\ref{eq:strong_cty_paths}), a routine change of variables, Fatou's lemma, and (\ref{eq:lem_bad_path_j}).

\end{proof}
\begin{remark}[Efficacy of our bounds]
In the proof of Lemma \ref{lem:bad_paths}, we  {have} bounded $\mathbb{W}_j(\mathbf{0},x_1,\ldots,x_j)$ from above by 
\[
\Exp{\kappa(\|x_j-x_{j-1}\|,W_0,W_1)}.
\] 
That is, we have bounded all  {except} the last term in the product inside the expectation \\$\mathbb{W}_j(\mathbf{0},x_1,\ldots,x_j)$ by $1$. This is usually a poor bound. However, as we see in the proof, this loss is well compensated by the strong double-exponential growth of $r=r(a,m,K)$ (recall \ref{eq:choice_of_r}). In particular, for our purpose, we have been able to successfully avoid the question of how correlated the random variables $\kappa(\|x_0-x_1\|,W^{(0)},W^{(1)})$ and $\kappa(\|x_1-x_2\|,W^{(1)},W^{(2)})$, where $x_0,x_1,x_2 \in \mathbb{R}^d$, $W^{(0)},W^{(1)},W^{(2)}$ are i.i.d.\ copies of the limiting weight distribution, are.  {We believe this question to be hard} to tackle in general,  {under our general assumptions on $\kappa$ as formulated} in Assumption \ref{sssec:assumption_connections}.   
\end{remark}

Next we go into the proof of Corollary \ref{cor:error_probs_small}:

\begin{proof}[Proof of Corollary \ref{cor:error_probs_small}]
Recall the definition of the distance of an edge from a vertex from (\ref{eq:defn_dist_edge_vertex}).  {Also recall the abbreviations in} (\ref{eq:imp_notations}).

We begin by analysing the event $\{BF_{n,r}\neq B_n\}$. Note that if $\|Y^{(n)}_i-Y^{(n)}_j\|<a^{m^{L+1}}$ for every edge $\{i,j\}$ in $B_n$  {that} is at graph distance $L$ ($0\leq L \leq K-1$) from the root $U_n$, then $B_n$ is a subgraph of $F_{n,r}$ with the same root ${U_n}$, which implies that $BF_{n,r}=B_n$.
Hence, the event $\{BF_{n,r}\neq B_n\}$ implies the event \begin{equation}\label{eq:defn_Bad_rn_event}
    \textbf{Bad}_{r,n}:=
    \{\text{there is some \emph{bad edge} in}\;B_n\},
\end{equation}  
where a \emph{bad edge} is an edge $\{i,j\}$ in $B_n$ with $\|Y^{(n)}_i-Y^{(n)}_j\|>a^{m^{L+1}}$, where $0\leq L \leq K-1$ is the graph distance of the edge $\{i,j\}$ from the root $U_n$ of $B_n$.

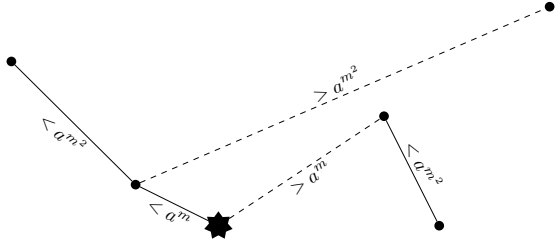
\begin{figure}[H]
    \centering
    \resizebox{0.5\textwidth}{!}{
        \begin{tikzpicture} 

\coordinate[rootnode] (O) at (3,0);
\coordinate[notneighborhood2] (goodpathnode1) at (1.5,0.75);
\coordinate[notneighborhood2] (goodpathnode2) at (-0.75,3);
\coordinate[notneighborhood2] (goodbadpathnode) at (9,4);
\coordinate[notneighborhood2] (badpathnode1) at (6,2);
\coordinate[notneighborhood2] (badpathnode2) at (7,0);
\draw (O) -- node[sloped,below, align=center] {$<a^m$} (goodpathnode1) -- node[sloped,below, align=center] {$<a^{m^2}$} (goodpathnode2);
\draw[dashed] (O) -- node[sloped,below, align=center] {$>a^m$} (badpathnode1);
\draw[dashed] (goodpathnode1) -- node[sloped,above, align=center] {$>a^{m^2}$} (goodbadpathnode);
\draw (badpathnode1) -- node[sloped,above, align=center] {$<a^{m^2}$} (badpathnode2);

\end{tikzpicture}
    }
    \caption{Illustration demonstrating bad edges. The star is the root. The dashed edges are bad, each of them connect pairs of vertices whose locations are at least $a^{m^{L+1}}$ apart, where $L$ is the distance of the edge from the root.  {The undashed edges are good.}}
    \label{fig:my_label}
\end{figure}

Therefore, 
    \begin{equation}\label{eq:UB_Bad_n}
    \Prob{BF_{n,r}\neq B_n}\leq \Prob{\textbf{Bad}_{r,n}}.
    \end{equation}
By a similar argument, 
\begin{equation}\label{eq:UB_Bad_r}
    \Prob{BF_{r}\neq B}\leq \Prob{\textbf{Bad}_r},
\end{equation}
where the event $\textbf{Bad}_r$ is similarly defined for the rooted graph $(\mathbb{G}_{\infty},0)$.

Define the event $\mathcal{I}_{n,j}$, for $n \in \mathbb{N},j \in [n]$, as
\begin{equation}
    \label{eq:defn_I_nj_event}
    \begin{split}
        \mathcal{I}_{n,j}:= &
    \{\exists\;v_1,\ldots,v_j \in V(\mathbb{G}_n)=[n]:\{v_0,v_1,\ldots,v_j\}\;\text{is a}\;j\text{-path in}\;(\mathbb{G}_n,U_n)\; \text{starting from} \\ & \text{the root}\;v_0=U_n,\; 
     \|Y^{(n)}_{v_{i-1}}-Y^{(n)}_{v_i}\|<a^{m^{i}} \forall  {i \in[j-1]}, \|Y^{(n)}_{v_{j-1}}-Y^{(n)}_{v_j}\|>a^{m^{j}} \}.
    \end{split}
\end{equation}
 {A} simple union bound gives 
    \begin{equation}
    \label{eq:UB_bad_path_lengths}
  \Prob{\textbf{Bad}_{r,n}}\leq \sum_{j=1}^K\Prob{ \mathcal{I}_{n,j}}.
    \end{equation}
Similarly,
\begin{equation}
    \label{eq:UB_Bad_r_sum_I_j events}
    \Prob{\textbf{Bad}_r}
 \leq \sum_{j=1}^K \Prob{\mathcal{I}_j},
    \end{equation}
where the event $\mathcal{I}_j$ for $j \in \mathbb{N}$ is similarly defined for the rooted graph $(\mathbb{G}_{\infty},0)$.

 {Note that since $K$ is fixed, it suffices to prove that $\lim_{m \to \infty} \limsup_{n \to \infty} \Prob{\mathcal{I}_{n,j}} = 0$ and $\lim_{m \to \infty} \limsup_{n \to \infty} \Prob{\mathcal{I}_{j}} = 0$ for $1 \le j \le K$.}

We proceed by bounding the probabilities $\Prob{\mathcal{I}_{n,j}}$ and $\Prob{\mathcal{I}_{j}}$. Recall  {that} we use $v_0$ to denote the typical vertex $U_n$ of $\mathbb{G}_n$ from the definition of the event $\mathcal{I}_{n,j}$ from (\ref{eq:defn_I_nj_event}). Note that, by Markov's inequality,
    \begin{equation}
    \label{eq:prob_I_nj}
    \begin{split}
    &    \Prob{\mathcal{I}_{n,j}}
    \leq 
    \\& \sum_{{v_1},\ldots,{v_j} \in [n]}
    \Prob{\{v_0,v_1,\ldots,v_j\}\;j\text{-path},\; \|Y^{(n)}_{v_{i-1}}-Y^{(n)}_{v_i}\|<a^{m^{i}}\; \forall  {i \in[j-1]}, \|Y^{(n)}_{v_{j-1}}-Y^{(n)}_{v_j}\|>a^{m^{j}}}.    
    \end{split}
    \end{equation}

For convenience, for $v_0,\ldots,v_j \in [n]$, and using the notation $\vec{v}=(v_0,v_1,\ldots,v_n)$, we denote the event 
\begin{equation}\label{eq:defn_bad_location_event_n}
    \mathcal{E}_{n,j}(\vec{v}):=\{\|Y^{(n)}_{v_{i-1}}-Y^{(n)}_{v_i}\|<a^{m^{i}}\; \forall  {i \in[j-1]}, \|Y^{(n)}_{v_{j-1}}-Y^{(n)}_{v_j}\|>a^{m^{j}}\}.
\end{equation}
We compute
\begin{align*}
    & \Prob{\{v_0,v_1,\ldots,v_j\}\;j\text{-path},\; \|Y^{(n)}_{v_{i-1}}-Y^{(n)}_{v_i}\|<a^{m^{i}}\; \forall  {i \in[j-1]}, \|Y^{(n)}_{v_{j-1}}-Y^{(n)}_{v_j}\|>a^{m^{j}}} 
    \\& = \Exp{\indE{\mathcal{E}_{n,j}(\vec{v})}\;\CProb{\{v_0,v_1,\ldots,v_j\}\;j\text{-path}}{Y^{(n)}_{v_0},\ldots,Y^{(n)}_{v_j},W^{(n)}_{v_0},\ldots,W^{(n)}_{v_j}}} 
    \\& = \Exp{\indE{\mathcal{E}_{n,j}(\vec{v})}\;\kappa_n\left(\|Y^{(n)}_{v_0}-Y^{(n)}_{v_1}\|,W^{(n)}_{v_0},W^{(n)}_{v_1}\right)\cdots\kappa_n\left(\|Y^{(n)}_{v_{j-1}}-Y^{(n)}_{v_j}\|,W^{(n)}_{v_{j-1}},W^{(n)}_{v_j}\right)}.\numberthis\label{eq:comput_prob_I_j_sum}
\end{align*}

Recall the notation $\mathbb{W}^{v_1,\ldots,v_j}_{n}$ from (\ref{eq:defn_W_n}). Using Fubini's theorem and that  $\{Y^{(n)}_{v_0},Y^{(n)}_{v_1},\ldots,Y^{(n)}_{v_j}\}$ is an i.i.d.\ collection of $j+1$ uniform random variables on $I_n$, recalling (\ref{eq:defn_bad_location_event_n}), and using the notation $\vec{x}=(x_0,x_1,\dots,x_j)$, (\ref{eq:comput_prob_I_j_sum}) becomes       
\begin{equation}
    \label{eq:comput_prob_I_nj_integral}
    \begin{split}
    \frac{1}{n^{j+1}} \int_{I_n}\cdots\int_{I_n} \mathbb{W}^{v_1,\ldots,v_j}_{n}(\vec{x}) \prod_{i=0}^{j-2}\ind{\|x_{i+1}-x_i\|<a^{m^{i+1}}}\ind{\|x_{j}-x_{j-1}\|>a^{m^j}} dx_{0}\cdots dx_j.   
    \end{split}
\end{equation}

Using (\ref{eq:prob_I_nj}) and (\ref{eq:comput_prob_I_nj_integral}), we have
\begin{equation}\label{eq:prob_I_nj_UB_integral}
\begin{aligned}
    \Prob{\mathcal{I}_{n,j}} 
        &\leq \frac{1}{n^{j+1}}\int_{I_n}\cdots\int_{I_n} \sum_{v_1,\ldots,v_j \in [n]}\mathbb{W}^{v_1,\ldots,v_j}_{n}(\vec{x}) \\ &\hspace{10pt} \times \prod_{i=0}^{j-2}\ind{\|x_{i+1}-x_i\|<a^{m^{i+1}}}\ind{\|x_{j}-x_{j-1}\|>a^{m^j}} dx_{0}\cdots dx_j. 
\end{aligned}
\end{equation}
Recall the notation $\mathbb{W}_j$ from (\ref{eq:defn_W}). Similarly,  {by the multivariate Mecke formula for Poisson processes \cite[Theorem 4.4]{Last_Penrose_LPP}}, 
\begin{equation}\label{eq:prob_I_j_UB_integral_3_Mecke}
    \begin{aligned}
        \Prob{\mathcal{I}_j}
        &\leq \int_{\mathbb{R}^d}\cdots \int_{\mathbb{R}^d} \mathbb{W}_j(\mathbf{0},x_1,\ldots,x_j)  \ind{\|x_1\|<a^m}\\
        &\hspace{10pt} \times \prod_{i=1}^{j-2}\ind{\|x_{i+1}-x_i\|<a^{m^{i+1}}}\ind{\|x_j-x_{j-1}\|>a^{m^j}} dx_1\dots dx_j.
    \end{aligned}
\end{equation}

Now we apply Lemma \ref{lem:bad_paths} for the bounds (\ref{eq:prob_I_j_UB_integral_3_Mecke}) and (\ref{eq:comput_prob_I_nj_integral}), and use the bounds (\ref{eq:UB_Bad_n}), (\ref{eq:UB_Bad_r}), (\ref{eq:UB_bad_path_lengths}) and (\ref{eq:UB_Bad_r_sum_I_j events}) to conclude Corollary \ref{cor:error_probs_small}.

\end{proof}

\begin{remark}[A general estimate]\label{rem:gen_estimate}
Recall the notations (\ref{eq:imp_notations}) and recall $r(a,m,K)$ from (\ref{eq:choice_of_r}). Note that the bound (\ref{eq:UB_Bad_n}) is true even when one has $a=a_n$ and $K=K_n$ in $r=r(a_n,m,K_n)$, and so is the simple union bound (\ref{eq:UB_bad_path_lengths}). The integral bound (\ref{eq:prob_I_nj_UB_integral}) also works in this generality. Recalling the definition of the function $\mathcal{W}^{(x_0,\ldots,x_j)}_n$ from (\ref{eq:defn_calW_n_prod_path}), the equality (\ref{eq:bb_Wn_to_cal_W_n}), and the display below it, bounding $\Exp{\mathcal{W}^{(x_0,\ldots,x_j)}_{n}(W^{(n)}_{U_{n,0}},\ldots,W^{(n)}_{U_{n,j}})}$ from above by $\Exp{\kappa_n(\|x_j-x_{j-1}\|,W^{(n)}_{U_{n,j}},W^{(n)}_{U_{n,j-1}})}$ (where $U_{n,0},\dots,U_{n,j}$ are i.i.d.\ uniformly distributed random variables on $[n]$), and making an easy change of variable, we obtain the general bound
\begin{align*}
        & \Prob{B^{F^{\mathbb{G}_n}_{U_{n}}(r_n)}_{U_{n}}(K_n) \neq B^{\mathbb{G}_n}_{U_{n}}(K_n)} \\
    &  \leq \sum_{j=1}^{K_n} \int_{\mathscr{B}^{a_n^m}_{\mathbf{0}}}\cdots\int_{\mathscr{B}^{a_n^{m^{j-1}}}_{\mathbf{0}}}\int_{\mathbb{R}^d \setminus \mathscr{B}^{a_n^{m^j}}_{\mathbf{0}}} \Exp{\kappa_n(\|z_j\|,W^{(n)}_{\mathrm{U}'_{n,1}},W^{(n)}_{\mathrm{U}'_{n,2}})} dz_j\cdots dz_1, \numberthis \label{eq:rem_gen_bd}     
\end{align*}
where $\mathrm{U}'_{n,1}$ and $\mathrm{U}'_{n,2}$ are two i.i.d.\ uniform elements in $[n]$.  {Below, we} will use the bound (\ref{eq:rem_gen_bd}) with suitable choices of $a=a_n$ and $K=K_n$ to prove Theorem \ref{thm:main_typical_dist}.
\end{remark}

\subsection{Proof of Theorem \ref{thm:main_LWC_SIRGs}}\label{ssec:proof_LWC}
Proposition \ref{prop:loc_trunc} implies  {that} the random rooted graph $(F^{\mathbb{G}_n}_{U_n}(r),U_n)$ converges in distribution to the random rooted graph $(F^{\mathbb{G}_{\infty}}_0(r),0)$ in the space $\mathcal{G}_{\star}$. The proof of this fact can be carried out in the same manner as \cite[Theorem 2.13]{RGCN_2} is proved assuming \cite[Definition 2.10]{RGCN_2}, and so we leave this for the reader to check. In particular, as a consequence of Proposition \ref{prop:loc_trunc}, 
\begin{equation}\label{eq:law_F_n_to_F}
    \Prob{(F^{\mathbb{G}_n}_{U_n}(r),U_n) \in A} \to \Prob{(F^{\mathbb{G}_{\infty}}_{0}(r),U_n) \in A},
\end{equation}
for any subset $A \subset \mathcal{G}_{\star}$.

\begin{proof}[Proof of Theorem \ref{thm:main_LWC_SIRGs}.]
Recall the definition of local weak convergence from Definition \ref{defn:LWC}. Recall the  {abbreviations in} (\ref{eq:imp_notations}).
Let $G_*=(G,g) \in \mathcal{G}_\star$.

Fix $K \in \mathbb{N}$, and $\varepsilon >0$. To conclude Theorem \ref{thm:main_LWC_SIRGs}, we need to find a $N \in \mathbb{N}$ such that for all $n>N$, 
\begin{equation}\label{eq:K_balls_close}
    \left|\Prob{B^{\mathbb{G}_n}_{{U_n}}(K)\cong (G,g)}-\Prob{B^{\mathbb{G}_{\infty}}_{0}(K)\cong (G,g)}\right|
    < \varepsilon,
\end{equation}
that is  
\[
\left| \Prob{B_n \cong (G,g)} - \Prob{B \cong (G,g)} \right|< \varepsilon.
\]

Note that
    \begin{equation}
    \label{eq:B=F+error}
    \left|\Prob{B_n\cong (G,g)}-\Prob{B\cong (G,g)}\right|
    \leq \left|\Prob{BF_{n,r}\cong (G,g)}-\Prob{BF_r \cong (G,g)}\right|
    +\left|\varepsilon_{n,r}\right|
    +\left|\varepsilon_r\right|,
    \end{equation}
where 
    \begin{equation}
    \label{eq:error_nr}
    \varepsilon_{n,r}= 
    \Prob{B_n \cong (G,g), BF_{n,r}\neq B_n}
    -\Prob{BF_{n,r} \cong (G,g), BF_{n,r}\neq B_n}, 
    \end{equation} 
and
    \begin{equation}
    \label{eq:error_r}
    \varepsilon_{r}= 
    \Prob{B \cong (G,g), BF_{r}\neq B}
    -\Prob{BF_{r} \cong (G,g), BF_{r}\neq B}. 
    \end{equation}

Clearly, 
\begin{equation}\label{eq:UB_error_nr}
    |\varepsilon_{n,r}|\leq  \Prob{BF_{n,r}\neq B_n},
\end{equation}
and 
\begin{equation}\label{eq:UB_error_r}
    |\varepsilon_r|\leq  \Prob{BF_{r}\neq B}.
\end{equation}

For the rest of the proof, we fix $m>0$ and $n_0 \in \mathbb{N}$ such that, for all $n \geq n_0$,
    \begin{equation}
    \label{eq:error_small}
    |\varepsilon_{n,r}|+|\varepsilon_r|< \varepsilon/2,
    \end{equation}
which is possible by Corollary \ref{cor:error_probs_small}.

Note that 
    \begin{align*}
    \Prob{BF_{n,r}\cong (G,g)}=\Prob{F_{n,r} \in A(K,(G,g))}, 
    \end{align*}
where $A(K,(G,g))\subset \mathcal{G}_{\star}$ is defined as
\begin{align*}
    A(K,(G,g)):=\{ (H,h) \in \mathcal{G}_{\star}: B^{H}_h(K) \cong (G,g)\} \numberthis \label{eq:defn_subsets_with_given_nbd}
\end{align*}

By (\ref{eq:law_F_n_to_F}), as $n \to \infty$, 
    \begin{equation}
    \label{eq:using_prop_loctrunc}
    \Prob{F_{n,r} \in A(K,(G,g))} 
    \to \Prob{F_{r} \in A(K,(G,g))}=\Prob{BF_{r}\cong (G,g)}.
\end{equation} 
Combining (\ref{eq:using_prop_loctrunc}) with (\ref{eq:error_small}) and (\ref{eq:B=F+error}), we can choose $n_1 \in \mathbb{N}$ such that for all $n>N=\max\{n_1,n_0\}$, (\ref{eq:K_balls_close}) holds.
This completes the proof of Theorem \ref{thm:main_LWC_SIRGs}.
\end{proof}

\subsection{Proof of Theorem \ref{thm:main_LCinP_SIRGs}}\label{ssec:proof_LCinP}

For $ (H,h) \in \mathcal{G}_{\star}$, and $r>0$, recall the empirical Euclidean graph neighborhood distribution as defined in (\ref{eq:defn_C_rn_rvs}). We first give the proof of Lemma \ref{Lem:sec_mom_C_rn}.

\begin{proof}[Proof of Lemma \ref{Lem:sec_mom_C_rn}]
To ease notation, let us write $C_{r,n}$ for $C_{r,n} (H,h)$.

Note that
\begin{align*}
  C_{r,n}^2 
  & = \frac{1}{n^2}\sum_{i,j=1}^n\ind{F^{\mathbb{G}_n}_i(r)\cong  (H,h)}\ind{F^{\mathbb{G}_n}_j(r)\cong  (H,h)} \\
  & = \CProb{F^{\mathbb{G}_n}_{U_{n,1}}(r)\cong  (H,h),F^{\mathbb{G}_n}_{U_{n,2}}(r)\cong  (H,h)}{\mathbb{G}_n},
\end{align*}
where $U_{n,1},U_{n,2}$ are i.i.d.\ uniformly distributed random variables on $[n]$. Therefore,

\begin{equation}\label{eq:nbhd_two_indep_nodes}
    \Exp{C_{r,n}^2}=\Prob{F^{\mathbb{G}_n}_{U_{n,1}}(r)\cong  (H,h),F^{\mathbb{G}_n}_{U_{n,2}}(r)\cong  (H,h)}.
\end{equation}

We introduce the following  {abbreviations} which we will use throughout this proof to keep  {notation concise} (recall the point process $\Gamma_n$ from (\ref{eq:defn_Gamma_n}) of the locations of the vertices of $\mathbb{G}_n$, the set $A^r_n$ from (\ref{eq:defn_set_Ar_n}) and the ball $\mathscr{B}^r_x$ from (\ref{eq:defn_ball_L_r_x})):
\begin{equation}
    \label{eq:imp_notations_2}
    \begin{split}
     \mathscr{E}&:=\left\{F^{\mathbb{G}_n}_{U_{n,1}}(r)\cong  (H,h)\right\},\quad \mathscr{F}:=\left\{F^{\mathbb{G}_n}_{U_{n,2}}(r)\cong  (H,h)\right\},\\ 
    \mathscr{U}&:=\left\{\Gamma_n\left(\mathscr{B}^{r}_{Y^{(n)}_{U_{n,1}}}\right)=|V(H)|\right\},\quad \mathscr{V}:=\left\{\Gamma_n\left(\mathscr{B}^{r}_{Y^{(n)}_{U_{n,2}}}\right)=|V(H)|\right\},\\
    \mathscr{J}&:=\left\{\mathscr{B}^r_{Y^{(n)}_{U_{n,1}}}\cap \mathscr{B}^r_{Y^{(n)}_{U_{n,2}}}=\varnothing\right\},\quad \mathscr{W}:=\left\{U_{n,1},U_{n,2} \in A^{r}_n\right\}.
    \end{split}
\end{equation}

Recall  {that} the target is to show
\begin{equation}\label{eq:LCinP_target}
  |\Prob{\mathscr{E} \cap \mathscr{F}}-\Prob{\mathscr{E}}\Prob{\mathscr{F}}|=|\Prob{\mathscr{E}\cap \mathscr{F}}-\Prob{\mathscr{E}}^2| \to0,  
\end{equation} 
as $n \to \infty$ 
(note that conditionally on $\mathbb{G}_n$, the random variables $\indE{\mathscr{E}}$ and $\indE{\mathscr{F}}$ are identically distributed, and hence they have the same conditional expectation, and hence same expectation).

 {First we write,}
\begin{align*}
  \Prob{\mathscr{E}\cap \mathscr{F}}&=\CProb{\mathscr{E}\cap \mathscr{F}}{\mathscr{J}\cap \mathscr{U}\cap \mathscr{V} \cap \mathscr{W}}\Prob{\mathscr{J}\cap \mathscr{U}\cap \mathscr{V}\cap \mathscr{W}} \numberthis \label{eq:LWP_Term1}\\
  &\qquad+\Prob{\{\mathscr{E}\cap \mathscr{F}\}\cap \{\mathscr{J}^c\cup \mathscr{U}^c \cup \mathscr{V}^c \cup \mathscr{W}^c\}} \numberthis \label{eq:LWP_Term2}.
\end{align*}

 {Note} that the term  {in} (\ref{eq:LWP_Term2}) is bounded from above by 
\[\Prob{\{\mathscr{E}\cap \mathscr{F}\}\cap \{\mathscr{U}^c\cup \mathscr{V}^c\}}+\Prob{\mathscr{J}^c}+\Prob{\mathscr{W}^c},\] 
and it is easily observed that the first term is equal to $0$, and the last term  {tends} to $0$ as $n \to \infty$ using (\ref{eq:aas_typical_in_Arn}). Also note that \[\Prob{\mathscr{J}^c}\leq \Prob{Y^{(n)}_{U_{n,1}} \in \mathscr{B}^{2r}_{Y^{(n)}_{U_{n,2}}} \cap I_n}.\] 
 
Clearly, 
\[\CProb{Y^{(n)}_{U_{n,1}} \in \mathscr{B}^{2r}_{Y^{(n)}_{U_{n,2}}} \cap I_n}{Y^{(n)}_{U_{n,2}}}\stackrel{\text{a.s.}}{\leq} \frac{\lambda_d(\{y \in \mathbb{R}^d: \|y\|< 2r\})}{n},\] 
which tends to $0$ as $n \to \infty$. Hence taking expectations of both sides in the last display and letting $n \to \infty$, we get $\Prob{\mathscr{J}^c}\to0$. Hence the term  {in} (\ref{eq:LWP_Term2})  {tends} to $0$ as $n \to \infty$.

 {To analyse the term in} (\ref{eq:LWP_Term1}), we observe that, conditionally on $\mathscr{J}\cap \mathscr{U}\cap \mathscr{V}\cap \mathscr{W}$,
the random variables $\indE{\mathscr{E}}$ and $\indE{\mathscr{F}}$ are independent, since they are just functions of the locations of the $|V(H)|$ many points falling in $\mathscr{B}^r_{Y^{(n)}_{U_{n,1}}}$ and $\mathscr{B}^r_{Y^{(n)}_{U_{n,2}}}$, and these locations are independent (since the locations of different vertices of $\mathbb{G}_n$ are independent).
Hence, 
\[(\ref{eq:LWP_Term1})=\CProb{\mathscr{E}}{\mathscr{J}\cap \mathscr{U}\cap \mathscr{V}\cap \mathscr{W}}\CProb{\mathscr{F}}{\mathscr{J}\cap \mathscr{U}\cap \mathscr{V}\cap \mathscr{W}} \Prob{\mathscr{J}\cap \mathscr{U}\cap \mathscr{V}\cap \mathscr{W}}.\]

Now note that $\mathscr{E}$ is independent of $\mathscr{V}$,  {conditionally} on $\{\mathscr{J}\cap \mathscr{U}\cap \mathscr{W}\}$. Similarly, $\mathscr{F}$ is independent of $\mathscr{U}$,  {conditionally} on $\{\mathscr{J}\cap \mathscr{V}\cap \mathscr{W}\}$.

Hence, 
\[(\ref{eq:LWP_Term1})=\CProb{\mathscr{E}}{\mathscr{J}\cap \mathscr{U}\cap \mathscr{W}}\CProb{\mathscr{F}}{\mathscr{J}\cap \mathscr{V}\cap \mathscr{W}} \Prob{\mathscr{J}\cap \mathscr{U}\cap \mathscr{V}\cap \mathscr{W}}.\]

As argued earlier $\Prob{\mathscr{J}^c}\to0$, so $\Prob{\mathscr{J}}\to1$, and it is easy to observe that $\Prob{\mathscr{W}}\to1$ using (\ref{eq:aas_typical_in_Arn}). Hence, we can forget about the `almost' certain events $\mathscr{W}$ and $\mathscr{J}$ for the first two terms from the last display, and condition on $\mathscr{J}\cap \mathscr{W}$ for the third term, to conclude that (recall (\ref{eq:LWP_Term1})) 
\begin{align*}
    & |\CProb{\mathscr{E}\cap \mathscr{F}}{\mathscr{J}\cap \mathscr{U}\cap \mathscr{V} \cap \mathscr{W}}\Prob{\mathscr{J}\cap \mathscr{U}\cap \mathscr{V}\cap \mathscr{W}} \\ &-  \CProb{\mathscr{E}}{\mathscr{U}}\CProb{\mathscr{F}}{\mathscr{V}}\CProb{\mathscr{U}\cap \mathscr{V}}{\mathscr{J}\cap \mathscr{W}}|\to 0,\numberthis \label{eq:almost_Term1}
\end{align*}
that is, the difference between (\ref{eq:LWP_Term1}) and $\CProb{\mathscr{E}}{\mathscr{U}}\CProb{\mathscr{F}}{\mathscr{V}}\CProb{\mathscr{U}\cap \mathscr{V}}{\mathscr{J}\cap \mathscr{W}}$  {tends} to $0$. 

We claim that to conclude the proof, it suffices to check that
\begin{equation}\label{LWP_finalbit}
  |\CProb{\mathscr{U}\cap \mathscr{V}}{\mathscr{J}\cap \mathscr{W}}-\Prob{\mathscr{U}}\Prob{\mathscr{V}}|\to0.  
\end{equation}
This is because (\ref{LWP_finalbit}) combined with (\ref{eq:almost_Term1}) and  the observations that $\mathscr{E} \subset \mathscr{U}$ and $\mathscr{F} \subset \mathscr{V}$ implies the difference between the expression (\ref{eq:LWP_Term1}) and $\Prob{\mathscr{E}}\Prob{\mathscr{F}}$ goes to $0$. Combining this with the fact that the expression  {in} (\ref{eq:LWP_Term2}) goes to $0$, we obtain (\ref{eq:LCinP_target}).

We now show (\ref{LWP_finalbit}).
Recall  {that} $|V(H)|$ denotes the size of the vertex set of the graph $H$. Observe that 
\[\CProb{\mathscr{U}\cap \mathscr{V}}{\mathscr{J}\cap \mathscr{W}}=\Prob{(M_1,M_2,M_3)=(|V(H)|-1,|V(H)|-1,n-2|V(H)|)},\]
where $\textbf{M}=(M_1,M_2,M_3)$ is a multinomial vector with parameters $(n-2;\frac{\lambda_d(\mathscr{B}^{r}_{\mathbf{0}})}{n},\frac{\lambda_d(\mathscr{B}^{r}_{\mathbf{0}})}{n},1-2\frac{\lambda_d(\mathscr{B}^{r}_{\mathbf{0}})}{n}).$ Hence,
\begin{align*}
 & \CProb{\mathscr{U}\cap \mathscr{V}}{\mathscr{J}\cap \mathscr{W}}\\
& =  \frac{n!}{(|V(H)|-1)!\;(|V(H)|-1)!\;(n-2|V(H)|)!}\\
& \times \left(\frac{\lambda_d(\mathscr{B}^{r}_{\mathbf{0}})}{n}\right)^{|V(H)|-1}\left(\frac{\lambda_d(\mathscr{B}^{r}_{\mathbf{0}})}{n}\right)^{|V(H)|-1}\left(1-2\frac{\lambda_d(\mathscr{B}^{r}_{\mathbf{0}})}{n}\right)^{n-2|V(H)|}. 
\end{align*}
It is an easy analysis to check that this converges to $\Prob{Y=|V(H)|-1}^2$, where $Y\sim \text{Poi}(\lambda_d(\mathscr{B}^{r}_{\mathbf{0}}))$.

Again,
\[\CProb{\mathscr{U}}{Y^{(n)}_{U_{n,1}}\in A^{r}_n}\CProb{\mathscr{V}}{Y^{(n)}_{U_{n,2}}\in A^{r}_n}=\Prob{Y'_n=|V(H)|-1}^2,\]
where $Y'_n \sim \text{Bin}(n-1, \frac{\lambda_d(\mathscr{B}^{r}_{\mathbf{0}})}{n})$.

Since $\Prob{Y'_n=|V(H)|-1}^2\to\Prob{Y=|V(H)|-1}^2$, and both \[\Prob{Y^{(n)}_{U_{n,1}}\in A^{r}_n}, \Prob{Y^{(n)}_{U_{n,2}}\in A^{r}_n} \geq \Prob{\mathscr{W}}\to 1\] 
(using (\ref{eq:aas_typical_in_Arn})), we have shown that
\begin{align*}
  &|\CProb{\mathscr{U}\cap \mathscr{V}}{\mathscr{J}\cap \mathscr{W}}-\Prob{\mathscr{U}}\Prob{\mathscr{V}}|\\
  &\leq \left|\CProb{\mathscr{U}\cap \mathscr{V}}{\mathscr{J}\cap \mathscr{W}}-\CProb{\mathscr{U}}{Y^{(n)}_{U_{n,1}}\in A^{r}_n}\CProb{\mathscr{V}}{Y^{(n)}_{U_{n,2}}\in A^{r}_n} \right|\\
  &+\left|\CProb{\mathscr{U}}{Y^{(n)}_{U_{n,1}}\in A^{r}_n}\CProb{\mathscr{V}}{Y^{(n)}_{U_{n,2}}\in A^{r}_n}-\Prob{\mathscr{U}}\Prob{\mathscr{V}} \right|\to0.
\end{align*}
This completes the proof of (\ref{LWP_finalbit}) and hence Lemma \ref{Lem:sec_mom_C_rn}.
\end{proof}

Note that using Lemma \ref{Lem:sec_mom_C_rn} with Proposition \ref{prop:loc_trunc}, a direct application of Chebyshev's inequality gives, for any $ (H,h) \in \mathcal{G}_{\star}$,
\[
C_{r,n} (H,h)\plim \Prob{F^{\mathbb{G}_{\infty}}_{0}(r)\cong  (H,h)}.
\]

In particular, this implies  {that} the empirical Euclidean graph neighborhood measure of $\mathbb{G}_n$ converges in probability to the measure induced by the random element $(F^{\mathbb{G}_{\infty}}_0(r),0)$ in $\mathcal{G}_{\star}$: for any subset $A \subset \mathcal{G}_{\star}$, and, for any $r>0$ (recall the Euclidean graph neighborhoods from Definition \ref{defn:F_graphs}),
\begin{equation}\label{eq:conv_Plim_C_n}
    \frac{1}{n}\sum_{i=1}^n \ind{(F^{\mathbb{G}_n}_i(r),i) \in A} \plim \Prob{(F^{\mathbb{G}_{\infty}}_{0}(r),0) \in A},
\end{equation} 
as $n \to \infty$. We now use (\ref{eq:conv_Plim_C_n}) to prove Theorem \ref{thm:main_LCinP_SIRGs}.

\begin{proof}[Proof of Theorem \ref{thm:main_LCinP_SIRGs}]
For $K \in \mathbb{N}$ and $(G,g) \in \mathcal{G}_{\star}$, recall the random variables $B_n(G,g)$ as defined in (\ref{eq:defn_B_n_rvs}).  {Further, recall} the definition of local convergence in probability from Definition \ref{defn:LPC}. 

Fix $\varepsilon>0$. Note that the target is to show that,  {for every $\varepsilon>0$,} there exists $N \in \mathbb{N}$ such that for all $n > N$, 
    \begin{equation}\label{eq:final_target}
     {\mathbb{P}\left(\Big|B_n(G,g)-\Prob{B^{\mathbb{G}_{\infty}}_{0}(K)\cong (G,g)}\Big|>\varepsilon\right)< \varepsilon.}
\end{equation}

We abbreviate the Euclidean graph neighborhoods $F^{n,r,i}:=F^{\mathbb{G}_n}_i(r)$, with $r=r(a,m,K)$ is as in (\ref{eq:choice_of_r}), and $i \in V(\mathbb{G}_n)=[n]$.
We note that 
\begin{align}
    B_n(G,g) & =\frac{1}{n}\sum_{i=1}^n \ind{B^{\mathbb{G}_n}_i(K)\cong (G,g)} \nonumber \\ &  =\frac{1}{n}\sum_{i=1}^n \ind{B^{\mathbb{G}_n}_i(K)\cong (G,g)}\ind{B^{\mathbb{G}_n}_i(K)=B^{F_{n,r,i}}_i(K)}\\
    & \hspace{10 pt}+\frac{1}{n}\sum_{i=1}^n \ind{B^{\mathbb{G}_n}_i(K)\cong (G,g)}\ind{B^{\mathbb{G}_n}_i(K)\neq B^{F_{n,r,i}}_i(K)}. \label{eq:B_n_decomp1} 
\end{align}

Writing the first term on the RHS of (\ref{eq:B_n_decomp1}) as \begin{equation} \label{eq:B_n_Term1_decomp}
    \frac{1}{n}\sum_{i=1}^n \ind{B^{F_{n,r,i}}_i(K)\cong (G,g)}-\frac{1}{n}\sum_{i=1}^n \ind{B^{F_{n,r,i}}_i(K)\cong (G,g)}\ind{B^{\mathbb{G}_n}_i(K)\neq B^{F_{n,r,i}}_i(K)},
\end{equation}

we note from (\ref{eq:B_n_decomp1}) that 
\begin{equation} \label{eq:B_n_decomp}
    B_n(G,g)=\frac{1}{n}\sum_{i=1}^n \ind{B^{F_{n,r,i}}_i(K)\cong (G,g)}+\varepsilon'_{n,r},
\end{equation}
where (recall the event $\textbf{Bad}_{r,n}$ from (\ref{eq:defn_Bad_rn_event})) the expectation $\Exp{{|\varepsilon'_{n,r}|}}$ of the {absolute value of the} error $\varepsilon'_{n,r}$ is bounded from above by $2\Prob{\textbf{Bad}_{r,n}}$.

It can be shown using a similar argument that (recall the notations (\ref{eq:imp_notations})) \begin{equation}\label{eq:limit_decomp}
\Prob{B^{\mathbb{G}_{\infty}}_{0}(K)\cong (G,g)}=\Prob{BF_r\cong (G,g)}+\varepsilon'_r,    
\end{equation}
where $\Exp{{|\varepsilon'_r|}}$ is bounded from above by $2\Prob{\textbf{Bad}_r}$.

Using Corollary \ref{cor:error_probs_small}, it is easy to see that there exists $n_0 \in \mathbb{N}$ such that for $m>0$ sufficiently large (recall $r=r(a,m,K)$), for all $n>n_0$, \begin{equation}\label{eq:error_estimate_final}
    2\Prob{\textbf{Bad}_{{r,n}}}+2\Prob{\textbf{Bad}_r}<\varepsilon/2.
\end{equation}

Recall the subset $A(K,(G,g)) \subset \mathcal{G}_{\star}$ from (\ref{eq:defn_subsets_with_given_nbd}). 
Note that using (\ref{eq:conv_Plim_C_n}),
\begin{align*}
 \frac{1}{n}\sum_{i=1}^n \ind{B^{F_{n,r,i}}_{i}(K)\cong (G,g)} 
 &=\frac{1}{n}\sum_{i=1}^n \ind{(F_{n,r,i},i) \in A(K,(G,g))} \\
 & \plim \Prob{(F^{\mathbb{G}_{\infty}}_0(r),0)\in A(K,(G,g))}\\
 &=\Prob{BF_r \cong (G,g)},
\end{align*} and so there exists $n_1 \in \mathbb{N}$ such that for all $n>n_1$ we have

\begin{equation}\label{eq:final_ineq}
     {\mathbb{P}\left(\Big|\frac{1}{n}\sum_{i=1}^n {\ind{B^{F_{n,r,i}}_{i}(K)\cong (G,g)}}-\Prob{BF_r \cong (G,g)}\Big|>\varepsilon/2\right)< \varepsilon/2.}
\end{equation}
Hence for all $n>N=\max\{n_0,n_1\}$, we note that using (\ref{eq:error_estimate_final}) and (\ref{eq:final_ineq}), (\ref{eq:final_target}) holds. This completes the proof of Theorem \ref{thm:main_LCinP_SIRGs}.
\end{proof}

\subsection{Proof of Theorem \ref{thm:main_typical_dist}}\label{ssec:proof_typ_dist}

The main idea in this proof is that the first result of Corollary \ref{cor:error_probs_small} can be pushed to the case when $K=K_n$ is allowed to grow in a doubly logarithmic manner, instead of being fixed. For this, we need a finer analysis of the error terms we encounter while proving Lemma \ref{lem:bad_paths}. We now go into the formal argument:

\begin{proof}

Recall $C \in \left( 0, \frac{1}{\log(\frac{\alpha}{\alpha-d})}\right)$.
Fix some $\overline{C} \in \left( C, \frac{1}{\log(\frac{\alpha}{\alpha-d})}\right)$,
and let 
\begin{equation}\label{eq:choice_m}
    m:= e^{1/\overline{C}}.
\end{equation}
Since $\overline{C}< \frac{1}{\log(\frac{\alpha}{\alpha-d})}$, 
\begin{align*}
    \frac{\alpha}{d}> \frac{m}{m-1}.
\end{align*}

Let
\begin{equation}\label{eq:dist_choice_an}
    a_n=\log n,
\end{equation}
and
\begin{equation}\label{eq:dist_choice_Kn}
    K_n = \frac{\log \log n + \log(\frac{1}{d}-\delta) - \log \log a_n}{\log m} = \frac{\log \log n + \log(\frac{1}{d}-\delta) - \log \log \log n}{\log m},
\end{equation}
where $\delta>0$ is such that $\frac{1}{d}-\delta >0$.
Note that by the choice of $K_n$, we have 
\begin{equation}\label{eq:choice_r_n}
    a_n^{m^{K_n}}=n^{\frac{1}{d}-\delta}.
\end{equation}
We also let $r_n=r(a_n,m,K_n)$ be as in (\ref{eq:choice_of_r}).

Note that, since $C < \frac{1}{\log m}$ by the choice of $m$,
for all large $n$,
\begin{align*}
    C \log \log n \leq K_n.
\end{align*}
Then for $U_{n,1}$ and $U_{n,2}$ two uniformly chosen vertices of $\mathbb{G}_n$, for all large $n$,
\begin{align*}
    & \Prob{d_{\mathbb{G}_n}(U_{n,1},U_{n,2})\leq C\log \log n}\\ 
    & \leq \Prob{d_{\mathbb{G}_n}(U_{n,1},U_{n,2})\leq K_n}\\
    &= \Prob{U_{n,2}\in B^{\mathbb{G}_n}_{U_{n,1}}(K_n)}\\
    & \leq \Prob{U_{n,2} \in B^{F^{\mathbb{G}_n}_{U_{n,1}}(r_n)}_{U_{n,1}}(K_n)} + \Prob{B^{F^{\mathbb{G}_n}_{U_{n,1}}(r_n)}_{U_{n,1}}(K_n) \neq B^{\mathbb{G}_n}_{U_{n,1}}(K_n)} \\
    & \leq \Prob{U_{n,2} \in F^{\mathbb{G}_n}_{U_{n,1}}(r_n)}+ \Prob{B^{F^{\mathbb{G}_n}_{U_{n,1}}(r_n)}_{U_{n,1}}(K_n) \neq B^{\mathbb{G}_n}_{U_{n,1}}(K_n)}\\
    & =\Prob{Y^{(n)}_{U_{n,2}}\in \mathscr{B}^{r_n}_{Y^{(n)}_{U_{n,1}}}} + \Prob{B^{F^{\mathbb{G}_n}_{U_{n,1}}(r_n)}_{U_{n,1}}(K_n) \neq B^{\mathbb{G}_n}_{U_{n,1}}(K_n)}, \numberthis \label{eq:dist_decomp}
\end{align*}
where  {we} recall the Euclidean graph neighborhoods $F^{\mathbb{G}_n}_i(r)$ from Definition \ref{defn:F_graphs}.

Note that it is sufficient to establish  {that} the RHS of (\ref{eq:dist_decomp})  {tends} to $0$ as $n \to \infty$, to conclude Theorem \ref{thm:main_typical_dist}.

Via a simple conditioning on $Y^{(n)}_{U_{n,1}}$, and using that $Y^{(n)}_{U_{n,2}}$ is uniformly distributed on $I_n$ and is independent of $Y^{(n)}_{U_{n,1}}$, we have
\begin{align*}
    & \Prob{Y^{(n)}_{U_{n,2}}\in \mathscr{B}^{r_n}_{Y^{(n)}_{U_{n,1}}}} \leq \frac{\lambda_d(\mathscr{B}^{r_n}_{\mathbf{0}})}{n}.
\end{align*}
 {Further,} note that (recall $r(a,m,K)$ from (\ref{eq:choice_of_r})) using the upper bound
\[
r(a_n,m,K_n)\leq K_n a_n^{m^{K_n}},
\]
we have 
\begin{align*}
    \frac{\lambda_d(\mathscr{B}^{r_n}_{\mathbf{0}})}{n} \leq \frac{w r_n^d}{n} \leq \frac{w K_n^d n^{d\delta-1}}{n}\to 0,
\end{align*}
as $n \to \infty$, by the choice of $K_n$ as in (\ref{eq:dist_choice_Kn}), where $w>1$ is some constant upper bound on $\lambda_d(\mathscr{B}^1_{\mathbf{0}})$, and in the last inequality we have used (\ref{eq:choice_r_n}). So the first term on the RHS of (\ref{eq:dist_decomp}) tends to $0$ as $n \to \infty$, and so we are left to show
\begin{equation}\label{eq:typ_dist_final_target}
    \Prob{B^{F^{\mathbb{G}_n}_{U_n}(r_n)}_{U_n}(K_n) \neq B^{\mathbb{G}_n}_{U_n}(K_n)} \to 0,
\end{equation}
as $n \to \infty$, where $U_n$ is uniformly distributed on $V(\mathbb{G}_n)=[n]$.

We now recall the general estimate from Remark \ref{rem:gen_estimate}:
\begin{align*}
        & \Prob{B^{F^{\mathbb{G}_n}_{U_n}(r_n)}_{U_n}(K_n) \neq B^{\mathbb{G}_n}_{U_n}(K_n)} \\
    &  \leq \sum_{j=1}^{K_n} \int_{\mathscr{B}^{a_n^m}_{\mathbf{0}}}\cdots\int_{\mathscr{B}^{a_n^{m^{j-1}}}_{\mathbf{0}}}\int_{\mathbb{R}^d \setminus \mathscr{B}^{a_n^{m^j}}_{\mathbf{0}}} \Exp{\kappa_n(\|z_j\|,W^{(n)}_{{U}'_{n,1}},W^{(n)}_{{U}'_{n,2}})} dz_j\cdots dz_1, \numberthis \label{eq:dist_gen_bd.}     
\end{align*}
where $a=a_n$ is as in (\ref{eq:dist_choice_an}), $m$ is as in (\ref{eq:choice_m}), $K=K_n$ is as in (\ref{eq:dist_choice_Kn}), $r_n=r(a_n,m,K_n)$ is as in (\ref{eq:choice_of_r}), and ${U}'_{n,1}, {U}'_{n,2}$ are two independent uniformly distributed random variables on $[n]$.

 {Using (\ref{eq:dist_gen_bd.}) and Assumption \ref{sssec:assumption_connections} (3), taking $h_n\colon\mathbb{R}^2 \to \mathbb{R}$ in (\ref{eq:conv_bdd_func_weghts}) being equal to $h_n(s,t)=\kappa_n(\|z_j\|,s,t,)$,} we note that for $n$ sufficiently large such that $a_n^m > t_0$ and
 {\[
 \Exp{\kappa_n(\|z_j\|,W^{(n)}_{{U}'_{n,1}},W^{(n)}_{{U}'_{n,2}})} \leq \|z_j\|^{-\alpha},
 \]}
    \begin{align*}
     \Prob{B^{F^{\mathbb{G}_n}_{U_n}(r_n)}_{U_n}(K_n) \neq B^{\mathbb{G}_n}_{U_n}(K_n)} & \leq \sum_{j=1}^{K_n} \int_{\mathscr{B}^{a_n^m}_{\mathbf{0}}}\cdots\int_{\mathscr{B}^{a_n^{m^{j-1}}}_{\mathbf{0}}}\int_{\mathbb{R}^d \setminus \mathscr{B}^{a_n^{m^j}}_{\mathbf{0}}} \frac{1}{\|z_j\|^{\alpha}} dz_j\cdots dz_1.\\
     &= \sum_{j=1}^{K_n} \frac{w^{j}}{a_n^{dm^j\left(\frac{\alpha}{d}-1-\frac{1}{m}-\cdots-\frac{1}{m^{j-1}} \right)}}, \numberthis \label{eq:dist_UB_error_1}
\end{align*}
where $w>1$ is some constant upper bound on $\lambda_d(\mathscr{B}^1_{\mathbf{0}})$.

Since $\frac{\alpha}{d}> \frac{m}{m-1}$, and if we let $C_0=d \left(\frac{\alpha}{d}- \frac{m}{m-1}\right)>0$, then from (\ref{eq:dist_UB_error_1})  {it follows that for all $n$ large enough} such that $\frac{w}{a_n^{C_0}}<1$, {and for some sufficiently large $J > 1$ such that for all $j>J$ we have $j^{1/j}<m$ (note that such a $J$ exists since $m>1$), we can write (assume $n$ is large so that $K_n>J+1$)}

\begin{align*}
         \Prob{B^{F^{\mathbb{G}_n}_{U_n}(r_n)}_{U_n}(K_n) \neq B^{\mathbb{G}_n}_{U_n}(K_n)} \leq  \sum_{j=1}^{K_n} \frac{w^{j}}{a_n^{C_0m^j}} {=\sum_{j=1}^{J} \frac{w^{j}}{a_n^{C_0m^j}}}{+\sum_{j=J+1}^{K_n}\frac{w^{j}}{a_n^{C_0m^j}}.} \numberthis \label{eq:dist_UB_error_2}
\end{align*}
{Note that the first term on the RHS of (\ref{eq:dist_UB_error_2}) clearly converge to $0$ as $n \to \infty$ since $a_n \to \infty$. For the second term on the RHS of (\ref{eq:dist_UB_error_2}), we note that since for all $j \geq J+1$, $m^j>j$, and since $\frac{w}{a_n^{C_0}}<1$, we have}
\begin{align*}
    {\sum_{j=J+1}^{K_n}\frac{w^{j}}{a_n^{C_0m^j}}}
     {\leq \sum_{j=J+1}^{K_n}\frac{w^{m^j}}{a_n^{C_0m^j}}}
     {\leq \sum_{j=J+1}^{K_n}\left(\frac{w}{a_n^{C_0}}\right)^{j}}
     {\leq \sum_{j=1}^{\infty}\left(\frac{w}{a_n^{C_0}}\right)^{j}}
     {=\frac{w}{a_n^{C_0}} \frac{1}{1-\frac{w}{a_n^{C_0}}} \to 0,}
\end{align*}
{as $n \to \infty$.}
{This} implies
(\ref{eq:typ_dist_final_target}), and completes the proof of Theorem \ref{thm:main_typical_dist}.

\end{proof}

\subsection{Proofs of results on examples}\label{ssec:proof_examples}
\begin{proof}[Proof of Lemma \ref{lem:PSIRG_regular_variation}]

We write
\begin{align*}
    &\hspace{-30pt}\Exp{\kappa(t,W^{(1)},W^{(2)})} \\
    &= \Exp{1 \wedge f(t)g(W^{(1)},W^{(2)})}\\
    &=  {\Prob{g(W^{(1)}, W^{(2)})> 1/f(t)}}+ f(t) \Exp{g(W^{(1)}W^{(2)})\ind{g(W^{(1)}W^{(2)})<1/f(t)}}. \numberthis \label{eq:PSIRG_exp_decomp_1}
\end{align*}

Let us denote $Y= {g(W^{(1)},W^{(2)})}\ind{g(W^{(1)}W^{(2)})<1/f(t)}$. Note that since $Y$ is a non-negative random variable, $\Exp{Y}=\int_{0}^{\infty} \Prob{Y \geq l} dl$. We note that

\begin{equation}\label{UB_Ygl_PSIRG}
    \begin{split}
        \Prob{Y\geq l} \leq 
        \begin{cases}
            0 & \text{if}\;l > 1/f(t),
            \\ \Prob{ {g(W^{(1)},W^{(2)})} > l} & \text{if}\; l \leq 1/f(t).
        \end{cases}
    \end{split}
\end{equation}
Hence, 
\begin{equation}\label{eq:expec_Y_PSIRG_UB}
\begin{split}
    \Exp{g(W^{(1)}W^{(2)})\ind{g(W^{(1)},W^{(2)})<1/f(t)}}&=\Exp{Y}=\int_{\mathbb{R}_+} \Prob{Y \geq l}dl\\
    & \leq \int_{0}^{1/f(t)} \Prob{g(W^{(1)}W^{(2)}) > l} dl.
\end{split}
\end{equation}

{Recall $t_1$ and $t_2$ from Assumption \ref{sssec:assump_PSIRG_connections}.} Let
\[
\overline{t}_1:=\inf\{t'>0:t>t'\implies f(t)^{-1}>t_2\}.
\]
Note from Assumption \ref{sssec:assump_PSIRG_connections} (2)  {that}, since $f(t)^{-1}$ increases to $\infty$ as $t \to \infty$, $\overline{t}_1$ is well defined.

Let 
\[
t_0:=\max\{\overline{t}_1,t_1\}.
\]
Then for any $t>t_0$, we bound
\begin{align*}
    \int_{0}^{1/f(t)} \Prob{g(W^{(1)}W^{(2)}) > l} dl
    &\leq t_2+\int_{t_2}^{f(t)^{-1}} l^{- \beta_p} dl\\
    &=t_2+\frac{1}{1-\beta_p}\left(f(t)^{\beta_p-1}-t_2^{1-\beta_p}\right).
\end{align*}
Hence from (\ref{eq:PSIRG_exp_decomp_1}), we have for $t>t_0$, 
\begin{align*}
    \Exp{\kappa(t,W^{(1)},W^{(2)})}
    &\leq f(t)^{\beta_p}+t_2f(t)+\frac{1}{1-\beta_p}\left(f(t)^{\beta_p}-t_2^{1-\beta_p}f(t)\right)\\
    &\leq t^{-\alpha_p\beta_p}+t_2t^{-\alpha_p}+\frac{1}{1-\beta_p}\left(t^{-\alpha_p\beta_p}-t_2^{1-\beta_p}t^{-\alpha_p}\right).\numberthis \label{eq:PSIRG_reg_var_last_step}
\end{align*}
Since {for any $\epsilon>0$,} the RHS of (\ref{eq:PSIRG_reg_var_last_step}) is dominated by $t^{-\min\{\alpha_p,\alpha_p\beta_p\}+\epsilon}$ outside a compact set {(depending on $\epsilon$)}, we are done.
\end{proof}

\begin{proof}[Proof of Corollary \ref{cor:LWP_GIRGs}]
We will rely upon Corollary \ref{cor:PSIRG} and Remark \ref{rem:dominance_PSIRG_connections} afterwards to conclude the proof. 

Note that by definition, $\mathrm{GIRG}_{n,\alpha_G,\beta_G,d}$ is the SIRG $G(\mathbf{X}^{(n)},\mathbf{W}^{(n)},\kappa_n^{\alpha_G})$, where the locations $(\mathbf{X}^{(n)})_{i \in [n]}$ satisfy Assumption \ref{sssec:assumption_location}, the weights $(\mathbf{W}^{(n)}_i)_{i \in [n]}$ satisfy Assumption \ref{sssec:assumption_weights}, and where $\kappa_n^{\alpha_G}:\mathbb{R}_+ \times \mathbb{R} \times \mathbb{R} \to [0,1]$ is defined as
\begin{equation}\label{GIRG_connection_fn_kappa_n}
    \begin{split}
        \kappa_n^{\alpha_G}(t,x,y):=
        \begin{cases}
            1 \wedge \left(\frac{xy}{\sum_{i\in[n]} W^{(n)}_i}\right)^{\alpha_G}\frac{1}{t^{d\alpha_G}} , &\text{if}\; 1<\alpha_G<\infty; \\
            \ind{\left(\frac{xy}{\sum_{i\in[n]} W^{(n)}_i}\right)^{1/d}>t}, & \text{if}\; \alpha_G = \infty.
        \end{cases}
    \end{split}
\end{equation}
It is not very difficult to check that $\kappa^{\alpha_G}_n$ satisfies Assumption \ref{sssec:assumption_connections} (1) with limiting connection function $\kappa^{\alpha_G}$ as defined in (\ref{eq:GIRG_limiting_connection_fn}), hence we only need to check Assumption \ref{sssec:assump_PSIRG_connections} for $\kappa^{(\alpha_G)}$ to directly apply Corollary \ref{cor:PSIRG}.

\paragraph{Case 1: $\alpha_G< \infty$.}
Note that for the $\alpha<\infty$ case, $\kappa^{\alpha_G}$ is a PSIRG connection function with $g(x,y)=\left(\frac{xy}{\Exp{W}}\right)^{\alpha_G}$, and $f(t)={{t^{-d\alpha_G}}}$. 

Since $W$ has a power-law tail with exponent $\beta_G-1$, using Brieman's Lemma \cite[Lemma 1.4.3]{Kulik_Soulier_2020}, for any $\epsilon>0$, the tail $\Prob{W_1W_2>t}$ of the product of two i.i.d.\ copies $W_1$ and $W_2$ of $W$  is dominated from above by a regularly varying function with exponent $\beta_G-1-\epsilon$ for all sufficiently large $t$. Hence if we choose $\epsilon>0$ sufficiently small such that $\beta_G-1-\epsilon>1$, we note that $g(x,y)$ satisfies Assumption \ref{sssec:assump_PSIRG_connections} (3) with $\beta_p=(\beta_G-1-\epsilon)/\alpha_G$. Also, clearly $f(t)$ satisfies Assumption \ref{sssec:assump_PSIRG_connections} (2) with $\alpha_p=d\alpha_G$. Hence in this case $\gamma_p=\min\{\alpha_p,\alpha_p\beta_p\}=\min\{d\alpha_G,d(\beta_G-1-\epsilon)\}>d$, since both $\alpha_G,(\beta_G-1-\epsilon)>1$. So, we can conclude the result in this case using Corollary \ref{cor:PSIRG}.

\paragraph{Case 2: $\alpha_G = \infty$.} Fix $\gamma>d$. When $\alpha_G=\infty$, we note from (\ref{eq:GIRG_limiting_connection_fn}) that the function $\kappa^{(\infty)}(t,x,y)$ can be bounded from above as
\[
    \kappa^{(\infty)}(t,x,y)= \ind{\left(\frac{xy}{\Exp{W}}\right)^{\gamma/d}>t^{\gamma}} \leq 1 \wedge \frac{\left(\frac{xy}{\Exp{W}}\right)^{\gamma/d}}{t^{\gamma}} =: h(t,x,y),
\]
and clearly $h(t,x,y)$ is a PSIRG connection function with $f(t)=\frac{1}{t^{\gamma}}$ satisfying Assumption \ref{sssec:assump_PSIRG_connections} (2) with $\alpha_p=\gamma$, and $g(x,y)=\left(\frac{xy}{\Exp{W}}\right)^{\gamma/d}$ satisfying Assumption \ref{sssec:assump_PSIRG_connections} (3) with $\beta_p= d(\beta_G-1-\epsilon)/\gamma$, for $\epsilon>0$ sufficiently small such that $\beta_G-1-\epsilon>1$, using again Brieman's Lemma \cite[Lemma 1.4.3]{Kulik_Soulier_2020}. Since in this case also we have $\gamma_p=\min\{\alpha_p,\alpha_p\beta_p\}=\min\{\gamma,d(\beta_G-1-\epsilon)\}>d$, we can conclude the proof using Remark \ref{rem:dominance_PSIRG_connections}.

\end{proof}

\begin{proof}[{Proof of Corollary \ref{cor:LWP_HRGs}}]


{We first transform the Hyperbolic Random Graph models into $1$-dimensional SIRGs with appropriate parameters. To do this, we follow the proof of \cite[Theorem 9.6]{JK_BL_explosions_HRGs_19}. Recall from Section \ref{sssec:examples_HRG} respectively the radial component vector $(r_i^{(n)})_{i=1}^n$ and the angular component vector $(\theta_i^{(n)})_{i=1}^n$ of the vertices $(u_i^{(n)})_{i=1}^n$ of the THRG and PHRG models. Consider the transformations 
\begin{align*}
    {X_i^{(n)}=\mathscr{X}(\theta_i^{(n)}):= \frac{\theta_i^{(n)}}{2\pi};\;\; W_i^{(n)}=\mathscr{W}(r_i^{(n)}):= \exp{\frac{R_n-r_i^{(n)}}{2}}.} \numberthis \label{eq:HRG_transform_1_GIRG}
\end{align*}
Clearly, $(X_i^{(n)})_{i=1}^n$ is then a vector with i.i.d.\ coordinates on $[-1/2,1/2]$, and using (\ref{eq:radial_CDF_HRG}), it can be shown the i.i.d.\ components of the vector $(W_i^{(n)})_{i=1}^n$ have a power-law distribution with parameter $2\alpha_H+1$, when $\alpha_H>1/2$ (see \cite[(9.8)]{JK_BL_explosions_HRGs_19}, and the text following it). 

Recall the connection functions $p^{(n)}_{\mathrm{THRG}}$ and $p^{(n)}_{\mathrm{PHRG}}$ respectively from (\ref{eq:THRG_conn}) and (\ref{eq:PHRG_conn}). The hyperbolic distance $d_{\mathbb{H}}(u_i^{(n)},u_j^{(n)})$ between $u_i^{(n)}=(r_i^{(n)},\theta_i^{(n)})$ and $u_j^{(n)}=(r_j^{(n)},\theta_j^{(n)})$ depends on the angular coordinates $\theta_i^{(n)}$ and $\theta_j^{(n)}$ through $\cos(\theta_i^{(n)}-\theta_j^{(n)})$ (see \cite[(9.1)]{JK_BL_explosions_HRGs_19}). Hence, it can be seen as a function of $|\theta_i^{(n)}-\theta_j^{(n)}|$ since $\cos(\cdot)$ is symmetric. Consequently, there exist functions $\overline{p}_{\mathrm{THRG}}$ and $\overline{p}_{\mathrm{PHRG}}$ such that
\begin{align*}
    &p^{(n)}_{\mathrm{THRG}}(u^{(n)}_i,u^{(n)}_j)=\overline{p}^{(n)}_{\mathrm{THRG} }(|\theta_i^{(n)}-\theta_j^{(n)}|,r_i^{(n)},r_j^{(n)}),\\
    &p^{(n)}_{\mathrm{PHRG}}(u^{(n)}_i,u^{(n)}_j)=\overline{p}^{(n)}_{\mathrm{PHRG} }(|\theta_i^{(n)}-\theta_j^{(n)}|,r_i^{(n)},r_j^{(n)}). \numberthis \label{eq:HRG_conn_fN_rad_ang}
\end{align*}
Writing
\begin{align*}
    &\kappa_{\mathrm{THRG},n}\left(t,x,y\right)=\overline{p}^{(n)}_{\mathrm{THRG} }(2\pi t,g_n(x),g_n(y)),\\
    &\kappa_{\mathrm{PHRG},n}\left(t,x,y\right)=\overline{p}^{(n)}_{\mathrm{PHRG} }(2\pi t,g_n(x),g_n(y)),\numberthis \label{eq:transformed_HRG_connections}
\end{align*}
where the function $g_n$ satisfies
\begin{align*}
    g_n(x) = R_n - 2\log(x) 
\end{align*}
it was shown in \cite[(9.17)]{JK_BL_explosions_HRGs_19} and \cite[(9.16)]{JK_BL_explosions_HRGs_19} respectively, that for fixed $(t,x,y)$ 
as $n \to \infty$,
\begin{align*}
    &\kappa_{\mathrm{THRG},n}\left(t,x,y\right)=\overline{p}^{(n)}_{\mathrm{THRG} }(2\pi t,g_n(x),g_n(y)) \to \kappa_{\mathrm{THRG},\infty}(t,x,y)\\
    &\kappa_{\mathrm{PHRG},n}\left(t,x,y\right)=\overline{p}^{(n)}_{\mathrm{PHRG} }(2\pi t,g_n(x),g_n(y)) \to \kappa_{\mathrm{THRG},\infty}(t,x,y),
\end{align*}
where 
\begin{align*}
    \kappa_{\mathrm{THRG},\infty}(t,x,y):= \ind{t \leq \frac{\nu xy}{\pi}}; \;\; \kappa_{\mathrm{PHRG},\infty}(t,x,y):= \left(1+\left(\frac{ \pi t}{\nu xy} \right)^{1/T_H} \right)^{-1}. \numberthis \label{eq:HRG_limiting_conn_fn}
\end{align*}
In particular, the THRG and the PHRG models can be seen as finite $1$-dimensional SIRGs, with the vertex locations $(X_i^{(n)})_{i=1}^n$ satisfying Assumption \ref{sssec:assumption_location}, the vertex weights $(W_i^{(n)})_{i=1}^n$ satisfying Assumption \ref{sssec:assumption_weights} with the function $F_{W}(x)$ being a power-law distribution function with exponent $2\alpha_H+1$, and with the connection functions as in (\ref{eq:transformed_HRG_connections}), converging pointwise to limiting connection functions as in (\ref{eq:HRG_limiting_conn_fn}). The pointwise convergence can in fact be improved to the case where one has sequences $x_n\to x,y_n\to y$ as in Assumption \ref{sssec:assumption_connections} (1). This is because the error terms are uniformly bounded, see \cite[(9.15)]{JK_BL_explosions_HRGs_19}, which implies $\kappa_{\mathrm{THRG},n}(t,x_n,y_n)\to \kappa_{\mathrm{THRG},\infty}(t,x,y)$ and $\kappa_{\mathrm{PHRG},n}(t,x_n,y_n)\to \kappa_{\mathrm{PHRG},\infty}(t,x,y)$, with $t$ avoiding a set of measure zero for the THRG case, namely the set $\{\frac{\nu xy}{\pi}\}$. So, the sequence of connection functions $\kappa_{\mathrm{THRG},n}$ and $\kappa_{\mathrm{PHRG},n}$ satisfy Assumption \ref{sssec:assumption_connections} (1) with limiting connection functions (\ref{eq:HRG_limiting_conn_fn}).

Finally, we need to check the limiting connection functions $\kappa_{\mathrm{THRG},n}$ and $\kappa_{\mathrm{PHRG},n}$ satisfy Assumption \ref{sssec:assumption_connections} (2) with some $\alpha>d$. For this, we use Corollary \ref{cor:PSIRG} and Remark \ref{rem:dominance_PSIRG_connections}.

\paragraph{Case 1: THRG.}
Let $\gamma>1$ be any constant. Note that the function $\kappa_{\mathrm{THRG},\infty}$ can be bounded from above as
\begin{align*}
    \kappa_{\mathrm{THRG},\infty}(t,x,y) \leq 1 \wedge \frac{\left(\frac{\nu xy}{\pi}\right)^{\gamma}}{t^{\gamma}}.
\end{align*}
Note that this is a PSIRG connection function with $f(t)=\frac{1}{t^{\gamma}}$ satisfying Assumption \ref{sssec:assump_PSIRG_connections} (1) with $\alpha_p=\gamma$, and $g(x,y)=\left(\frac{\nu xy}{\pi}\right)^{\gamma}$. In addition, from \eqref{eq:HRG_transform_1_GIRG} it follows that the limiting weights are i.i.d. and have a power-law distribution with exponent $2\alpha_H+1$. Hence by using Brieman's Lemma \cite[Lemma 1.4.3]{Kulik_Soulier_2020}, if $W^{(1)}$ and $W^{(2)}$ are i.i.d.\ copies of the limiting weight distribution, $g(W^{(1)},W^{(2)})$ is regularly varying with exponent $2 \alpha_H/\gamma$. Applying Potter's bounds we conclude that $g(W_1,W_2)$ satisfies Assumption \ref{sssec:assump_PSIRG_connections} (2) with $\beta_p= (2 \alpha_H-\epsilon)/\gamma$, for $\epsilon>0$ sufficiently small such that $2\alpha_H-\epsilon>1$. 

Since in this case we have $\gamma_p=\min\{\alpha_p,\alpha_p\beta_p\}=\min\{\gamma,(2\alpha_H-\epsilon)\}>1$, we can conclude the proof using Corollary \ref{cor:PSIRG} and Remark \ref{rem:dominance_PSIRG_connections}.


\paragraph{Case 2: PHRG.}
The function $\kappa_{\mathrm{PHRG},\infty}$ can be bounded from above as,
\begin{align*}
    \kappa_{\mathrm{PHRG},\infty}(t,x,y)\leq C_1\left(1 \wedge a_1\left(\frac{xy}{t}\right)^{1/T_H}\right), \numberthis \label{eq:PHRG_conn_UB_PSIRG_1}
\end{align*}
for some constants $C_1,a_1>0$. To see this, combine \cite[(9.14)]{JK_BL_explosions_HRGs_19}, \cite[(9.16)]{JK_BL_explosions_HRGs_19} and Assumption \ref{sssec:assumption_connections} (1). 
Using the fact that $\kappa_{\mathrm{PHRG},\infty}(t,x,y)$ is a probability, and hence $\leq 1$, we can further get an upper bound from (\ref{eq:PHRG_conn_UB_PSIRG_1}) as,
\begin{align*}
    \kappa_{\mathrm{PHRG},\infty}(t,x,y)\leq 1 \wedge C_1a_1\left(\frac{xy}{t}\right)^{1/T_H}.
\end{align*}
Note that this is a PSIRG connection function, with $f(t)=\frac{C_1a_1}{t^{1/T_H}}$, and $g(x,y)=\left(xy\right)^{1/T_H}$. Recall from the statement of Corollary \ref{cor:LWP_HRGs} that $0<T_H<1$. Note that $f(t)$ then satisfies Assumption \ref{sssec:assump_PSIRG_connections} (1) with $\alpha_p=1/(T_H-\epsilon_1)$, for some $\epsilon_1>0$ sufficiently small such that $1/T_H-\epsilon_1>1$. Also note that from (\ref{eq:HRG_transform_1_GIRG}), the limiting weights are i.i.d.\ and have a power-law distribution with exponent $2\alpha_H+1$. So if we let $W^{(1)}$ and $W^{(2)}$ to be i.i.d.\ copies of the limiting weight distribution, using Brieman's Lemma \cite[Lemma 1.4.3]{Kulik_Soulier_2020}, $g(W^{(1)},W^{(2)})$ is regularly varying with exponent $2\alpha_HT_H$. Hence, applying Potter's bounds, we note that $g(W^{(1)},W^{(2)})$  satisfies Assumption \ref{sssec:assump_PSIRG_connections} (2) with $\beta_p= (2\alpha_H-\epsilon_2)T_H$, for some $\epsilon_2>0$ sufficiently small such that $2\alpha_H-\epsilon_2>1$. Since in this case also we have $\gamma_p=\min\{\alpha_p,\alpha_p\beta_p\}=\min\{1/(T_H-\epsilon_1),(2\alpha_H-\epsilon_2)T_H/(T_H-\epsilon_1)\}>1$, we can conclude the proof using Corollary \ref{cor:PSIRG} and Remark \ref{rem:dominance_PSIRG_connections}.
 
}

\end{proof}

\begin{proof}[Proof of Corollary \ref{cor:LCinP_CSFP}]

We apply Corollary \ref{cor:PSIRG}, using Remark \ref{rem:dominance_PSIRG_connections}. Note that Assumptions \ref{sssec:assumption_location} and \ref{sssec:assumption_weights} are immediate. Since Assumption \ref{sssec:assumption_connections} (1) is immediate for $\kappa_n$ with limit $\kappa$, we only need to check that $\kappa$ is dominated by a PSIRG connection function which satisfies Assumption \ref{sssec:assump_PSIRG_connections}.

We use the easy bound 
\[
1-\exp\left({-\frac{\lambda xy}{t^{\alpha}}}\right) \leq 1 \wedge \frac{\lambda xy}{t^{\alpha}}
\]
to observe that the limiting connection function $\kappa$ is dominated by the PSIRG connection function $1 \wedge f(t)g(x,y)$, where $f(t)=1/t^{\alpha}$ and $g(x,y)=\lambda xy$. 

Note that for i.i.d.\ copies $W_1$ and $W_2$ of the weight distribution $W$ as in (\ref{eq:wt_law_CSFP}), for any $\epsilon>0$, the tail $\Prob{g(W_1,W_2)>t}$ of the random variable $g(W_1,W_2)$ is dominated by a regularly varying function with exponent $\beta-\epsilon>0$ by Breiman's Lemma \cite[Lemma 1.4.3]{Kulik_Soulier_2020}, which implies  {that} $g$ satisfies Assumption \ref{sssec:assump_PSIRG_connections} (2) with $\beta_p=\beta-\epsilon$. Also, clearly $f$ satisfies Assumption \ref{sssec:assump_PSIRG_connections} (1) with $\alpha_p = \alpha$. For $\epsilon>0$ sufficiently small such that $\beta-\epsilon>1$, since we have $\gamma_p=\min\{\alpha_p\beta_p,\alpha_p\}=\min\{\alpha(\beta-\epsilon),\alpha\}>d$, the proof of Corollary \ref{cor:LCinP_CSFP} is complete using Remark \ref{rem:dominance_PSIRG_connections} and Corollary \ref{cor:PSIRG}.    
\end{proof}

\subsection{Proofs of degree results}\label{ssec:proof_degree}

\begin{proof}[Proof of Proposition \ref{prop:degree_law}]
We argue by showing the moment generating function converges to the moment generating function of the claimed limit. Also, we write interchangeably the vertex set $V(\mathbb{G}_{\infty})$, and $\mathbb{N} \cup \{0\}$. In particular, the set of all vertices of $\mathbb{G}_{\infty}$ other than $0$, is $\mathbb{N}$.  

Note that we have to show that, for all $t \in \mathbb{R}$, 
\begin{equation}\label{eq:deg_law_target}
    \CExp{e^{tD}}{W_0} \stackrel{\text{a.s.}}{=} \exp{\left((e^t-1)\int_{\mathbb{R}^d} \CExp{\kappa(\|z\|,W_0,W^{(1)})}{W_0}dz\right)}.
\end{equation}

Recall that $E(\mathbb{G}_{\infty})$ is the edge set of $\mathbb{G}_{\infty}$. For any $r>0$, we define $D_{\leq r}$ as
\begin{equation}\label{eq:defn_Dlr}
    D_{\leq r}:= \sum_{i \in \mathbb{N}}\ind{\{0,i\}\in E(\mathbb{G}_{\infty})}\ind{\|Y_i\| \leq r}.
\end{equation}

Clearly, $D_{\leq r}$ increases to $D$ as $r \to \infty$. So applying conditional monotone convergence, when $t \geq 0$, and conditional dominated convergence, when $t<0$,  for any $t \in \mathbb{R}$,
\begin{equation}\label{eq:deg_law_lim_r}
    \lim_{r \to \infty} \CExp{e^{t D_{\leq r}}}{W_0} \stackrel{\text{a.s.}}{=} \CExp{e^{tD}}{W_0}.
\end{equation}

Let $\mathscr{B}^r_{\mathbf{0}}=\{y \in \mathbb{R}^d:\|y\|<r\}$ denote the open Euclidean ball of radius $r$ in $\mathbb{R}^d$, centered at $\mathbf{0}$.

Let $Q_r$ be a $\text{Poi}\left(\lambda_d(\mathscr{B}^r_{\mathbf{0}})\right)$ random variable, and conditionally on $Q_r$, let
\begin{itemize}
    \item[a.] $\{\mathcal{R}_i\}_{i=1}^{Q_r}$ be a collection of $Q_r$ i.i.d.\ uniform random variables on $\mathscr{B}^r_{\mathbf{0}}$;
    \item[b.] $\{\mathcal{U}_i\}_{i=1}^{Q_r}$ be a collection of $Q_r$ i.i.d.\ uniform random variables on $[0,1]$;
    \item[c.] $\{W^{(i)}\}_{i=1}^{Q_r}$ be $Q_r$ i.i.d.\ copies of $W_0$, independent of $\{\mathcal{R}_i\}_{i=1}^{Q_r}$, $\{\mathcal{U}_i\}_{i=1}^{Q_r}$ and $W_0$. 
\end{itemize}

Note then that 
\begin{equation}\label{eq:decomp_law_Dlr}
    e^{t D_{\leq r}}\Big|W_0 \stackrel{d}{=} \prod_{i=1}^{Q_r} \exp{\left(t\ind{\mathcal{U}_i\leq \kappa(\|\mathcal{R}_i\|,W_0,W^{(i)})}\right)}\Big|W_0,
\end{equation}
and observe that the product on the RHS of (\ref{eq:decomp_law_Dlr}) is, conditionally on $W_0$, a product of  {(conditionally)} independent random variables.

We compute that
\begin{align}
    &\CExp{\exp{\left(t\ind{\mathcal{U}_i\leq \kappa(\|\mathcal{R}_i\|,W_0,W^{(i)})}\right)}}{W_0,W^{(i)},\mathcal{R}_i}\nonumber\\
    &\qquad \stackrel{\text{a.s.}}{=} 1-\kappa(\|\mathcal{R}_i\|,W_0,W^{(i)})+e^t \kappa(\|\mathcal{R}_i\|,W_0,W^{(i)}),
    \label{eq:deg_law_1}
\end{align}
so that
\begin{align*}
    &\CExp{\exp{\left(t\ind{\mathcal{U}_i\leq \kappa(\|\mathcal{R}_i\|,W_0,W^{(i)})}\right)}}{W_0} \\
    & \stackrel{\text{a.s.}}{=} 1-\frac{1}{\lambda_d(\mathscr{B}^r_{\mathbf{0}})}\int_{\mathscr{B}^r_{\mathbf{0}}}\CExp{\kappa(\|z\|,W_0,W^{(i)})}{W_0}dz+e^t \frac{1}{\lambda_d(\mathscr{B}^r_{\mathbf{0}})}\int_{\mathscr{B}^r_{\mathbf{0}}}\CExp{\kappa(\|z\|,W_0,W^{(i)})}{W_0}dz\\
    & \stackrel{\text{a.s.}}{=} (e^t-1)\frac{1}{\lambda_d(\mathscr{B}^r_{\mathbf{0}})}\int_{\mathscr{B}^r_{\mathbf{0}}}\CExp{\kappa(\|z\|,W_0,W^{(1)})}{W_0}dz+1.
    \numberthis \label{eq:deg_law_2}
\end{align*}
Hence, by (\ref{eq:decomp_law_Dlr}),
\begin{align*}
    &\CExp{e^{tD_{\leq r}}}{W_0} \\
    & \stackrel{\text{a.s.}}{=} \CExp{\left((e^t-1)\frac{1}{\lambda_d(\mathscr{B}^r_{\mathbf{0}})}\int_{\mathscr{B}^r_{\mathbf{0}}}\CExp{\kappa(\|z\|,W_0,W^{(1)})}{W_0}dz+1\right)^{Q_r}}{W_0}\\
    & \stackrel{\text{a.s.}}{=} \exp{\left((e^t-1)\int_{\mathscr{B}^r_{\mathbf{0}}} \CExp{\kappa(\|z\|,W_0,W^{(1)})}{W_0}dz\right)},
    \numberthis \label{eq:deg_law_3}
\end{align*}
 {since $Q_r$ has a Poisson distribution with parameter $\lambda_d(\mathscr{B}^r_{\mathbf{0}})$.}
Now, we let $r \to \infty$ in both sides of (\ref{eq:deg_law_3}), and use (\ref{eq:deg_law_lim_r}), to establish (\ref{eq:deg_law_target}) and conclude the proof of Proposition \ref{prop:degree_law}.
\end{proof}

\begin{proof}[Proof of Proposition \ref{prop:UI_typical_degs}]
Recall that $D_n$ is the degree of the uniformly chosen vertex $U_n$ of $\mathbb{G}_n$. Fix $\varepsilon>0$, and note that the target is to show that there  {exist} $M_0,N \in \mathbb{N}$, such that for all $M>M_0$ and $n> N=N(M_0)$, 
\begin{equation}\label{eq:deg_UI_target}
    \Exp{D_n \ind{D_n>M}}<\varepsilon.
\end{equation}

For any $r>0$, we can write 
\begin{equation}\label{eq:D_n_decomp}
    D_n=D_{n,<r}+D_{n,\geq r},
\end{equation}
where (recall  {that} $E(\mathbb{G}_n)$ is the edge set of $\mathbb{G}_n$)
\eqan{
\label{eq:defn_D_nlr}
    D_{n,<r}&:= \sum_{j \in [n]}\ind{\{j,U_n\}\in E(\mathbb{G}_n),\|Y^{(n)}_{U_n}-Y^{(n)}_j\|<r},\\
\label{eq:defn_D_ngr}
    D_{n,\geq r}&:= \sum_{j \in [n]}\ind{\{j,U_n\}\in E(\mathbb{G}_n),\|Y^{(n)}_{U_n}-Y^{(n)}_j\|\geq r}.
}

 {By applying the case $j=1$ of (\ref{eq:lem_bad_path_nj}), it} is not hard to see that we can choose and fix $r_0=r_0(\varepsilon)>0$, and $n_0 \in \mathbb{N}$, such that, whenever $r\geq r_0$ and $n>n_0$,
\begin{equation}\label{eq:exp_Dgr_small}
    \Exp{D_{n,\geq r}}\leq \varepsilon/4.
\end{equation}

Splitting depending on whether $D_{n,<r_0}$ or $D_{n,\geq r_0}$ is larger, we obtain
\begin{align*}
    \Exp{D_n \ind{D_n>M}}&\leq 
    \Exp{D_n \ind{D_n>M}\big(\ind{D_{n,<r_0}\leq D_{n,\geq r_0}}
    +\ind{D_{n,<r_0}>D_{n,\geq r_0}}\big)}\\
    &\leq 2\Exp{D_{n,\geq r_0} \ind{D_{n,\geq r_0}>M/2}}+
    2\Exp{D_{n,<r_0} \ind{D_{n,<r_0}>M/2}}\\
    & \leq \varepsilon/2+2\Exp{D_{n,<r_0} \ind{D_{n,<r_0}>M/2}}, \numberthis \label{eq:final_term_UI}
\end{align*}
where in the last step we have used (\ref{eq:exp_Dgr_small}).
    
Now observe that $D_{n,<r_0}$ is stochastically dominated by $Y_n$, where $Y_n$  {is} a $\text{Bin}\left(n-1,\frac{\lambda_d(\mathscr{B}^{r_0}_{\mathbf{0}})}{n}\right)$ random variable,  {so that} $(Y_n)_{n \geq 1}$ is uniformly integrable. Hence, there  {exist} $M_0 \in \mathbb{N}$ and $n_1=n_1(M_0) \in \mathbb{N}$ such that, whenever $M>M_0$, $n>n_1$, 
\begin{equation}\label{eq:UI_SD}
    2\Exp{D_{n,< r_0}\ind{D_{n,< r_0}>M/2}}< \varepsilon/2.
\end{equation}
Hence using (\ref{eq:UI_SD}) and (\ref{eq:final_term_UI}), we note that (\ref{eq:deg_UI_target}) holds for $M>M_0$, and $n>N=\max\{n_1(M_0),n_0\}$.
\end{proof}

\subsection{Proofs of clustering results}\label{ssec:proof_clustering}

\begin{proof}[Proof of Corollary \ref{cor:conv_clustering}]

Parts (2) and (3) follow directly using Theorem \ref{thm:main_LCinP_SIRGs}, with \cite[Theorem 2.22]{RGCN_2} and \cite[Exercise 2.31]{RGCN_2} respectively.

For Part (1), using \cite[Theorem 2.21]{RGCN_2}, we only have to verify the uniformly integrability of $(D_n^2)_{n \geq 1}$ and that $\Prob{D>1}>0$, where, as before, $D_n$ is the degree of the uniformly chosen vertex $U_n$ of $\mathbb{G}_n$, and $D$ is the degree of $0$ in $\mathbb{G}_{\infty}$.
By Proposition \ref{prop:degree_law}, $\Prob{D>1}>0$ is trivial (we  {do not} focus on the pathological case where $\int_{\mathbb{R}^d} \CExp{\kappa(\|z\|,W_0,W^{(1)})}{W_0}dz=0$, in which case it is not hard to see that $\mathbb{G}_n$ is an empty graph). So we need only verify that $\alpha>2d$ (where $\alpha$ is as in Assumption \ref{sssec:assumption_connections} (3)) implies the  {uniformly integrability of the sequence $(D_n^2)_{n \geq 1}$}.

Fix $\varepsilon>0$. We want to show there is $M_0 \in \mathbb{N}$, and $N=N(M_0) \in \mathbb{N}$, such that whenever $M>M_0$ and $n>N$,
\begin{equation}\label{eq:target_Dn2_UI}
    \Exp{D_n^2 \ind{D_n^2>M}}< \varepsilon.
\end{equation}

Recall the decomposition (\ref{eq:D_n_decomp}). Note that
\begin{align*}
    D_{n,\geq r}^2
    &= 2\sum_{i,j \in [n], \; i<j} \ind{\{U_n,i\} \in E(\mathbb{G}_{\infty}),\|Y^{(n)}_{U_n}-Y^{(n)}_i\|\geq r}\ind{\{U_n,j\} \in E(\mathbb{G}_{\infty}),\|Y^{(n)}_{U_n}-Y^{(n)}_j\|\geq r}
    \\& \hspace{10 pt}+\sum_{i=1}^n\ind{\{U_n,i\} \in E(\mathbb{G}_{\infty}),\|Y^{(n)}_{U_n}-Y^{(n)}_i\|\geq r}. \numberthis \label{eq:D_ngr2}
\end{align*}
Taking expectations on both sides of (\ref{eq:D_ngr2}), and after applying some routine change of variables, we get the bound
\begin{align*}
    &\Exp{D_{n,\geq r}^2}
    \\& \leq \int_{\mathbb{R}^d \setminus \mathscr{B}^r_{\mathbf{0}}}\int_{\mathbb{R}^d \setminus \mathscr{B}^r_{\mathbf{0}}} \Exp{\kappa_n\left(\|x\|,W^{(n)}_{U_{n,1}},W^{(n)}_{U_{n,2}}\right)\kappa_n\left(\|y\|,W^{(n)}_{U_{n,1}},W^{(n)}_{U_{n,3}}\right)}dx dy
    \\& \hspace{10 pt} + \int_{\mathbb{R}^d \setminus \mathscr{B}^r_{\mathbf{0}}} \Exp{\kappa_n\left(\|z\|,W^{(n)}_{U_{n,1}},W^{(n)}_{U_{n,2}}\right)}dz,
\end{align*}
where $U_{n,1},U_{n,2},U_{n,3}$ are i.i.d. uniformly distributed random variables on $[n]$, and recall $(W^{(n)}_i)_{i\in [n]}$ is the weight sequence corresponding to the random graph $\mathbb{G}_n$.

Applying  {the Cauchy-Schwarz} inequality on the integrand of the first term, and noting that $\kappa_n \leq 1$, we have the bound
\begin{align*}
    &\Exp{D_{n,\geq r}^2}
    \\& \leq \left(\int_{\mathbb{R}^d \setminus \mathscr{B}^r_{\mathbf{0}}} \Exp{\kappa_n\left(\|x\|,W^{(n)}_{U_{n,1}},W^{(n)}_{U_{n,2}}\right)}^{1/2}dx\right) \left( \int_{\mathbb{R}^d \setminus \mathscr{B}^r_{\mathbf{0}}} \Exp{\kappa_n\left(\|y\|,W^{(n)}_{U_{n,1}},W^{(n)}_{U_{n,3}}\right)}^{1/2} dy \right)
    \\& \hspace{10 pt} + \int_{\mathbb{R}^d \setminus \mathscr{B}^r_{\mathbf{0}}} \Exp{\kappa_n\left(\|z\|,W^{(n)}_{U_{n,1}},W^{(n)}_{U_{n,2}}\right)}dz. \numberthis \label{eq:exp_Dngr2_int_UB}
\end{align*}
Since $\alpha>2d$, we can choose and fix $r_0=r_0(\varepsilon)>0$, such that there is $n_0=n_0(r_0) \in \mathbb{N}$, such that whenever $n>n_0$,  {by} (\ref{eq:exp_Dngr2_int_UB}),
\begin{equation}\label{exp_Dngr2_small}
    \Exp{D_{n,\geq r_0}^2}< \varepsilon/ {8}.
\end{equation}

Hence, by splitting according to which one of $D_{n,<r_0}$ and $D_{n,\geq r_0}$ is  {largest},
\begin{align*}
    & \Exp{D_n^2 \ind{D_n^2>M}}\\
    & = \Exp{D_n^2 \ind{D_n^2>M}\left(\ind{D_{n,<r_0}>D_{n,\geq r_0}}+\ind{D_{n,<r_0}\leq D_{n,\geq r_0}}\right)}\\
    & \leq \Exp{4D_{n,<r_0}^2 \ind{4 D_{n,<r_0}^2 >M}}+ {4}\Exp{D_{n,  {\geq} r_0}^2}\\
    & \leq \Exp{4D_{n,<r_0}^2 \ind{4 D_{n,<r_0}^2 >M}}+\varepsilon/2. \numberthis \label{eq:Dn2_UI_conc_1}
\end{align*}
As argued previously in the proof of Proposition \ref{prop:UI_typical_degs}, $D_{n,<r_0}$ is stochastically dominated by a $\text{Bin}\left(n-1, \frac{\lambda_d(\mathscr{B}^{r_0}_{\mathbf{0}})}{n}\right)$ distributed random variable. Hence by standard arguments, the sequence $\left(4D_{n,<r_0}^2\right)_{n \geq 1}$ is uniformly integrable. Hence there  {exist} $M_0 \in \mathbb{N}$, and $n_1=n_1(M_0) \in \mathbb{N}$, such that, whenever $M>M_0$ and $n>n_1$, the first term on the RHS of (\ref{eq:Dn2_UI_conc_1}) is smaller than $\varepsilon/2$.  {We conclude that} (\ref{eq:target_Dn2_UI}) holds whenever $M>M_0$, and $n>N(M_0)=\max\{n_0,n_1\}$. This finishes the proof of part (1).  
\end{proof}

\paragraph{\bf Acknowledgements.} The work of RvdH was supported in part by the Netherlands Organisation for Scientific Research (NWO) through Gravitation-grant {\sc NETWORKS}-024.002.003. NM thanks Joost Jorritsma and Suman Chakraborty for helpful discussions and pointers to the literature, and Martijn Gösgens for help with the pictures.{The authors thank the anonymous reviewers for their helpful comments which greatly improved the presentation of the paper, and specially for pointing out the confusions regarding the statement of Corollary \ref{cor:LWP_HRGs} and the error in Conjecture \ref{conj:UB_typ_dist} in the first version of the paper. We also thank Peter Mörters for his inputs on the reformulation of Conjecture \ref{conj:UB_typ_dist}.}

\bibliographystyle{abbrvnat}
\bibliography{references}

\begin{thebibliography}{33}
\providecommand{\natexlab}[1]{#1}
\providecommand{\url}[1]{\texttt{#1}}
\expandafter\ifx\csname urlstyle\endcsname\relax
  \providecommand{\doi}[1]{doi: #1}\else
  \providecommand{\doi}{doi: \begingroup \urlstyle{rm}\Url}\fi

\bibitem[spa(1973)]{spatial_friendship_nets}
The spatial character of friendship formation.
\newblock \emph{Environment and Behavior}, 5\penalty0 (1):\penalty0 43--65,
  1973.
\newblock \doi{10.1177/001391657300500103}.
\newblock URL \url{https://doi.org/10.1177/001391657300500103}.

\bibitem[Aldous and Steele(2004)]{AS04}
D.~Aldous and J.~M. Steele.
\newblock The objective method: probabilistic combinatorial optimization and
  local weak convergence.
\newblock In \emph{Probability on discrete structures}, volume 110 of
  \emph{Encyclopaedia Math. Sci.}, pages 1--72. Springer, Berlin, 2004.
\newblock \doi{10.1007/978-3-662-09444-0_1}.

\bibitem[Benjamini and Schramm(2001)]{BSc01}
I.~Benjamini and O.~Schramm.
\newblock Recurrence of distributional limits of finite planar graphs.
\newblock \emph{Electron. J. Probab.}, 6:\penalty0 no. 23, 13, 2001.
\newblock ISSN 1083-6489.
\newblock \doi{10.1214/EJP.v6-96}.

\bibitem[Benjamini et~al.(2011)Benjamini, Kesten, Peres, and
  Schramm]{Geom_USF_BKPS}
I.~Benjamini, H.~Kesten, Y.~Peres, and O.~Schramm.
\newblock Geometry of the uniform spanning forest: transitions in dimensions
  {$4,8,12,\dots$} [mr2123930].
\newblock In \emph{Selected works of {O}ded {S}chramm. {V}olume 1, 2}, Sel.
  Works Probab. Stat., pages 751--777. Springer, New York, 2011.
\newblock \doi{10.1007/978-1-4419-9675-6\_25}.
\newblock URL \url{https://doi.org/10.1007/978-1-4419-9675-6_25}.

\bibitem[Bollob\'{a}s et~al.(2007)Bollob\'{a}s, Janson, and Riordan]{IRGs_BRJ}
B.~Bollob\'{a}s, S.~Janson, and O.~Riordan.
\newblock The phase transition in inhomogeneous random graphs.
\newblock \emph{Random Structures Algorithms}, 31\penalty0 (1):\penalty0
  3--122, 2007.
\newblock ISSN 1042-9832.
\newblock \doi{10.1002/rsa.20168}.
\newblock URL \url{https://doi.org/10.1002/rsa.20168}.

\bibitem[Bringmann et~al.(2016)Bringmann, Keusch, and Lengler]{KB_Avg_dist_16}
K.~Bringmann, R.~Keusch, and J.~Lengler.
\newblock Average distance in a general class of scale-free networks with
  underlying geometry.
\newblock \emph{ArXiv}, abs/1602.05712, 2016.

\bibitem[Bringmann et~al.(2017)Bringmann, Keusch, and
  Lengler]{KB_SamplingGIRGs_17}
K.~Bringmann, R.~Keusch, and J.~Lengler.
\newblock Sampling geometric inhomogeneous random graphs in linear time.
\newblock In \emph{25th {E}uropean {S}ymposium on {A}lgorithms}, volume~87 of
  \emph{LIPIcs. Leibniz Int. Proc. Inform.}, pages Art. No. 20, 15. Schloss
  Dagstuhl. Leibniz-Zent. Inform., Wadern, 2017.

\bibitem[Bringmann et~al.(2019)Bringmann, Keusch, and Lengler]{KB_GIRGs_19}
K.~Bringmann, R.~Keusch, and J.~Lengler.
\newblock Geometric inhomogeneous random graphs.
\newblock \emph{Theoret. Comput. Sci.}, 760:\penalty0 35--54, 2019.
\newblock ISSN 0304-3975.
\newblock \doi{10.1016/j.tcs.2018.08.014}.
\newblock URL \url{https://doi.org/10.1016/j.tcs.2018.08.014}.

\bibitem[Chung and Lu(2002{\natexlab{a}})]{CL02_1}
F.~Chung and L.~Lu.
\newblock The average distances in random graphs with given expected degrees.
\newblock \emph{Proc. Natl. Acad. Sci. USA}, 99\penalty0 (25):\penalty0
  15879--15882, 2002{\natexlab{a}}.
\newblock ISSN 0027-8424.
\newblock \doi{10.1073/pnas.252631999}.
\newblock URL \url{https://doi.org/10.1073/pnas.252631999}.

\bibitem[Chung and Lu(2002{\natexlab{b}})]{CL02_2}
F.~Chung and L.~Lu.
\newblock Connected components in random graphs with given expected degree
  sequences.
\newblock \emph{Ann. Comb.}, 6\penalty0 (2):\penalty0 125--145,
  2002{\natexlab{b}}.
\newblock ISSN 0218-0006.
\newblock \doi{10.1007/PL00012580}.
\newblock URL \url{https://doi.org/10.1007/PL00012580}.

\bibitem[Dalmau and Salvi(2021)]{CSFP_clustering}
J.~Dalmau and M.~Salvi.
\newblock Scale-free percolation in continuous space: quenched degree and
  clustering coefficient.
\newblock \emph{J. Appl. Probab.}, 58\penalty0 (1):\penalty0 106--127, 2021.
\newblock ISSN 0021-9002.
\newblock \doi{10.1017/jpr.2020.76}.
\newblock URL \url{https://doi.org/10.1017/jpr.2020.76}.

\bibitem[Deijfen et~al.(2013)Deijfen, {\swap{Hofstad}{~van~der~}}, and
  Hooghiemstra]{SFP}
M.~Deijfen, R.~{\swap{Hofstad}{~van~der~}}, and G.~Hooghiemstra.
\newblock Scale-free percolation.
\newblock \emph{Ann. Inst. Henri Poincar\'{e} Probab. Stat.}, 49\penalty0
  (3):\penalty0 817--838, 2013.
\newblock ISSN 0246-0203.
\newblock \doi{10.1214/12-AIHP480}.
\newblock URL \url{https://doi.org/10.1214/12-AIHP480}.

\bibitem[Deprez and W\"{u}thrich(2019)]{CSFP_D_W}
P.~Deprez and M.~V. W\"{u}thrich.
\newblock Scale-free percolation in continuum space.
\newblock \emph{Commun. Math. Stat.}, 7\penalty0 (3):\penalty0 269--308, 2019.
\newblock ISSN 2194-6701.
\newblock \doi{10.1007/s40304-018-0142-0}.
\newblock URL \url{https://doi.org/10.1007/s40304-018-0142-0}.

\bibitem[{\swap{Esker}{~van~den~}} et~al.(2005){\swap{Esker}{~van~den~}},
  {\swap{Hofstad}{~van~der~}}, Hooghiemstra, and
  Znamenski]{RvdH_GH_vdE_CM_const_dist}
H.~{\swap{Esker}{~van~den~}}, R.~{\swap{Hofstad}{~van~der~}}, G.~Hooghiemstra,
  and D.~Znamenski.
\newblock Distances in random graphs with infinite mean degrees.
\newblock \emph{Extremes}, 8\penalty0 (3):\penalty0 111--141 (2006), 2005.
\newblock ISSN 1386-1999.
\newblock \doi{10.1007/s10687-006-7963-z}.
\newblock URL \url{https://doi.org/10.1007/s10687-006-7963-z}.

\bibitem[Fountoulakis et~al.(2021)Fountoulakis, van~der Hoorn, Müller, and
  Schepers]{FHMS_clustering_HRGs}
N.~Fountoulakis, P.~van~der Hoorn, T.~Müller, and M.~Schepers.
\newblock {Clustering in a hyperbolic model of complex networks}.
\newblock \emph{Electronic Journal of Probability}, 26\penalty0
  (none):\penalty0 1 -- 132, 2021.
\newblock \doi{10.1214/21-EJP583}.
\newblock URL \url{https://doi.org/10.1214/21-EJP583}.

\bibitem[Gilbert(1961)]{Gilbert_61}
E.~N. Gilbert.
\newblock Random plane networks.
\newblock \emph{J. Soc. Indust. Appl. Math.}, 9:\penalty0 533--543, 1961.
\newblock ISSN 0368-4245.

\bibitem[Gracar et~al.(2019)Gracar, Heydenreich, Mönch, and
  Mörters]{WDRCM_Rec_trans}
P.~Gracar, M.~Heydenreich, C.~Mönch, and P.~Mörters.
\newblock Recurrence versus transience for weight-dependent random connection
  models.
\newblock 2019.
\newblock URL \url{https://arxiv.org/abs/1911.04350v2}.

\bibitem[Gracar et~al.(2021)Gracar, Grauer, and Mörters]{chem_dist}
P.~Gracar, A.~Grauer, and P.~Mörters.
\newblock Chemical distance in geometric random graphs with long edges and
  scale-free degree distribution.
\newblock 2021.
\newblock URL \url{https://arxiv.org/abs/2108.11252}.

\bibitem[Heydenreich et~al.(2017)Heydenreich, Hulshof, and
  Jorritsma]{SC_SFP_JHH}
M.~Heydenreich, T.~Hulshof, and J.~Jorritsma.
\newblock Structures in supercritical scale-free percolation.
\newblock \emph{Ann. Appl. Probab.}, 27\penalty0 (4):\penalty0 2569--2604,
  2017.
\newblock ISSN 1050-5164.
\newblock \doi{10.1214/16-AAP1270}.
\newblock URL \url{https://doi.org/10.1214/16-AAP1270}.

\bibitem[{\swap{Hofstad}{~van~der~}}(2017)]{RGCN_1}
R.~{\swap{Hofstad}{~van~der~}}.
\newblock \emph{Random graphs and complex networks. {V}ol. 1}.
\newblock Cambridge Series in Statistical and Probabilistic Mathematics, [43].
  Cambridge University Press, Cambridge, 2017.
\newblock ISBN 978-1-107-17287-6.
\newblock \doi{10.1017/9781316779422}.
\newblock URL \url{https://doi.org/10.1017/9781316779422}.

\bibitem[{\swap{Hofstad}{~van~der~}}(2021+{\natexlab{a}})]{Giant_local}
R.~{\swap{Hofstad}{~van~der~}}.
\newblock The giant in random graphs is almost local.
\newblock 2021+{\natexlab{a}}.
\newblock URL \url{https://arxiv.org/abs/2103.11733}.

\bibitem[{\swap{Hofstad}{~van~der~}}(2021+{\natexlab{b}})]{RGCN_2}
R.~{\swap{Hofstad}{~van~der~}}.
\newblock {Random Graphs and Complex Networks, Vol. 2 (in preparation)}.
\newblock 2021+{\natexlab{b}}.
\newblock URL \url{http://www.win.tue.nl/\~rhofstad/NotesRGCNII.pdf}.

\bibitem[Kallenberg(2017)]{Kallenberg_RMTA}
O.~Kallenberg.
\newblock \emph{Random measures, theory and applications}, volume~77 of
  \emph{Probability Theory and Stochastic Modelling}.
\newblock Springer, Cham, 2017.
\newblock ISBN 978-3-319-41596-3; 978-3-319-41598-7.
\newblock \doi{10.1007/978-3-319-41598-7}.
\newblock URL \url{https://doi.org/10.1007/978-3-319-41598-7}.

\bibitem[Komj\'{a}thy and Lodewijks(2020)]{JK_BL_explosions_HRGs_19}
J.~Komj\'{a}thy and B.~Lodewijks.
\newblock Explosion in weighted hyperbolic random graphs and geometric
  inhomogeneous random graphs.
\newblock \emph{Stochastic Process. Appl.}, 130\penalty0 (3):\penalty0
  1309--1367, 2020.
\newblock ISSN 0304-4149.
\newblock \doi{10.1016/j.spa.2019.04.014}.
\newblock URL \url{https://doi.org/10.1016/j.spa.2019.04.014}.

\bibitem[Krioukov et~al.(2010)Krioukov, Papadopoulos, Kitsak, Vahdat, and
  Bogu\~{n}\'{a}]{KPKVB_HRG}
D.~Krioukov, F.~Papadopoulos, M.~Kitsak, A.~Vahdat, and M.~Bogu\~{n}\'{a}.
\newblock Hyperbolic geometry of complex networks.
\newblock \emph{Phys. Rev. E (3)}, 82\penalty0 (3):\penalty0 036106, 18, 2010.
\newblock ISSN 1539-3755.
\newblock \doi{10.1103/PhysRevE.82.036106}.
\newblock URL \url{https://doi.org/10.1103/PhysRevE.82.036106}.

\bibitem[Kulik and Soulier(2020)]{Kulik_Soulier_2020}
R.~Kulik and P.~Soulier.
\newblock \emph{Regularly varying random variables}, pages 3--21.
\newblock Springer New York, New York, NY, 2020.
\newblock ISBN 978-1-0716-0737-4.
\newblock \doi{10.1007/978-1-0716-0737-4_1}.
\newblock URL \url{https://doi.org/10.1007/978-1-0716-0737-4_1}.

\bibitem[Last and Penrose(2018)]{Last_Penrose_LPP}
G.~Last and M.~Penrose.
\newblock \emph{Lectures on the {P}oisson process}, volume~7 of \emph{Institute
  of Mathematical Statistics Textbooks}.
\newblock Cambridge University Press, Cambridge, 2018.
\newblock ISBN 978-1-107-45843-7; 978-1-107-08801-6.

\bibitem[Lee and Campbell(2018)]{Nbr_nets_B_W}
B.~Lee and K.~Campbell.
\newblock \emph{Neighbor Networks of Black and White Americans}, pages
  119--146.
\newblock 10 2018.
\newblock ISBN 9780429498718.
\newblock \doi{10.4324/9780429498718-4}.

\bibitem[Penrose(2003)]{Penrose_RGGs}
M.~Penrose.
\newblock \emph{Random geometric graphs}, volume~5 of \emph{Oxford Studies in
  Probability}.
\newblock Oxford University Press, Oxford, 2003.
\newblock ISBN 0-19-850626-0.
\newblock \doi{10.1093/acprof:oso/9780198506263.001.0001}.
\newblock URL \url{https://doi.org/10.1093/acprof:oso/9780198506263.001.0001}.

\bibitem[Potter(1940)]{Potter}
H.~S.~A. Potter.
\newblock The mean values of certain {D}irichlet series, {II}.
\newblock \emph{Proc. London Math. Soc. (2)}, 47:\penalty0 1--19, 1940.
\newblock ISSN 0024-6115.
\newblock \doi{10.1112/plms/s2-47.1.1}.
\newblock URL \url{https://doi.org/10.1112/plms/s2-47.1.1}.

\bibitem[Wellman and Wortley(1990)]{diff_strokes_friendhsip}
B.~Wellman and S.~Wortley.
\newblock Different strokes from different folks: Community ties and social
  support.
\newblock \emph{American Journal of Sociology}, 96\penalty0 (3):\penalty0
  558--588, 1990.
\newblock ISSN 00029602, 15375390.
\newblock URL \url{http://www.jstor.org/stable/2781064}.

\bibitem[Wellman et~al.(1988)Wellman, Carrington, and
  Hall]{wellman1988networks}
B.~Wellman, P.~Carrington, and A.~Hall.
\newblock Networks as personal communities.
\newblock \emph{Social structures: A network approach}, 2:\penalty0 130--184,
  1988.

\bibitem[Wong et~al.(2006)Wong, Pattison, and Robins]{Spatial_social_nets}
L.~H. Wong, P.~Pattison, and G.~Robins.
\newblock A spatial model for social networks.
\newblock \emph{Physica A: Statistical Mechanics and its Applications},
  360\penalty0 (1):\penalty0 99--120, 2006.
\newblock ISSN 0378-4371.
\newblock \doi{https://doi.org/10.1016/j.physa.2005.04.029}.
\newblock URL
  \url{https://www.sciencedirect.com/science/article/pii/S0378437105004334}.

\end{thebibliography}

\end{document}